\documentclass{pspum-l}

\usepackage{amscd,amsmath,amssymb,amsfonts,stmaryrd}
\usepackage[cmtip, all]{xy}
\usepackage{enumerate}
\usepackage{amscd,amsmath,amssymb,amsfonts}

\usepackage{graphicx}

\usepackage{appendix}
\usepackage{mathrsfs}
\usepackage[usenames,dvipsnames]{xcolor}
\usepackage{mathrsfs}
\usepackage{tikz-cd}
\usepackage{framed, hyperref}

\numberwithin{equation}{section}

\newcounter{keepeqno}

\newenvironment{num}
{\setcounter{keepeqno}{\value{equation}}%
	\begin{list}{(\theequation)}{\usecounter{equation}}%
		\setcounter{equation}{\value{keepeqno}}}
	{\end{list}}

\newtheorem{theorem}{Theorem}
\newtheorem{thm}[theorem]{Theorem}
\newtheorem{hyp}[theorem]{Hypothesis}

\newtheorem{lemma}[theorem]{Lemma}

\newtheorem{proposition}[theorem]{Proposition}
\newtheorem{prop}[theorem]{Proposition}
\newtheorem{corollary}[theorem]{Corollary}

\newtheorem{ideal}[theorem]{Ideal Theorem}

\theoremstyle{remark}

\newtheorem{remark}[theorem]{Remark}

\numberwithin{theorem}{section}
\numberwithin{equation}{section}

\DeclareMathOperator{\Ad}{Ad}
\DeclareMathOperator{\Bun}{Bun}

\DeclareMathOperator{\Out}{Out}
\DeclareMathOperator{\GL}{GL}
\DeclareMathOperator{\SL}{SL}
\DeclareMathOperator{\PGL}{PGL}

\DeclareMathOperator{\Supp}{Supp}

\DeclareMathOperator{\Frob}{Frob}

\DeclareMathOperator{\Hom}{Hom}

\DeclareMathOperator{\Gal}{Gal}

\DeclareMathOperator{\diag}{diag}

\DeclareMathOperator{\tr}{tr}

\DeclareMathOperator{\Lie}{Lie}

\DeclareMathOperator{\Spec}{Spec}

\DeclareMathOperator{\vol}{vol}

\newcommand{\isoarrow}{{~\overset\sim\longrightarrow~}}

\newcommand{\ad}{\mathbf{A}}

\newcommand{\cA}{{\mathcal A}}
\newcommand{\CA}{{\mathcal A}}

\newcommand{\cC}{{\mathcal C}}

\newcommand{\cE}{{\mathcal E}}
\newcommand{\cF}{{\mathcal F}}

\newcommand{\cH}{{\mathcal H}}
\newcommand{\cI}{{\mathcal I}}

\newcommand{\cL}{{\mathcal L}}

\newcommand{\CL}{{\mathcal L}}

\newcommand{\cO}{{\mathcal O}}
\newcommand{\cP}{{\mathcal P}}

\newcommand{\ra}{{\rightarrow}}

\newcommand{\bbA}{{\mathbb A}}

\newcommand{\bbC}{{\mathbb C}}

\newcommand{\bbF}{{\mathbb F}}

\newcommand{\bbG}{{\mathbb G}}
\newcommand{\bbH}{{\mathbb H}}

\newcommand{\bbK}{{\mathbb K}}
\newcommand{\bbL}{{\mathbb L}}

\newcommand{\bbN}{{\mathbb N}}

\newcommand{\bbQ}{{\mathbb Q}}

\newcommand{\bbT}{{\mathbb T}}

\newcommand{\bbZ}{{\mathbb Z}}

\newcommand{\hG}{{\widehat{G}}}
\newcommand{\lG}{{{}^LG}}

\DeclareMathOperator{\disc}{disc}

\newcommand{\nc}{\newcommand}
\nc{\renc}{\renewcommand}
\nc{\ssec}{\subsection}
\nc{\sssec}{\subsubsection}
\nc{\on}{\operatorname}

\nc\ol{\overline}
\nc\wt{\widetilde}
\nc\tboxtimes{\wt{\boxtimes}}
\nc{\alp}{\alpha}

\nc{\BunBb}{\overline{\Bun}_B}

\newcommand{\Z}{\mathbb{Z}}

\DeclareMathOperator{\Tr}{Tr}

\DeclareMathOperator{\Ker}{Ker}

\title{Inductive construction of supercuspidal $L$-packets}
\author{Rapha\"el Beuzart-Plessis, Michael Harris, Jack A. Thorne}
\thanks{R. B-P. was partially supported by the European Union (ERC CoG-101044930).  M.H. was partially supported by NSF Grant DMS-2302208.  J.T.'s work was funded by the European Union (ERC CoG-101169866).}

\begin{document}
\begin{abstract} Genestier--Lafforgue and Fargues--Scholze have constructed a semi-simple local Langlands paramterization for reductive groups over equicharacteristic local fields.  Assuming a version of the stable twisted trace formula for function fields, we prove the surjectivity of this parameterization for unramified quasi-split semisimple groups in characteristic not dividing the order of the Weyl group. 	\end{abstract}
	\maketitle
	
	\setcounter{tocdepth}{1}
	\tableofcontents

    \section{Introduction}

    \textbf{Context.} The subject of this paper is the conjectural local Langlands correspondence over a non-archimedean local field of equal characteristic. Let $F$ be such a field, let $G$ be a reductive group over $F$, and let ${}^L G = \widehat{G} \rtimes W_F$ be its associated $L$-group. In this case, the local Langlands correspondence $\mathcal{L}_F$ should associate to irreducible admissible $\mathbb{C}[G(F)]$-module $\pi$ a semisimple $L$-homomorphism
    \[ \mathcal{L}_F(\pi) : W_F \times \SL_2(\mathbb{C}) \to {}^L G(\mathbb{C}). \]
    In the case $G = \GL_n$, the correspondence is known to exist \cite{LRS93}. Moreover, it is a bijection (between the set of isomorphism classes of representations $\pi$ and the set of $\widehat{G}$-conjugacy classes of semisimple $L$-homomorphisms). 
    
    For a general reductive group $G$, something close to the map $\mathcal{L}_F$ has been constructed (see the next paragraph), but the properties of the map $\mathcal{L}_F$ are expected to be more subtle: in particular, it should not be bijective, but rather its fibres should partition the isomorphism classes of irreducible admissible representations $\pi$ into $L$-packets (see \cite{Ka16} for a survey of their expected properties). The main question we aim to address in this paper is the surjectivity of $\mathcal{L}_F$ (equivalently: the non-emptiness of certain $L$-packets).

    Following work of Genestier--Laffogure, Li-Huerta, and Scholze, it is known that there is a well-defined semisimplified local Langlands correspondence $\mathcal{L}_F^{ss}$, which associates to any irreducible, admissible $\mathbb{C}[G(F)]$-module $\pi$ a semisimple homomorphism
    \[ \mathcal{L}^{ss}_F(\pi) : W_F \to {}^L G(\mathbb{C}), \]
    in a manner compatible with Lafforgue's global Langlands correspondence, parabolic induction, and twist by unramified characters. (We discuss these results in more detail in \S \ref{sec_review_of_Galois} below.) The expected relation between $\mathcal{L}_F$ and $\mathcal{L}^{ss}_F$ is as follows: for any $w \in W_F$, we should have 
    \[ \mathcal{L}_F^{ss}(w) = \mathcal{L}_F(w, \left( \begin{array}{cc} \| w \|^{1/2} & 0 \\ 0 & \| w \|^{-1/2} \end{array}\right) ). \]
    In this paper, we focus on $L$-homomorphisms $\phi : W_F \to {}^L  G(\mathbb{C})$ that are irreducible, in the sense that they have finite centralizer modulo $Z(\widehat{G})^{W_F}$ (in fact, in the body of the paper, we restrict to the case that $G$ is semisimple, so we need not work modulo the Galois-fixed points in the center of $\widehat{G}$). If $\mathcal{L}_F^{ss}(\pi) = \phi$, then one expects that $\mathcal{L}_F(\pi) = \mathcal{L}^{ss}_F(\pi)$ (i.e. that $\mathcal{L}_F(\pi)$ is trivial on the $\SL_2(\mathbb{C})$-factor) and that the $L$-packet $\mathcal{L}_F^{-1}(\phi)$ consists of supercuspidal representations. 
    
    \textbf{Results of this paper.} Our main theorem shows that, under some conditions, the $L$-packet of an irreducible parameter is non-empty and contains a generic supercuspidal representation: 
    \begin{theorem}[Theorem \ref{existenceofpi}]\label{introthm_descent}
        Let $G$ be an unramified semisimple group over the non-archimedean, equicharacteristic local field $F$. Let $\tau : W_F \to {}^L G(\mathbb{C})$ be an irreducible $L$-homomorphism. Suppose that the following conditions are satisfied:
        \begin{enumerate}
            \item The characteristic of $F$ does not divide the order of the (absolute) Weyl group of $G$.
            \item The hypotheses of Theorem \ref{hyp rho traces} are satisfied.
        \end{enumerate}
        Then there exists an irreducible, admissible, generic, supercuspidal $\mathbb{C}[G(F)]$-module $\pi$ such that $\mathcal{L}_F^{ss}(\pi) = \tau$. 
    \end{theorem}
    The hypotheses in the statement of Theorem \ref{hyp rho traces} amount to assuming a version of the stable twisted trace formula for function fields analogous to known formulas over number fields. We require this only for test functions supported on strongly regular semisimple elements. We believe that this version of the stable twisted trace formula is within reach, as the main difficulty in generalising the arguments of the number field case to the case of a global field of positive characteristic arises from the contribution of the non-regular semisimple elements to the geometric side of the trace formula.

    \textbf{Methods.} As one might imagine, based on the presence of assumptions about the trace formula, the argument that we use to prove Theorem \ref{introthm_descent} is global, and follows the broad lines of the conjectural sketch given in \cite{Boc19}: first globalise $\tau$ to a global parameter $\tau_0$, that becomes potentially automorphic over some extension of the base field; then use cyclic descent (together with the solvability of local Galois groups) to find an automorphic representation whose parameter localises to $\tau$, therefore (by local-global compatibility) with a local component $\pi$ satisfying $\mathcal{L}^{ss}_F(\pi) = \tau$. The contribution of this paper is to show that this strategy does work, assuming a good enough understanding of the stable twisted trace formula.

    The main new ingredient that makes this possible is the knowledge that the automorphic representation witnessing the potential automorphy of $\tau_0$ occurs with automorphic multiplicity 1 (an essential hypothesis in applications of the twisted trace formula). This multiplicity-1 property has been proved by Gaitsgory--Raskin \cite{GR25} as a consequence of the categorical Langlands correspondence for curves. We need to show that this property is preserved by cyclic descent at each step of the argument; this requires some local character identities.  When $p$ is prime to the order of the Weyl group, the image under the parameter of wild inertia (and perhaps a little more) is contained in a maximal torus. The induction steps at this stage are based on a known identity for the twisted trace of a principal series representation.  In the remaining stages the cyclic subquotient is of order prime to $p$.  To continue the induction  we use an argument going back to Rodier, Moeglin, and Waldspurger to pick out a {\it type} -- a $1$-dimensional subspace of the local representation at the primes of interest -- on which an appropriately chosen compact open subgroup acts by a character.  
	
	The local properties of cyclic  descent, as developed in Henniart's proof of the numerical local Langlands correspondence for $\GL_n$, are indispensable in every known proof of the local Langlands correspondence for $\GL_n$.   In Arthur's study of the automorphic spectrum of classical groups, the construction of supercuspidal $L$-packets is carried out in a single descent relating the twisted trace of $\GL_n$ with respect to an outer automorphism and distributions on related classical groups.  The key steps in this construction are contained in \cite[\S 6]{Art13}.  In a future work we hope to apply Arthur's ideas to the inductive construction developed here in order to construct supercuspidal $L$-packets generally, and to relate them to the stable packets constructed in \cite{Ka19b}.  

    \textbf{Structure of this paper.} In \S \ref{sec_review_of_Galois}, we recall in more detail the existing state of knowledge concerning the local correspondence $\mathcal{L}_F^{ss}$, and give some basic results about local $L$-parameters. In \S \ref{tracehypo}, we make explicit and precise the properties of the stable twisted trace formula over global fields of positive characteristic that we need to assume in order to prove our main theorems. 

    The next few sections are preparatory. In \S \ref{types}, we construct the types required for the tame part of the descent argument, and show the required character identities. In \S \ref{sec_globalization}, we construct globalizations of the local parameters we consider with the required properties, including their potential automorphy over a suitable extension of the base field. In \S \ref{sec_mult_one}, we summarise what we need from the work of Gaitsgory--Raskin concerning multiplicity one.

    We then begin the main arguments. First, in \S \ref{secKV}, we give an idealised version of the descent argument, based on the results of \cite{Ko, Va17}, which are however known to hold only for local fields of characteristic 0. Then, in \S \ref{sec_descent_arguments}, we carry out the argument to prove Theorem \ref{introthm_descent} as outlined above. Finally, in \S \ref{sec_p-adic_apps} we combine our results with the work of \cite{LH24} to draw consequences for the surjectivity of the local Langlands correspondence over $p$-adic local fields constructed in \cite{FS}.   
    
   We thank David Hansen for providing a short Appendix that explains how our main existence result extends to inner forms of quasi-split groups.  We also thank Jean-Loup Waldspurger for the argument used in the proof of Proposition \ref{ptw}, and the referee for a careful reading and for suggesting a much simpler proof of Lemma \ref{lem_connectedness_of_centralizer}.
	
	\subsection{Notation} If $K$ is a field, then we will write $\Gamma_K = \Gal(K^s / K)$ for its Galois group with respect to some fixed choice of separable closure. If $K$ is a global field and $v$ is a place of $K$, then we will suppose fixed a compatible embedding $K^s \to K_v^s$ extending the inclusion $K \to K_v$, which therefore determines a map $\Gamma_{K_v} \to \Gamma_K$, well-defined up to conjugacy in $\Gamma_K$. If $K$ is a non-archimedean local field, then we write $W_K \leq \Gamma_K$ for the Weil group, $I_K \leq W_K$ for the inertia group, and $P_K \leq I_K$ for the wild inertia group. We write $\Frob_v \in W_K$ for a choice of lift of geometric Frobenius.

	\section{Review of Galois parameters}\label{sec_review_of_Galois}


    We begin with the global case. Let $\mathbb{K}$ be a global field of positive characteristic, and let $G$ be a semisimple group over $\mathbb{K}$. We write $\mathcal{A}_0(G)$ for the space of cuspidal automorphic forms on $G$, which is a semisimple admissible $\mathbb{C}[G(\mathbb{A}_{\mathbb{K}})]$-module. We write ${}^L G = \widehat{G} \rtimes \Gamma_{\mathbb{K}}$ for the (Galois form of) the $L$-group of $G$. 

    By definition, a cuspidal automorphic representation of $G$ is an irreducible $\mathbb{C}[G(\mathbb{A}_{\mathbb{K}})]$-module that occurs as a subrepresentation of $\mathcal{A}_0(G)$. Any such representation may be defined over $\overline{\mathbb{Q}} \subset \mathbb{C}$. We will shortly define $\ell$-adic Galois representations associated to cuspidal automorphic representations. Their definition depends on a choice of square root of $q$ (the cardinality of the field of constants of $\mathbb{K}$) in $\overline{\mathbb{Q}}_\ell$; we will address this by defining our representations over $\overline{\mathbb{Q}}_\lambda$, where $\lambda$ is a place of $\overline{\mathbb{Q}}$ of characteristic $\ell$ prime to $q$, and taking the square root of $q$ that corresponds to the positive square root in $\mathbb{C}$. 

    The following theorem is \cite[Th\'eor\`eme 0.1]{Lafa}.
    \begin{theorem}
        Let $\Pi$ be a cuspidal automorphic representation of $G$, and let $\lambda$ be a prime-to-$q$ place of $\overline{\mathbb{Q}}$. Then we can find a continuous semisimple representation 
        \[ r_\lambda(\Pi) : \Gamma_{\mathbb{K}} \to {}^L G (\overline{\mathbb{Q}}_\lambda) \]
        with the following property: for any place $v$ of $\mathbb{K}$ such that $\Pi_v$ is unramified, $r_\lambda(\Pi)|_{W_{\mathbb{K}_v}}$ is unramified, and the semisimple part of $r_\lambda(\Pi)(\Frob_v)$ is conjugate under  $\widehat{G}(\overline{\mathbb{Q}}_\lambda)$ to the Satake parameter of $\Pi_v$.
    \end{theorem}
    The representation $r_\lambda(\Pi)$ is not uniquely determined without additional conditions. The following proposition, the first part of which extends \cite[Proposition 6.4]{Boc19}, gives one such condition: 
	\begin{proposition}\label{prop ext BHKT}
		Suppose that $G$ splits over a cyclic extension $\mathbb{K}' / \mathbb{K}$. Then: 
		\begin{enumerate}
			\item Let $\rho,\rho': \Gamma_{\bbK}\to {}^L G(\overline{\bbQ}_{\ell})$ be $L$-homomorphisms. Assume that $\rho(\Gamma_{\mathbb{K}'})$ is Zariski dense in $\widehat{G}$ and that there exists a finite set $S$ of places such that for every $v\notin S$, $\rho$, $\rho'$ are unramified at $v$ and the semi-simple parts of $\rho(\Frob_v)$, $\rho'(\Frob_v)$ are $\widehat{G}(\overline{\bbQ}_{\ell})$-conjugate. Then, $\rho$ and $\rho'$ are $\widehat{G}(\overline{\bbQ}_\ell)$-conjugate, that is there exists an element $g\in \widehat{G}(\overline{\bbQ}_\ell)$ such that $\rho'(\sigma)=g\rho(\sigma)g^{-1}$ for every $\sigma\in \Gamma_{\bbK}$.
			\item We have $\ker^1(\bbK,Z(\widehat{G}))=0$.
		\end{enumerate}
	\end{proposition}
	Here, $\ker^1(\bbK,Z(\widehat{G}))$ denotes the kernel of the localization morphism between Galois cohomology groups
	$$\displaystyle H^1(\bbK,Z(\widehat{G})(\overline{\mathbb{Q}}_\lambda))\to \prod_v H^1(\bbK_v, Z(\widehat{G})(\overline{\mathbb{Q}}_\lambda)),$$
	the product being taken over all places of $\bbK$.
	\begin{proof}
		(i) By \cite[Proposition 6.4]{Boc19}, applied to $\rho|_{\Gamma_{\mathbb{K}'}}$ there exists $g\in \widehat{G}$ such that $g\rho'(\sigma)g^{-1}=\rho(\sigma)$ for every $\sigma\in \Gamma_{\bbK'}$. Up to conjugating $\rho'$ by $g$, we may therefore assume that $\rho$ and $\rho'$ have the same restriction to $\Gamma_{\bbK'}$.
			
			Let $\tau$ be a generator of $\Gal(\bbK'/\bbK)$ and set $a:=\rho'(\tilde{\tau})\rho(\tilde{\tau})^{-1}$, where we have chosen an arbitrary lift $\tilde{\tau}$ of $\tau$ in $\Gamma_{\bbK}$. It is an element of $\widehat{G}(\overline{\bbQ}_\ell)$ that commutes with $\rho(\Gamma_{\bbK'})$. Thus, by Zariski density of the latter, we have $a\in Z(\widehat{G})$ and it suffices to show that there exists $z\in Z(\widehat{G})$ with $a=z\tau(z)^{-1}$ as then $g':=gz$ would conjugate $\rho$ to $\rho'$.
			
			Let $S$ be a finite set of places of $\bbK$ such that for $v\notin S$, $\rho$, $\rho'$, $\bbK'$ are unramified at $v$ and the semisimple parts of $\rho(\Frob_v)$, $\rho'(\Frob_v)$ are conjugate by $\widehat{G}(\overline{\bbQ}_\ell)$; where $\Frob_v\in \Gamma_{\bbK,S}$ is a chosen Frobenius element. Let $T$ be the set of places $v\notin S$ such that the projection of $\Frob_v$ to $\Gal(\bbK'/\bbK)$ is $\tau$. In what follows, for $v\in T$, we will identify $\rho(\Frob_v)$, $\rho'(\Frob_v)$ with their image in $\widehat{G}(\overline{\bbQ}_\ell)\rtimes \Gal(\bbK'/\bbK)$; in particular these are elements of the coset $\widehat{G}(\overline{\bbQ}_\ell)\tau$. 
			
			By the Cebotarev density theorem, the union of the $\Gamma_{\bbK'}$-conjugacy classes of the elements $\Frob_v\in \Gamma_{\bbK,S}$, for $v\in T$, is dense in the coset $\Gamma_{\bbK',S}\tau$. Therefore, $\rho(\Gamma_{\bbK',S})$ being Zariski dense in $\widehat{G}(\overline{\bbQ}_\ell)$, we have that for every semi-simple element $\gamma\in \widehat{G}(\overline{\bbQ}_\ell)\tau$, $a\gamma$ and $\gamma$ are $\widehat{G}(\overline{\bbQ}_\ell)$-conjugate. In particular, there exists $h\in \widehat{G}(\overline{\bbQ}_\ell)$ such that $a\tau=h \tau h^{-1}$. Let $\mathcal{P}=(\widehat{T},\widehat{B}, \{X_\alpha \})$ be the canonical pinning of $\widehat{G}$, the one that is fixed by the Galois action. Then, the previous equality implies that $\tau$ also preserves the pinning $\Ad(h^{-1})\mathcal{P}$. By \cite[Corollary 1.7]{K84}, it follows that we can find $h_0\in \widehat{G}^{\tau}$ such that $\Ad(h^{-1})\mathcal{P}=\Ad(h_0)\mathcal{P}$. Set $z=hh_0$. Then, $a \tau=z \tau z^{-1}$, that is $a=z\tau(z)^{-1}$, and $z$ preserves the pînning $\mathcal{P}$, hence $z\in Z(\widehat{G})$.
			
(ii) Let $a: \Gamma_{\bbK}\to Z(\widehat{G})$ be a $1$-cocycle that is everywhere locally cohomologically trivial. Then, $a$ restricts to a morphism $\Gamma_{\bbK'}\to Z(\widehat{G})$ that is locally trivial, hence trivial by Cebotarev. Thus, $a$ can be replaced by a cohomologous cocycle that factors through $\Gal(\bbK'/\bbK)$. By Cebotarev again, the extension $\bbK'/\bbK$ being cyclic, we can find a place $v$ that is inert in $\bbK'$, hence $\Gal(\bbK'_v/\bbK_v)=\Gal(\bbK'/\bbK)$. By assumption, the localization $a_v$ of $a$ is cohomologically trivial and therefore so does $a$.
		
	\end{proof}
    In particular, we may say that $r_\lambda(\Pi)$ is unambiguously associated to $\Pi$ whenever its image is Zariski dense in the finite $L$-group.

    We now discuss the local correspondence. Let $K$ be a non-archimedean local field of positive characteristic, and let $G$ be a semisimple group over $K$. The following theorem combines results of \cite{GLa}, \cite{FS}, \cite{LH23}, and \cite{Sc25}.
	\begin{thm}\label{thm_loc_correspondence} Associated to any irreducible admissible $\mathbb{C}[G(K)]$-module $\pi$ that is definable over $\overline{\bbQ}$ is a semisimple $L$-homomorphism 
    \[ \mathcal{L}^{ss}_K(\pi) : W_K \to {}^L G(\overline{\mathbb{Q}}) \]
    with the following properties:
    \begin{enumerate}
        \item Local-global compatibility: let $\mathbb{K}$ be a global field with a place $v$ such that $\mathbb{K}_v \cong K$, and let $\mathbb{G}$ be a semisimple group over $\mathbb{K}$ such that $\mathbb{G}_K \cong G$. Let $\Pi$ be a cuspidal automorphic representation of $\mathbb{G}$, and let $\lambda$ be a prime-to-$q$ place of $\overline{\mathbb{Q}}$. Then the homomorphisms $r_\lambda(\Pi)|_{W_K}^{ss}$ and $\mathcal{L}^{ss}_K(\Pi_v)$ are $\widehat{G}(\overline{\mathbb{Q}}_\lambda)$-conjugate. 
        \item Compatibility with parabolic induction: Let $P \subset G$
		be a $K$-rational parabolic subgroup with Levi factor $M$ and let $\sigma$ be an  irreducible admissible representation of $M(K)$, that is definable over $\overline{\bbQ}$.  If $\pi$ is an irreducible subquotient of the normalized induction $I_{P(K)}^{G(K)} \sigma$ and $i_M:  {}^LM \ra {}^LG$ is the inclusion of the $L$-group of $M$
		as a Levi factor  of a parabolic subgroup in the $L$-group of $G$, then 
		$\CL_K^{ss}(\pi)$ and $i_M\circ\CL_K^{ss}(\sigma)$ are $\widehat{G}$-conjugate.
    \end{enumerate} 
	\end{thm}
    We recall that a semisimple $L$-homomorphism $\tau : W_K \to {}^L G(\mathbb{C})$ is said to be irreducible if its image is  not contained in any proper parabolic subgroup of ${}^L G$; equivalently (since $G$ is semisimple), if its centralizer is finite.  
	\begin{lemma}[Kaletha]\label{irrparam} Let 
    \[ \tau : W_K \to {}^L G(\mathbb{C}) \]
    be an irreducible $L$-homomorphism.   Suppose that the characteristic $p$ of $K$ is prime to the order of the absolute Weyl group of $G$, and that $G$ is unramified, so that in particular we can take ${}^L G=\widehat{G}\rtimes \Gal(K' / K)$ where $K' / K$ is a finite unramified extension of degree prime to $p$\footnote{The assumption implies that $p$ doesn't divide the order of the group of outer automorphisms of $\widehat{G}$ either, provided $G$ is almost simple.}. Then, after replacing $\tau$ by a $\widehat{G}$-conjugate, the following statements hold:
		\begin{itemize}
			\item[(a)]  The image $\tau(P_F)$ is contained in the maximal torus $\widehat{T} \subset \widehat{G}$.  
			\item[(b)]  The image $\tau(W_F)$ is contained in the normalizer $N_{{}^L G}(\widehat{T}) \subset {}^L G$ of $\widehat{T}$ in ${}^L G$.  In particular, the $p$-Sylow subgroup $\tau(W_F)_p$ of $\tau(W_F)$ is contained in $\widehat{T}$, hence is abelian and is normal in $\tau(W_F)$.
			
		\end{itemize}
	\end{lemma}
	\begin{proof}
		Point (a) is \cite[Lemma 4.1.3 1]{Ka19b}. Since $p$ is prime to the order of the Weyl group, in order to prove point (b)  it suffices to prove that the image of $\tau$ is contained in $N_{{}^L G}(\hat{T})$. But this follows from \cite[Lemma 4.1.3 2]{Ka19b} and \cite[Lemma 5.2.2 (2)]{Ka19a}.  See also the proof of \cite[Proposition 4.1.7]{Ka19b}.
	\end{proof}
	
	\begin{lemma}\label{irrprin} Let $\tau : W_K \to {}^L G(\mathbb{C})$ be an irreducible $L$-homomorphism. Suppose that $p$ is prime to the order of the absolute Weyl group of $G$, and that $G$ is unramified. Let $K_\tau/K$ be the extension of $K$ cut out by $\tau$, and let $K_{prin} \subset K_\tau$ denote the fixed field of $\tau(W_K)_p$, in the notation of Lemma \ref{irrparam}.  For any extension $K'/K_{prin}$ contained in $K_\tau$ let $\chi':  T(K') \ra \mathbb{C}^\times$ be the character attached by the Langlands correspondence for tori to the homomorphism $\Gal(K_\tau/K') \ra \hat{T}$ given by restricting $\tau$.   Then the (normalized) principal series representation $I(\chi') = I_{B(K')}^{G(K')} \chi'$ is irreducible.
	\end{lemma}
	
	\begin{proof}   Let $C(K') \subset \hat{T}$ denote the image $\tau(W_{K'})$.  By the theory of the $R$-group\footnote{More precisely, we use here the Galois theoretic description of the $R$-group of $I(\chi')$, see \cite[Section 5]{Sol25}.},	the irreducible components of  the normalized induction $ I(\chi')$ of $\chi'$ are in bijection with the set of connected components of the centralizer $Z_{\hat{G}}(C(K'))$ of $C(K')$ in $\hat{G}$.  Thus, we need to see that this centralizer is connected.  By hypothesis $C(K')$ is an abelian $p$-group, so it suffices to show that, for any finite abelian group $C \subset \hat{T}$ of order prime to the order of the Weyl group, its centralizer $Z(C) \subset \hat{G}$ is a Levi subgroup.   Thus the result follows from Lemma \ref{lem_connectedness_of_centralizer} below.
	\end{proof}

	\begin{lemma}\label{lem_connectedness_of_centralizer}  Let $\widehat{G}$ be a complex connected reductive group with Weyl group $W$ and let $C\subset \widehat{G}$ be a finite abelian subgroup.
		If $\lvert C\rvert$ is prime to $\lvert W\rvert$ then the centralizer $Z_{\widehat{G}}(C)$ of $C$ in $\widehat{G}$ is a Levi subgroup.  In particular it is connected.	
	\end{lemma}
	
	\begin{proof}
        (We thank the referee for pointing out a much shorter proof than the one we gave in the first version of this paper. The following proof uses similar ideas.) Inducting on $\dim \widehat{G}$, it suffices to see that the centralizer of any non-central element $c\in C$ is contained in a proper Levi subgroup. For this, we may assume that $\widehat{G}$ is adjoint and, up to conjugating $c$, that $c\in\widehat{T}$. Let $\Phi$ be the set of all roots of $\widehat{T}$ in $\widehat{G}$ and set
		$$\Phi_c=\{\alpha\in \Phi\mid \alpha(c)=1\}.$$
		Then, $c$ is contained in a proper Levi subgroup if and only if the $\mathbb{Q}$-span of $\Phi_c$ is different from $X^*(\widehat{T})\otimes \mathbb{Q}=\mathbb{Q}\Phi$. (This follows easily from the description of Levi subgroups as centralizers of subtori.)
		Therefore, we need to show that if $c\in \widehat{T}$ has finite order prime to $\lvert W\rvert$ then $\mathbb{Q}\Phi_c\neq \mathbb{Q}\Phi$. Since $\widehat{G}$ is adjoint, the map $t\mapsto (\chi\mapsto \chi(t))$ induces an isomorphism between the torsion subgroup $\widehat{T}_{tors}$ of $\widehat{T}$ and $\Hom(\mathbb{Z}\Phi, \mu_\infty(\mathbb{C}))$, where $\mu_\infty(\mathbb{C})\subset \mathbb{C}^\times$ stands for the subgroup of roots of unity. By this isomorphism, $c$ is sent to an element in $\Hom(\mathbb{Z}\Phi/\mathbb{Z}\Phi_c, \mu_\infty(\mathbb{C}))$. Thus, if $\mathbb{Q}\Phi_c=\mathbb{Q}\Phi$, the order of $c$ divides the finite index $[\mathbb{Z}\Phi:\mathbb{Z}\Phi_c]$. By definition \cite[E I \S 4]{Spr70}, the primes dividing the latter index are all bad primes for $\widehat{G}$. However, by \cite[E I \S 4.3 and \S 4.10]{Spr70}, the set of bad primes for $\widehat{G}$ is included in the set of prime divisors of $\lvert W\rvert$. By assumption the order of $c$ is prime to $\lvert W\rvert$ and therefore $\mathbb{Q}\Phi=\mathbb{Q}\Phi_c$ implies $c=1$ and we are done.
	\end{proof}

	\section{Hypotheses about the trace formula}\label{tracehypo} 
	
	Let $\bbK$ be a global field of characteristic $p$ and $G$ be a semi-simple group over $\bbK$. Let $\lambda$ be a place of $\overline{\mathbb{Q}}$ of characteristic prime to $p$. We fix an $L$-homomorphism
	$$\displaystyle \rho: \Gamma_{\mathbb{K}}\to {}^LG(\overline{\mathbb{Q}}_{\lambda})$$
	with Zariski dense image. 
	We will suppose throughout this section that $G$ splits over a cyclic extension $\mathbb{K}' / \mathbb{K}$,  in order to be able to apply Proposition \ref{prop ext BHKT}. 
	
	Let $\bbL/\bbK$ be a cyclic extension of order $d$. (We emphasize that $\bbL$ has a priori nothing to do with the splitting cyclic extension $\mathbb{K}'$.) We choose a generator $\theta$ of $\Gal(\bbL/\bbK)$. Then, $\theta$ naturally induces an automorphism (denoted also by $\theta$) of the Weil restriction of scalars $R_{\bbL/\bbK} G_{\bbL}$. The goal of this section is to make explicit our hypotheses concerning the stabilization of the trace formulas for $G_{\mathbb{K}}$ and the base-change ``twisted space'' $G_{\mathbb{L}}\theta$. The main output will be an identity relating the trace on the ``$\rho$-subspace'' $\CA_{0,\rho}(G_{\mathbb{K}})$ of the space of cuspidal automorphic forms for $G_{\mathbb{K}}$ with the $\theta$-twisted trace  on the ``$\rho_{\mathbb{L}}$-subspace'' $\CA_{0,\rho_{\mathbb{L}}}(G_{\mathbb{L}})$ (where $\rho_{\mathbb{L}}:=\rho\mid_{\Gamma_{\mathbb{L}}}$). It turns out that this identity doesn't involve stable trace formulas whatsoever but only requires the base-change transfer $f_{\mathbb{L}}\in C_c^\infty(G(\mathbb{A}_{\mathbb{L}}))\to f_{\mathbb{K}}\in C_c^\infty(G(\mathbb{A}_{\mathbb{K}}))$ for its formulation. Therefore, we will begin with the statement of that identity in the first subsection. In the remaining subsections, we discuss stabilization of the trace formulas properly speaking and how if available, this  leads to the aforementioned identity.
	
	Note that Labesse and Lemaire have recently developed a general twisted trace formula over function fields in \cite{LL21}.  Lemaire continues to make progress on its stabilization, as stated in Hypothesis \ref{hyp2} below, but for the moment this is still conjectural.
	
	Some notation and convention for this section are as follows:
	\begin{itemize}
		\item For $G$ a connected reductive group over $\bbK$, we denote by $\CA = \CA(G)$ the space of automorphic forms on $G(\bbA_{\bbK})$ with values in $\mathbb{C}$ and by $\CA_0 = \CA_{0}(G)$ the subspace of cuspidal automorphic forms.  When we wish to emphasise the base field, we will write $\CA_0(G_{\bbK})$ for this space. We also write  $\CA_{\disc}(G)$ for the space of forms $\phi\in \CA(G)$ such that
		$$\displaystyle \int_{G(K)\backslash G(\bbA_{\bbK})} \lvert \iota\phi(g)\rvert^2 dg<\infty.$$
		(We recall that $G$ is assumed to be semisimple throughout.)
		
		\item We  fix an isomorphism $\iota: \overline{\mathbb{Q}}_\lambda \simeq \mathbb{C}$ extending the tautological embedding of $\overline{\mathbb{Q}}$ in both fields. 
		
		
		\item For every place $v$ of $\bbK$, we let $C_c^\infty(G(\bbA_{\bbK}))$ (resp. $C_c^\infty(G(\bbK_v))$) be the space of locally constant and compactly supported $\overline{\bbQ}_{\ell}$-valued functions on $G(\bbA_{\bbK})$ (resp. on $G(\bbK_v)$). 
		
		\item We let $G_{sreg}\subset G$ be the open subset of strongly regular semisimple elements\footnote{We recall that an element $\gamma\in G$ is said to be strongly regular if its centralizer $T_\gamma$ is a (necessarily maximal) torus.}. In a similar way, we let $(R_{\bbL/\bbK} G_{\bbL})_{\theta-sreg}$ be the open subset of $\theta$-strongly regular elements in $R_{\bbL/\bbK} G_{\bbL}$, that is the subset of elements $\delta\in R_{\bbL/\bbK} G_{\bbL}$ whose {\em twisted centralizer}
		\begin{equation}\label{tw centralizer}
		\displaystyle T_{\delta,\theta}=\{g\in R_{\bbL/\bbK} G_{\bbL}\mid g\delta \theta(g)^{-1}=\delta \}
		\end{equation}
		is a torus. For any place $v$ of $\bbK$ we will write 
        \[ G(\bbL_v)_{\theta-sreg} = (R_{\bbL/\bbK} G_{\bbL})_{\theta-sreg}(\bbK_v) \]
        (it is an open-dense subset of $G(\bbL_v)$).
	\end{itemize}

	\subsection{Comparison of $\rho$-traces}\label{subsection rhotrace}
	
	First, we recall the definition of the base-change transfer. Let $v$ be a place of $\mathbb{K}$. We recall that two strongly regular elements $\gamma,\gamma'\in G(\mathbb{K}_v)$ are said to be {\em stably conjugate} if they are conjugate inside $G(\mathbb{K}_v^{sep})$. Let $f'_v\in C_c^\infty(G(\mathbb{K}_v))$ be a test function. Then, for every strongly regular element $\gamma\in G_{sreg}(\mathbb{K}_v)$, we can form the {\em stable orbital integral}
	$$\displaystyle SO_\gamma(f'_v)=\sum_{\gamma'} O_{\gamma'}(f'_v)$$
	where the sum runs over the (finite set of) $G(\mathbb{K}_v)$-conjugacy classes of strongly regular elements $\gamma'$ that are stably conjugate to $\gamma$ and
	$$\displaystyle O_{\gamma'}(f'_v)=\int_{T_{\gamma'}(\mathbb{K}_v)\backslash G(\mathbb{K}_v)} f'_v((g')^{-1}\gamma' g') \frac{dg'}{dt_{\gamma'}}$$
	denotes the orbital integral of $f'_v$ over the conjugacy class of $\gamma'$. Here $T_{\gamma'}$ denotes the centralizer of $\gamma'$ and $dg'$ (resp. $dt_{\gamma'}$) stands for a Haar measure on $G(\mathbb{K}_v)$ (resp. on the torus $T_{\gamma'}(\mathbb{K}_v)$). The choice of $dg'$ does not really matter; the Haar measures $dt_{\gamma'}$, on the other hand, should be chosen in a compatible way in the following sense: for every $\gamma'\in G(\mathbb{K}_v)$ stably conjugate to $\gamma$, there is a canonical isomorphism $T_\gamma\simeq T_{\gamma'}$ defined over $\mathbb{K}_v$, and we ask that it sends the Haar measure $d t_\gamma$ on $T_\gamma(\mathbb{K}_v)$ to the Haar measure $d t_{\gamma'}$ on $T_{\gamma'}(\mathbb{K}_v)$.

	Let $\Ad_\theta$ denote the $\theta$-conjugation action of $R_{\mathbb{L}/\mathbb{K}} G_{\mathbb{L}}$ on itself, that is:
	$$\displaystyle \Ad_\theta(g)h:=gh\theta(g)^{-1},\;\;\; g,h\in R_{\mathbb{L}/\mathbb{K}} G_{\mathbb{L}}.$$
	Two $\theta$-strongly regular elements $\delta,\delta'\in (R_{\mathbb{L}/\mathbb{K}} G)(\mathbb{K}_v)=G(\mathbb{L}_v)$ (where $\mathbb{L}_v:=\mathbb{L}\otimes_{\mathbb{K}}\mathbb{K}_v$) are called {\em $\theta$-stably conjugate} if there exists $g\in (R_{\mathbb{L}/\mathbb{K}} G)(\mathbb{K}^{sep}_v)$ such that $\delta'=\Ad_\theta(g)\delta$. Let $f_v\in C_c^\infty(G(\mathbb{L}_v))$ be a test function. Then, for every $\theta$-strongly regular element $\delta\in G(\mathbb{L}_v)_{\theta-sreg}$ we can form the {\em $\theta$-twisted stable orbital integral}
	$$\displaystyle SO_{\delta,\theta}(f_v)=\sum_{\delta'} O_{\delta',\theta}(f_v)$$
	where the sum runs over the (finite set of) $\theta$-conjugacy classes of $\theta$-strongly regular elements $\delta'$ that are stably conjugate to $\delta$ and
	$$\displaystyle O_{\delta',\theta}(f_v)=\int_{T_{\delta',\theta}(\mathbb{K}_v)\backslash G(\mathbb{L}_v)} f_v(g^{-1}\delta' \theta(g)) \frac{dg}{dt_{\delta'}}$$
	denotes the corresponding $\theta$-twisted orbital integral of $f_v$. Here, $T_{\delta',\theta}$ denotes the twisted centralizer of $\delta'$ \eqref{tw centralizer} whereas $dg$ and $dt_{\delta'}$ are Haar measures on $G(\mathbb{L}_v)$ and $T_{\delta',\theta}(\mathbb{K}_v)$, respectively. Moreover, we ask that the measures $dt_{\delta'}$ are compatible in a similar way as above.

	Define the norm of $\delta\in G(\bbL_v)$ by
	$$\displaystyle N(\delta)=\delta \theta(\delta)\ldots \theta^{d-1}(\delta)\in G(\mathbb{L}_v).$$
	The following lemma is well-known over perfect fields \cite[\S 5]{K82}. For completeness, we provide a proof in general.

    \begin{lemma}
    Let $\delta\in G(\bbL_v)_{\theta-sreg}$ be $\theta$-strongly regular. Then, there exists $\gamma\in G_{sreg}(\bbK_v)$ such that $\gamma$ and $N(\delta)$ are conjugate inside $G(\bbK_v^{sep})$.
    \end{lemma}

    \begin{proof}
    The base-change of the twisted centralizer $T_{\delta,\theta}$ to $\bbK_v^{sep}$ can be naturally identified with the centralizer of $N(\delta)$ in $G_{\bbK_v^{sep}}$. Thus, the fact that $\delta$ is $\theta$-strongly regular implies that $N(\delta)\in G(\bbL_v)\subset G(\bbK_v^{sep})$ is strongly regular. Let $C\subset G_{\bbK_v^{sep}}$ be the conjugacy class of $x$. Since $\theta(N(\delta))=\delta^{-1}N(\delta)\delta$, $C$ is stable under the action of $\Gamma_{\bbK_v^{sep}}$ and is hence defined over $\bbK_v$. Therefore, it suffices to know that every strongly regular conjugacy class of $G$ that is defined over $\bbK_v$ has a point over $\bbK_v$. In the case of a perfect base-field this is proved in \cite[Lemma 2.1]{K82}. Moreover, the proof of {\em op. cit.} proceeds by a reduction to the case of $G$ semi-simple and simply-connected, where the result is a consequence of Steinberg's work \cite[Theorem 9.8]{Ste65} under the assumption that the base-field is perfect. However, the reduction argument works verbatim for non-perfect fields and the extension of Steinberg result to non-perfect fields was done in \cite[\S 8.6]{BS68}.
    \end{proof}
    
    For $\delta\in G(\bbL_v)_{\theta-sreg}$ and $\gamma\in G_{sreg}(\bbK_v)$, we say that $\gamma$ is {\em a norm of $\delta$}, which will be abbreviated by $\gamma\sim N(\delta)$, when $\gamma$ and $N(\delta)$ are conjugate inside $G(\bbK_v^{sep})$. When $\gamma$ is a norm of $\delta$ there is a canonical isomorphism of tori $T_\gamma\simeq T_{\delta,\theta}$ which we will use to normalize $SO_\gamma$ and $SO_{\delta,\theta}$ in a compatible way.
	
	We now come to the notion of transfer. Let $f_v\in C_c^\infty(G(\mathbb{L}_v))$. A {\em base-change transfer} of $f_v$ is any test function $f'_v\in C_c^\infty(G(\mathbb{K}_v))$ satisfying the following condition: for every strongly regular element $\gamma\in G_{sreg}(\mathbb{K}_v)$ we have the identity
	$$\displaystyle SO_\gamma(f'_v)=\left\{\begin{array}{ll}
		SO_{\delta,\theta}(f_v) & \mbox{ if } \gamma\sim N(\delta) \mbox{ for some } \delta\in G(\mathbb{L}_v)_{\theta-sreg}, \\
		0 & \mbox{ otherwise.}
	\end{array} \right.$$
	
	At some point the following reformulation of the notion of base-change transfer will be useful. Let
	$$\mathcal{Q}=G\sslash \Ad(G):=\Spec \bbK[G]^{\Ad(G)},\;\;\; \mathcal{Q}_\theta=R_{\bbL/\bbK}G_{\bbL}\sslash \Ad_\theta(R_{\bbL/\bbK}G_\bbL)$$
	be the GIT quotients of $G$ (resp. $R_{\bbL/\bbK}G_{\bbL}$) by the conjugation action (resp. $\theta$-conjugation action) on itself. Let $\pi: G\to \mathcal{Q}$ and $\pi_\theta: R_{\bbL/\bbK}G_{\bbL}\to \mathcal{Q}_\theta$ be the natural projections. Then, there is a natural isomorphism
	$$\mathcal{Q}\simeq \mathcal{Q}_\theta$$
	satisfying that a strongly regular element $\gamma\in G_{reg}(\bbK_v)$ is the norm of a $\theta$-strongly regular element $\delta\in (R_{\bbL/\bbK}G_\bbL)_{\theta-reg}(\bbK_v)$ if and only if $\pi(\gamma)=\pi_\theta(\delta)$. Furthermore, there exists an open subset $\mathcal{Q}_{sreg}\subset \mathcal{Q}$ such that $\pi^{-1}(\mathcal{Q}_{sreg})$ (resp. $\pi_\theta^{-1}(\mathcal{Q}_{sreg})$) is the subset $G_{sreg}$ (resp. $(R_{\bbL/\bbK} G_{\bbL})_{\theta-sreg}$) of strongly regular (resp. $\theta$-strongly regular) elements in $G$ (resp. $R_{\bbL/\bbK}G_\bbL$). Moreover, the maps between $\bbK_v$-points
	$$\pi_{\bbK_v}: G_{sreg}(\bbK_v)\to \mathcal{Q}_{sreg}(\bbK_v),\;\;\; \pi_{\theta,\bbK_v}: (R_{\bbL/\bbK}G_\bbL)_{\theta-sreg}(\bbK_v)\to \mathcal{Q}_{sreg}(\bbK_v)$$
	are smooth in the analytic topology (i.e. are locally trivial fibrations) over the regular locus $\mathcal{Q}_{sreg}(\mathbb{K}_v)$. In particular, for every pair of test functions $(f'_v,f_v)\in C_c^\infty(G(\bbK_v))\times C_c^\infty(G(\bbL_v))$, we can push forward the measures $f'_v(g'_v)dg'_v$,  $f_v(g_v) dg_v$ (where $dg'_v$, $dg_v$ are our chosen Haar measures on $G(\mathbb{K}_v)$ and $G(\bbL_v)$) through $\pi_{\mathbb{K}_v}$, $\pi_{\theta,\bbK_v}$ to obtain smooth measures $\pi_{*}(f'_vdg'_v)$ and $\pi_{\theta*}(f_vdg_v)$ on $\mathcal{Q}_{sreg}(\mathbb{K}_v)$. We then have the following characterization of base-change transfer.
	
	\begin{prop}\label{Prop char bc GIT}
		Let $(f_v,f'_v)\in C_c^\infty(G(\bbL_v))\times C_c^\infty(G(\bbK_v))$. Then, $f'_v$ is a base-change transfer of $f_v$ if and only if we have the following equality of measures on $\mathcal{Q}_{sreg}(\mathbb{K}_v)$:
		\begin{equation}
			\displaystyle \pi_{*}(f'_vdg'_v)=\pi_{\theta*}(f_vdg_v).
		\end{equation}
	\end{prop}

	Later, when discussing the $p$-descent, we will need a hypothesis on the description of the image of the base-change transfer. Such a result is known over local fields of characteristic zero (see \cite[Proposition 4.11]{MWstab1} for a general description of the image of endoscopic transfers) and should in all probability be part of a proof of the stabilization of the twisted trace formula for $(R_{\bbL/\bbK} G_{\bbL},\theta)$.
	
	\begin{hyp}\label{hyp image bc}
		Let $f'_v\in C_c^\infty(G(\bbK_v))$. Then, $f'_v$ is the base-change transfer of some function $f_v\in C_c^\infty(G(\bbL_v))$ if and only if the following condition is satisfied: for every strongly regular element $\gamma\in G_{sreg}(\mathbb{K}_v)$ that is not a norm (i.e. there does not exist $\delta\in G(\bbL_v)$ such that $\gamma \sim N(\delta)$), $SO_\gamma(f'_v)=0$.
	\end{hyp}
	
	We extend the notion of base-change transfer globally in the obvious way: we say that a factorizable test function $f'=\prod_v f'_v\in C_c^\infty(G(\mathbb{A}_{\mathbb{K}}))$ is a base-change transfer of $f=\prod_v f_v\in C_c^\infty(G(\mathbb{A}_{\mathbb{L}}))$ if for every place $v$, $f'_v$ is a base-change transfer of $f_v$. We recall that, by the fundamental lemma for base-change \cite{K86}, for almost all places $v$ of $\bbK$, $\mathbf{1}_{G(\mathcal{O}_{\mathbb{K}_v})}$ is a base-change transfer of $\mathbf{1}_{G(\mathcal{O}_{\mathbb{L}_v})}$. 
	
	Let $\rho: \Gamma_{\bbK}\to {}^L G(\overline{\bbQ}_{\ell})$ be our fixed $L$-parameter with Zariski dense image. We define a distribution $I_{\rho}: C_c^\infty(G(\bbA_{\bbK}))\to \overline{\bbQ}_{\ell}$ by
	$$\displaystyle I_\rho(f')=\Tr(R(f')\mid \CA_{0,\rho}(G_{\mathbb{K}}))$$
	where $\CA_{0,\rho}(G_{\mathbb{K}})$ stands for the direct sum of the cuspidal representations $\pi$ of $G_{\mathbb{K}}$ whose global Lafforgue parameter is conjugate to $\rho$ and $R(f')$ denotes the operator of right convolution by $f'$ i.e.
	$$\displaystyle (R(f')\phi)(x)=\int_{G(\mathbb{A}_K)} f'(g') \phi(xg') dg'.$$
	Note that this requires the choice of a global Haar measure $dg'$ and we implicitly assume henceforth that this is the product of the local Haar measures that were used in the definition of the stable orbital integrals. Note also that this $\rho$-subspace is most naturally defined over $\overline{\mathbb{Q}}_\lambda$, and transported to $\mathbb{C}$ using $\iota$. However, Proposition \ref{prop ext BHKT} shows that $\CA_{0,\rho}(G_{\mathbb{K}}))$ is in fact a sum of isotypic components of $\CA_{0}(G_{\mathbb{K}}))$ with respect to spherical Hecke operators, hence is defined over $\overline{\mathbb{Q}}$.
	
	Similarly, for any test function $f\in C_c^\infty(G(\mathbb{A}_{\mathbb{L}}))$, we define
	$$\displaystyle I_{\rho_{\mathbb{L}},\theta}(f)=\Tr(R(f)\times \theta\mid \CA_{0,\rho_{\mathbb{L}}}(G_{\mathbb{L}}))$$
	where $\CA_{0,\rho_{\mathbb{L}}}(G_{\mathbb{L}})$ denotes the direct sum of the cuspidal representations $\pi$ of $G_{\mathbb{L}}$ whose global Lafforgue parameter is conjugate to $\rho_{\mathbb{L}}=\rho\mid_{\Gamma_{\mathbb{L}}}$, $R(f)$ denotes the operator of right convolution by $f$ (using, once again, a global Haar measure that factors into the product of the chosen Haar measures on the $G(\mathbb{L}_v)$) and we have identified $\theta$ with its natural action on automorphic forms given by\footnote{Note that the subspace $\CA_{0,\rho_{\mathbb{L}}}(G_{\mathbb{L}})$ is invariant by $\theta$. This is because $\rho_{\mathbb{L}}$ comes from restriction from $\Gamma_{\mathbb{K}}$ so that the twisted $L$-parameter $\rho^\theta_{\mathbb{L}}$ is conjugate to $\rho_{\mathbb{L}}$.} 
	$$\displaystyle (\theta \phi)(x)=\phi(\theta(x)).$$ 
	
	We can now state our main result on comparison of trace formulas. It is conditional on the stabilization of (twisted) trace formulas, at least for functions supported in the strongly regular locus, as made explicit in the subsequent subsections. Recall that we are assuming that $\rho$ has Zariski dense image.
	
	\begin{theorem}\label{hyp rho traces}
		Assume the hypotheses \ref{hyp transfer}, \ref{hyp2}, \ref{hyp3} and \ref{hyp4} concerning the stabilization of twisted trace formulas. Let $f=\prod_v f_v\in C_c^\infty(G(\mathbb{A}_{\mathbb{L}}))$ and $f'=\prod_v f'_v\in C_c^\infty(G(\mathbb{A}_{\mathbb{K}}))$. Then, if $f'$ is a base-change transfer of $f$ we have
		\begin{equation}\label{iden rho traces}
		\displaystyle I_{\rho_{\mathbb{L}},\theta}(f)=I_\rho(f').
		\end{equation}
	\end{theorem}
	
	\subsection{Discrete parts of the trace formulas}
	
	In the rest of this section, we explain how Theorem \ref{hyp rho traces} can be deduced from the expected stabilizations of the trace formulas for $G_{\mathbb{K}}$ and of the twisted trace formula for $(R_{\mathbb{L}/\mathbb{K}}G_{\mathbb{L}},\theta)$. 
	More precisely, we aim to formulate analogs of the main results of \cite{Artstab} and \cite{MWstab2}. However, as we will recall, stabilization is an inductive process, so that to deduce any explicit information, we also need to assume the stabilization of the trace formula for all the (elliptic) endoscopic groups of $G_{\mathbb{K}}$, as well as for all the endoscopic groups of those endoscopic groups and so on. 
	Therefore, we start by discussing the stabilization of an arbitrary pair $(H,\Theta)$ where $H$ is a connected reductive group over $\mathbb{K}$ and $\Theta$ is an automorphism of $H$. 
	To simplify the discussion, we assume that $\Theta$ is of finite order and preserves a $\bbK$-pinning of $H$\footnote{We recall that a $\bbK$-pinning of $H$ is in a triple $(B,T,\{X_\alpha\})$ where $(B,T)$ is a Borel pair of $H_{\bbK^{sep}}$ and the $X_\alpha$'s are basis elements for the corresponding simple root lines in the Lie algebra and that is invariant by the obvious action of $\Gamma_{\bbK}$.} as well as that the center of $H$ is anisotropic (i.e. $A_H=1$). These assumptions will be satisfied in all the cases we will need to consider. Namely, we will either take $\Theta=1$ or $(H,\Theta)=(R_{\bbL/\bbK}G_{\bbL},\theta)$; since $G$ is quasi-split over $\bbK$, $\theta$  preserves a $\bbK$-pinning of $R_{\bbL/\bbK}G_{\bbL}$. Moreover, since $G$ is semisimple, the center of $G$ is finite, hence in particular anisotropic; and it is well-known that the center of an elliptic endoscopic group of a group with anisotropic center  again has anisotropic center.
	
	Let $\widehat{H}$ be the dual group of $H$ over $\overline{\bbQ}_{\ell}$. It comes with a pinning $(\widehat{B},\widehat{T}, \{ X_\alpha\})$ as well as an algebraic action of $\Gamma_{\mathbb{K}}$ fixing this pinning, or equivalently a continuous morphism $\Gamma_{\mathbb{K}}\to Out(\widehat{H})$ where $Out(\widehat{H})$ denotes the group of algebraic outer automorphisms of $\widehat{H}$. The automorphism $\Theta$ also induces, by functoriality, an automorphism $\widehat{\Theta}$ of $\widehat{H}$ fixing the pinning.
	
	Assume first that $\Theta=1$. Then the distribution to be stabilized is the discrete part of the trace formula for $H$ that can be written as (see \cite[\S X.5.1]{MWstab2} or \cite[\S 3.1]{Art13})
	\begin{equation}\label{spectral} 
		\displaystyle I^{H}_{\rm disc}(f^H) =
		\sum_{(M,w)} \lvert W(M)_w\rvert^{-1}|det(w-1)_{\mathfrak{a}_M }|^{-1}
		\Tr(M_{P}(w)\cI_{P,\disc,w}(f^H))
	\end{equation}
	$(f^H\in C_c^\infty(H(\bbA_\bbK)))$, where:
	\begin{itemize}
		\item the sum runs over the $H(\mathbb{K})$-conjugacy classes of pairs $(M,w)$ with $M\subset H$ a Levi subgroup (defined over $\mathbb{K}$) and $w$ an element of the relative Weyl group $W(M):=Norm_{H(\mathbb{K})}(M)/M(\mathbb{K})$ that is {\em regular} in the sense that it acts without nonzero fixed point on the vector space $\mathfrak{a}_M=X_*(A_M)\otimes \mathbb{Q}$;
		
		\item for every $(M,w)$, we have chosen a parabolic subgroup $P\subset H$ with Levi factor $M$ and $\cI_{P,\disc,w}$ denotes the direct sum of the (normalized) parabolically induced representations
		$$\displaystyle \cI_{P,\disc,w}=\bigoplus_\chi I_{P(\mathbb{A}_{\mathbb{K}})}^{H(\mathbb{A}_{\mathbb{K}})}(\CA_{\disc,\chi}(M))$$
		where $\chi$ runs over characters of $A_M(\mathbb{K})\backslash A_M(\mathbb{A}_{\mathbb{K}})$ that are invariant under $w$ and $\CA_{\disc,\chi}(M)$ is the subspace of square-integrable automorphic forms on $M$ with central character $\chi$.
		
		\item $M_{P}(w)$ denotes the standard intertwining operator
		$$\displaystyle M_{P}(w): I_{P(\mathbb{A}_{\mathbb{K}})}^{H(\mathbb{A}_{\mathbb{K}})}(\CA_{\disc}(M))\to I_{P(\mathbb{A}_{\mathbb{K}})}^{H(\mathbb{A}_{\mathbb{K}})}(\CA_{\disc}(M)).$$
	\end{itemize}
	
	When $\Theta$ is nontrivial, we will consider the discrete part of the $\Theta$-twisted trace formula for $H$. It is defined by a similar expression:
	\begin{equation}\label{spectral2} 
		\displaystyle I^{H,\Theta}_{\rm disc}(f^H) =
		\sum_{(M,\tilde{w})}\lvert W(M)_{\tilde{w}}\rvert^{-1} |\det(\tilde{w}-1)_{\mathfrak{a}_M }|^{-1}
		\Tr(M_{P}(\tilde{w})\cI_{P,\disc,\tilde{w}}(f^H))
	\end{equation}
	$(f^H\in C_c^\infty(H(\bbA_{\bbK})))$, where this time:
	\begin{itemize}
		\item The sum runs over $G(\bbK)$-conjugacy classes of pairs $(M,\tilde{w})$ with $M\subset H$ a Levi subgroup and $\tilde{w}$ an element of the relative Weyl set $\tilde{W}(M):=Norm_{G(\bbK)\Theta}(M)/M(\bbK)$ acting without nonzero fixed points on $\mathfrak{a}_M$.
		
		\item $\cI_{P,\disc,\tilde{w}}=\bigoplus_\chi I_{P(\mathbb{A}_{\mathbb{K}})}^{H(\mathbb{A}_{\mathbb{K}})}(\CA_{\disc,\chi}(M))$, $\chi$ runs over the set of characters of $A_M(\mathbb{K})\backslash A_M(\mathbb{A}_{\mathbb{K}})$ that are invariant under $\tilde{w}$ and 
        \[ M_{P}(\tilde{w}):I_{P(\mathbb{A}_{\mathbb{K}})}^{H(\mathbb{A}_{\mathbb{K}})}(\CA_{\disc}(M))\to I_{P(\mathbb{A}_{\mathbb{K}})}^{H(\mathbb{A}_{\mathbb{K}})}(\CA_{\disc}(M)) \]
        is again a standard intertwining operator.
	\end{itemize}
	
	To make the inductive process of stabilization more transparent, it is convenient to also consider pairs $\bbH=(H,\mathcal{H})$ where $H$ is, as before, a connected reductive group over $\bbK$ whereas $\mathcal{H}$ stands for an extension
	$$\displaystyle 1\to \widehat{H}\to \mathcal{H}\to W_{\mathbb{K}}\to 1,$$
	with $W_{\mathbb{K}}$ the Weil group of $\mathbb{K}$, such that the corresponding morphism $W_{\mathbb{K}}\to Out(\widehat{H})$ is the restriction of the morphism $\Gamma_{\mathbb{K}}\to Out(\widehat{H})$ coming from the (pinned) action of $\Gamma_{\bbK}$ on $\widehat{H}$. We will henceforth call such a pair {\em an $L$-pair}. The stabilization of the usual trace formula for $H$ will then correspond to the case where $\mathcal{H}={}^L H=\widehat{H}\rtimes W_{\mathbb{K}}$ is the $L$-group of $H$ (in its Weil form). The case of a general extension $\mathcal{H}$ can essentially be reduced to this case thanks to the theory of $z$-extensions (see \cite[\S 2.2]{KS}). Namely, it is known that there always exists a central extension $1\to C_1\to H_1\to H\to 1$ with $C_1$ an induced torus (i.e. a product of restrictions of scalars $R_{\mathbb{K}_i/\bbK} \bbG_m$ for some set of finite separable extensions $\bbK_i$ of $\bbK$) as well as an $L$-morphism $\mathcal{H}\subset {}^L H_1$ that restricts to the inclusion $\widehat{H}\subset \widehat{H}_1$ dual to the projection $H_1\to H$. The composition of this inclusion with the natural $L$-morphism ${}^L H_1\to {}^L C_1$ factors through $\mathcal{H}/\widehat{H}=W_{\bbK}$ and therefore induces an $L$-parameter $W_{\bbK}\to {}^L C_1$  which in turn corresponds to an automorphic character $\chi_1: C_1(\bbK)\backslash C_1(\bbA_{\bbK})\to \overline{\bbQ}_{\ell}^\times$. We then set $C_c^\infty(\bbH(\bbA_{\bbK})):=C_c^\infty(H_1(\bbA_{\bbK}),\chi_1^{-1})$ for the space of test functions on $H_1(\bbA_{\bbK})$ with central character $\chi_1^{-1}$, and define the discrete part of the corresponding trace formula by
	\begin{equation}\label{spectral3} 
		\displaystyle I^{\bbH}_{\rm disc}(f^{\bbH}) =
		\sum_{(M,w)} \lvert W(M)_w\rvert^{-1}|\det(w-1)_{\mathfrak{a}_M }|^{-1}
		\Tr(M_{P}(w)\cI_{P,\disc,w,\chi_1}(f^{\bbH}))
	\end{equation}
	$(f^{\bbH}\in C_c^\infty(\bbH(\bbA_\bbK)))$. Here the notation is as before, with respect to the group $H_1$ rather than $H$, except for $\cI_{P,\disc,w,\chi_1}$, which denotes the subrepresentation of $\cI_{P,\disc,w}$ with central character $\chi_1$. It can be shown (see \cite[I.2.5]{MWstab1}) that the pair $( C_c^\infty(\bbH(\bbA_\bbK)),I^{\bbH}_{\rm disc})$ is independent of the auxilliary datum $(H_1\to H, \cH\subset {}^L H_1)$ up to a unique isomorphism.
	
	Finally, we will also need to consider a variant of the above trace formulas where we project to a fixed family of Satake parameters. More precisely, let us pick Frobenius elements $\Frob_v\in \Gamma_{\bbK}$ for every place $v$ of $\bbK$. Then, for every finite set of places $S$ of $\bbK$ outside of which $H$ is unramified, we let $\mathcal{C}_S(H)$ be the set of all families $c=(c_v)_{v\notin S}$ of Satake parameters outside of $S$ -- i.e. of semisimple $\widehat{H}$-conjugacy classes $c_v\in \widehat{H}\Frob_v$ for all $v\notin S$. We then define $\mathcal{C}(H)=\varinjlim_S \mathcal{C}_S(H)$, the transition maps being the natural projections. For every $c\in \mathcal{C}(H)$ and $P$, $w$ as in \eqref{spectral} we let $\cI_{P,\disc,w,c}$ be the direct sum of all the irreducible summands of $\cI_{P,\disc,w}$ that are unramified with Satake parameters $c_v$ outside of some finite set of places $S$ (that may depend on the representation). We then define the {\em $c$-part} $I^{H}_{\rm disc,c}$ of $I^{H}_{\rm disc}$ by the same formula \eqref{spectral} up to replacing $\cI_{P,\disc,w}$ by $\cI_{P,\disc,w,c}$. We define similarly the $c$-part $I^{H,\Theta}_{\rm disc,c}$ of $I^{H,\Theta}_{\rm disc}$. Note that we have a decomposition
	$$\displaystyle I^{H,\Theta}_{\rm disc}=\sum_{c\in \mathcal{C}(H)} I^{H,\Theta}_{\rm disc,c}.$$
	We emphasize that the right hand side typically contains infinitely many terms but that for every test function $f\in C_c^\infty(H(\bbA_K))$ only finitely many of the $I^{H,\Theta}_{\rm disc,c}(f)$ are nonzero. More precisely, this only depends on the level of $f$: for every compact-open subgroup $U\subset H(\bbA_K)$, there exists a finite subset $\mathcal{C}_U\subset \mathcal{C}(H)$ such that $I^{H,\Theta}_{\rm disc,c}(f)=0$ for every $f\in C_c(U\backslash H(\bbA_{\bbK})/U)$ and $c\notin \mathcal{C}_U$. 
	
	Such a decomposition also exists for general pairs $\bbH=(H,\cH)$ as above. More precisely, we define $\mathcal{C}_S(\bbH)$ as the set of families $c=(c_v)_{v\notin S}$ of $\widehat{H}$-conjugacy classes of semisimple elements $c_v\in \mathcal{H}$ that map to the chosen Frobenius elements $\Frob_v\in \Gamma_{\bbK}$. Then it can be shown that the colimit $\mathcal{C}(\bbH)=\varinjlim_S \mathcal{C}_S(\bbH)$ does not depend on the chosen Frobenius representatives (up to a unique bijection) and that we have a similar decomposition
	$$\displaystyle I^{\bbH}_{\rm disc}=\sum_{c\in \mathcal{C}(H)} I^{\bbH}_{\rm disc,c}$$
	coming from the natural embedding $\cC(\bbH)\subset \cC(H_1)$ for any chosen auxilliary datum $(H_1\to H,\cH\subset {}^L H_1)$.

	\begin{prop}\label{prop Idiscc}
		$I^{H}_{\disc,c}=0$ unless there exists an $L$-parameter $\phi: \Gamma_{\bbK}\to {}^L H(\bbQ_\ell)$ such that $c_v$ is conjugate to $\phi(\Frob_v)$ for almost all $v$. More generally, $I^{\bbH}_{\rm disc,c}=0$ unless there exists an $L$-parameter $\phi: \Gamma_{\bbK}\to \mathcal{H}$ such that $c_v$ is conjugate to $\phi(\Frob_v)$ for almost all $v$.
	\end{prop}
	
	\begin{proof}
		This is essentially a direct consequence of V. Lafforgue's work. More precisely, let 
        \[ (H_1\to H,\cH\subset {}^L H_1) \]
        be an auxiliary datum for $\bbH$ and 
        \[ \chi_1: C_1(\bbK)\backslash C_1(\bbA_{\bbK})\to \overline{\bbQ}_{\ell}^\times \]
        the corresponding automorphic character of $C_1=\Ker(H_1\to H)$. Then, all the automorphic representations appearing in the spectral expansion \eqref{spectral3} of $I^{\bbH}_{\rm disc}$ are parabolically induced from discrete automorphic representations of Levi subgroups of $H_1$. Moreover, by Langlands spectral theorem, discrete automorphic representations are themselves quotients of parabolic inductions of cuspidal representations. Therefore, $I^{\bbH}_{\rm disc}$ can be written as a sum of characters of automorphic representations, each appearing in the parabolic induction of some cuspidal representation of a Levi subgroup $M_1$ of $H_1$. By V. Lafforgue's automorphic to Galois result (i.e.\ \cite[Th\'eor\`eme 11.11]{Lafa}), we can attach to each of these cuspidal representations an $L$-parameter $W_{\bbK}\to {}^L M_1$ with a global-local compatibility at almost all unramified places. Composing with the inclusions ${}^L M_1\subset {}^L H_1$, we see that we can also attach global $L$-parameters to their parabolic inductions (and their subquotients). More precisely, if $\pi_1$ is an automorphic representation of $H_1$ appearing in the spectral decomposition of $I^{\bbH}_{\rm disc}$, there exists a global $L$-parameter $\phi:W_{\bbK}\to {}^L H_1$ such that the Satake parameter of $\pi_{1,v}$ is $\widehat{H}_1$-conjugate to $\phi(\Frob_v)$ for almost all $v$. Moreover, by construction the central character of $\pi_1$ restricts to $\chi_1$ on $C_1(\bbK)\backslash C_1(\bbA_{\bbK})$. This fact is reflected on the dual side by the following property: the composition of $\phi$ with the $L$-morphism ${}^L H_1\to {}^L C_1$, dual to the central inclusion $C_1\subset H_1$, is the $L$-parameter ${}^L \chi_1$ of $\chi_1$. Indeed, this can be deduced from a purely local unramified computation, namely that for almost all $v$ the projection of the Satake parameter of $\pi_{1,v}$ to ${}^L C_1$ is the Satake parameter of the central character of $\pi_{1,v}$ on $C_1(\bbK_v)$, since two global $L$-parameters $\psi,\psi':W_{\bbK}\to {}^L C_1$ that are conjugate at almost all places are actually globally conjugate (a consequence of Cebotarev theorem since $C_1$ is an induced torus). Noting that the embedding $\cH\subset {}^L H_1$ identifies $\cH$ with the fiber product ${}^L H_1\times_{{}^L C_1} W_{\bbK}$, where the morphism $W_{\bbK}\to {}^L C_1$ is ${}^L \chi_1$, it follows that up to conjugation we may assume that $\phi$ has image in $\cH$. The proposition then readily follows.
	\end{proof}
	
	\subsection{Elliptic endoscopic data}

	Recall (see e.g. \cite[VI.3.1]{MWstab2}) that an {\em endoscopic datum} for $(H,\Theta)$ is a quadruple $\underline{\bbH}'=(H',\mathcal{H}',s,\eta)$ where:
	\begin{itemize}
		\item $H'$ a quasi-split reductive group over $\bbK$;
		
		\item $\mathcal{H}'$ a closed subgroup of ${}^L H$ that surjects onto the Weil group $W_{\bbK}$;
		
		\item $\eta$ is an algebraic isomorphism $\widehat{H}'\simeq (\mathcal{H}')^0$ (where $(.)^0$ stands for the neutral component) so that $\mathcal{H}'$ sits in an exact sequence $1\to \widehat{H}'\to \mathcal{H}'\to W_{\bbK}\to 1$. We require that the morphism $W_{\bbK}\to Out(\widehat{H}')$ corresponding to this extension coincides with the one coming from the action of $\Gamma_{\bbK}$ on $\widehat{H}'$. In other words, $\bbH':=(H',\cH')$ is an $L$-pair in the terminology introduced in the previous subsection;
		
	\item $s$ is a semisimple element in the coset $\widehat{H}\widehat{\Theta}$ (recall that $\widehat{\Theta}$ denotes the dual automorphism of $\widehat{H}$) such that $(\mathcal{H}')^0=Z_{\widehat{H}}(s)^0$ is the neutral component of the centralizer of $s$ in $\widehat{H}$ and there exists a $1$-cocycle $a\in Z^1(W_{\bbK},Z(\widehat{H}))$ that is everywhere locally cohomologically trivial such that
		\begin{equation}
			\displaystyle (g,w)s(g,w)^{-1}=a(w)s \mbox{ for every } (g,w)\in \mathcal{H}'.
		\end{equation}
	\end{itemize}

	Moreover, for two endoscopic data $\underline{\bbH}'_1=(H'_1,\mathcal{H}'_1,s_1,\eta_1)$ and $\underline{\bbH}'_1=(H'_2,\mathcal{H}'_2,s_2,\eta_2)$ of $(H,\Theta)$, we say that:
	\begin{itemize}
		\item $\underline{\bbH}'_1$ and $\underline{\bbH}'_2$ are {\em equivalent} if there exists $g\in \widehat{H}$ such that $g\cH'_1g^{-1}=\cH'_2$ and $gs_1g^{-1}\in Z(\widehat{H})s_2$;
		
		\item $\underline{\bbH}'$ is {\em elliptic} if $Z(\widehat{H}')^{\Gamma_{\bbK}}/Z(\widehat{H})^{\Gamma_{\bbK},\widehat{\Theta}}$ is finite.
	\end{itemize}
	We denote by $\mathcal{E}_{ell}(H,\Theta)$ the set of equivalence classes of elliptic endoscopic data for $(H,\Theta)$.
	
	The above definition (when $\Theta=1$) naturally extends to an $L$-pair $\bbH=(H,\mathcal{H})$. Indeed, we just have to replace the $L$-group of $H$ by $\cH$ in the definition. We again denote by $\mathcal{E}_{ell}(\bbH)$ the set of equivalence classes of elliptic endoscopic data for $\bbH$.
	
	We will mostly be interested in a particular endoscopic group for $(R_{\bbL/\bbK}G_{\bbL},\theta)$. Namely, we may identify the $L$-group of $R_{\bbL/\bbK}G_{\bbL}$ with $\widehat{G}^{\Gal(\bbL/\bbK)}\rtimes W_{\bbK}$, where $W_{\bbK}$ acts on $\widehat{G}^{\Gal(\bbL/\bbK)}$ by
	$$\displaystyle \tau\cdot (g_\sigma)_{\sigma\in \Gal(\bbL/\bbK)}=(\tau g_{\sigma \tau})_{\sigma\in \Gal(\bbL/\bbK)}.$$
	Then, the {\em principal endoscopic group} of $(R_{\bbL/\bbK}G_{\bbL},\theta)$ is given by the quadruple $(G_{\bbK}, {}^L G_{\bbK},\widehat{\theta},b_{\bbL/\bbK})$ where $\widehat{\theta}$ is the dual automorphism of $\widehat{G}^{\Gal(\bbL/\bbK)}$, given by permutation of the factors through the translation by $\theta$ on $\Gal(\bbL/\bbK)=\langle \theta\rangle$, and
	$$\displaystyle b_{\bbL/\bbK}: {}^L G_{\bbK}=\widehat{G}\rtimes W_{\bbK}\to {}^L R_{\bbL/\bbK}G_{\bbL}=\widehat{G}^{\Gal(\bbL/\bbK)} \rtimes W_{\bbK}$$
	is the base-change homomorphism, which is the identity on $W_{\bbK}$ and otherwise restricts to the diagonal embedding $\widehat{G}\hookrightarrow \widehat{G}^{\Gal(\bbL/\bbK)}$.

	\begin{lemma}\label{lem endobc}
		Let $\rho: \Gamma_{\bbK}\to {}^LG(\overline{\bbQ}_{\ell})$ be an $L$-parameter with Zariski dense image. Then, up to equivalence, the only (twisted) elliptic endoscopic group $\underline{\mathbb{H}}=(H,\cH,s,\eta)\in \mathcal{E}_{ell}(R_{\bbL/\bbK}G_{\bbL},\theta)$ such that $\cH$ contains a conjugate of the image of $b_{\bbL/\bbK}\circ \rho$ is the principal endoscopic group $(G_{\bbK},{}^L G_{\bbK}, \theta, b_{\bbL/\bbK})$.
	\end{lemma}
	
	\begin{proof}
		Up to replacing $\underline{\mathbb{H}}$ by an equivalent endoscopic group, we may assume that $b_{\bbL/\bbK}\circ \rho$ has image in $\mathcal{H}$. As $\rho$ has Zariski dense image, this then implies that $\mathcal{H}$ contains the image of $b_{\bbL/\bbK}$, that is $\widehat{G}(\overline{\bbQ}_{\ell})^{\diag}\rtimes W_{\bbK}$. In particular, the element $s$ must commute with $\widehat{G}(\overline{\bbQ}_{\ell})^{\diag}$. However, the centralizer of this latter group in $\widehat{G}^{\Gal(\bbL/\bbK)}\widehat{\theta}$ is readily seen to be $Z(\widehat{G})^{\Gal(\bbL/\bbK)} \widehat{\theta}$. Thus, up to translating $s$ by an element in the center, we may assume that $s=\widehat{\theta}$. However, the centralizer of $\widehat{\theta}$ in $\widehat{G}^{\Gal(\bbL/\bbK)}$ is precisely $\widehat{G}^{\diag}$ which forces $\cH$ to be equal to $\widehat{G}(\overline{\bbQ}_{\ell})^{\diag}\rtimes W_{\bbK}$.
	\end{proof}
	
	\subsection{Assumptions on the stabilization of trace formulas}
	
	We continue with the setting of the previous two subsections. Before stating properly our hypotheses about the stabilization of the (discrete part of the) trace formula, it remains to introduce one last ingredient:  namely the Langlands-Kottwitz-Shelstad (henceforth abbreviated to LKS) transfer. More precisely, let $(H,\Theta)$ be a twisted pair and $\underline{\bbH}'=(H',\cH',s,\eta)\in \cE_{ell}(H,\Theta)$ be an elliptic endoscopic group of $(H,\Theta)$. Let $\bbH':=(H',\cH')$ be the corresponding $L$-pair. In particular, we can pick an auxilliary datum $(H'_1\to H',\cH'\subset {}^L H'_1)$ relative to this $L$-pair. The LKS transfer for the endoscopic group $\underline{\bbH}'$ then takes the form of a correspondence
	$$\displaystyle f\in C_c^\infty(H(\bbA_{\bbK}))\to f^{\underline{\bbH}'}\in C_c^\infty(\bbH'(\bbA_{\bbK}))=C_c^\infty(H_1'(\bbA_{\bbK}),\chi_1^{-1}),$$
	that is characterized by identities of the form
	\begin{equation}\label{transfer eq}
		\displaystyle SO_{\gamma}(f^{\underline{\bbH}'})=\sum_{\delta} \Delta(\gamma,\delta) O_{\delta,\Theta}(f)
	\end{equation}
	where $\gamma$ is any strongly regular stable conjugacy class in $H'_1(\bbA_{\bbK})$, that is the non-empty intersection of a strongly regular $H'_1(\bbA_{\bbK}\otimes_{\bbK}\overline{\bbK})$-conjugacy class with $H'_1(\bbA_{\bbK})$,
	$$\displaystyle SO_{\gamma}(f^{\underline{\bbH}'})=\sum_{\gamma'} O_{\gamma'}(f^{\underline{\bbH}'})=\sum_{\gamma'} \int_{H_{\gamma'}(\bbA_{\bbK})\backslash H(\bbA_{\bbK})} f^{\underline{\bbH}'}(h^{-1}\gamma' h) \frac{dh}{dh_{\gamma'}}$$
	denotes the corresponding stable orbital integral of $f^{\underline{\bbH}'}$, that is the sum of its orbital integrals over the $H'_1(\bbA_{\bbK})$-conjugacy classes $\gamma'$ inside $\gamma$, the sum in the RHS of \eqref{transfer eq} is over $\Theta$-strongly regular conjugacy classes in $H(\bbA_{\bbK})$, $O_{\delta,\Theta}(f)$ stands for the $\Theta$-twisted orbital integral of $f$ at $\delta$ and $\Delta(\gamma,\delta)$ is the {\em global transfer factor} defined by Langlands, Shelstad and Kottwitz.
	
	Note that the transfer $f^{\underline{\bbH}'}$ is only defined through its stable orbital integrals $SO_{\gamma}(f^{\underline{\bbH}'})$ for $\gamma$ strongly regular. In other terms, the transfer, if it exists, defines a unique element of $SI(\bbH'(\bbA_{\bbK}))$ where the latter denotes the quotient of $C_c^\infty(\bbH'(\bbA_{\bbK}))$ by the subspace of functions all of whose strongly regular stable orbital integrals vanish. 
	
	Note also that here we have directly defined LKS transfer at a global level. Usually, the transfer is defined first locally (with the important caveat that the transfer factors are not canonical anymore, thus requiring the choice of auxiliary rigidifying data at all places in a coherent way) and the global transfer is a consequence of the local transfer together with the endoscopic fundamental lemma. Thus, in some sense, our presentation subsumes both the local transfer and fundamental lemma in one statement.
	
	In the particular case where $(H,\Theta)=(R_{\bbL/\bbK} G_{\bbL},\theta)$ and $\underline{\bbH}'$ is the principal endoscopic group $(G_{\bbK},{}^L G_{\bbK}, \theta, b_{\bbL/\bbK})$, the transfer map $f\mapsto f^{\underline{\bbH}'}$ coincides with the base-change transfer reviewed in Subsection \ref{subsection rhotrace}.
	
	Thanks to the work of Waldspurger \cite{Waltransfer} and Ng\^o \cite{Ngo}, the existence of LKS transfers is now a theorem over number fields. However, it is currently not known over function fields and should be assumed as a separate hypothesis:
	
	\begin{hyp}\label{hyp transfer}
		For every twisted pair $(H,\Theta)$ with either $\Theta=1$ and $H$ quasi-split with anisotropic center or $(H,\Theta)=(R_{\bbL/\bbK} G_{\bbL}, \theta)$, and every elliptic endoscopic group $\underline{\bbH}'\in \cE_{ell}(H,\Theta)$, there exists a LKS transfer map
		$$\displaystyle f\in C_c^\infty(H(\bbA_{\bbK}))\mapsto f^{\underline{\bbH}'}\in SI(\bbH'(\bbA_{\bbK}))$$
		i.e. a map satisfying the identity \eqref{transfer eq}.
	\end{hyp} 
	
	In the case of an $L$-pair $\bbH=(H,\cH)$, we can formulate a similar hypothesis on the existence of a LKS transfer $f\in C_c^\infty(\bbH(\bbA_{\bbK}))\mapsto f^{\underline{\bbH}'}\in SI(\bbH'(\bbA_{\bbK}))$ for every $\underline{\bbH}'\in \cE_{ell}(\bbH)$. However, it can be shown that this seemingly more general hypothesis is implied by the existence of LKS transfers for $H_1$ where $(H_1\to H,\cH\subset {}^L H_1)$ is any auxilliary datum, essentially by the process of ``projecting everything to the central character $\chi_1$''. This variant of endoscopic transfer will however be used in the formulation of our next hypothesis which concerns the stabilization of the trace formula. 
	
	Before the statement, let us introduce the following notation: for every $L$-pair $\bbH=(H,\cH)$ over $\bbK$, with $H$ quasi-split and having anisotropic center, and $\underline{\bbH}'=(H',\cH',s,\eta)\in \cE_{ell}(\bbH)$, we set
	$$\displaystyle \iota(\bbH,\underline{\bbH}'):=\lvert \Out(\underline{\bbH}')\rvert^{-1} \tau(H)\tau(H')^{-1}$$
	where 
    \[ \tau(H)=\lvert Z(\widehat{H})^{\Gamma_{\bbK}}\rvert \lvert \ker^1(\bbK,Z(\widehat{H}))\rvert^{-1} \]
    and 
    \[ \tau(H')=\lvert Z(\widehat{H}')^{\Gamma_{\bbK}}\rvert \lvert \ker^1(\bbK,Z(\widehat{H}'))\rvert^{-1} \]
    are the Tamagawa numbers of $H$ and $H'$ respectively and $\Out(\underline{\bbH}')$ is the group of outer automorphisms of the endoscopic datum defined by
	$$\displaystyle \Out(\underline{\bbH}')=\left\{ h\in \widehat{H}\mid h\cH' h^{-1}=\cH',  hsh^{-1}\in Z(\widehat{H})s\right\}/\widehat{H}'.$$
	
	\begin{hyp}\label{hyp2}
		For every $L$-pair $\bbH=(H,\cH)$ over $\bbK$ with $H$ quasi-split and having anisotropic center, there exists a stable distribution $S^{\bbH}_{\disc}: SI(\bbH(\bbA_{\bbK}))\to \overline{\bbQ}_{\ell}$ such that the following identities
		\begin{equation}\label{stabilization1}
			\displaystyle I_{\disc}^{\bbH}(f)=\sum_{\underline{\bbH}'\in \mathcal{E}_{ell}(\bbH)} \iota(\bbH, \underline{\bbH}') S^{\bbH'}_{\disc}(f^{\underline{\bbH}'})
		\end{equation}
		are satisfied for every $L$-pair $\bbH=(H,\cH)$ as before and $f\in C_c^\infty(\bbH(\bbA_{\bbK}))$.
	\end{hyp}
	
	\begin{remark}
		\begin{itemize}
			\item Implicit in this statement is a recursive definition of the stable distributions $S^{\bbH}_{\disc}$. Indeed, identity \eqref{stabilization1} can be rewritten as
			\begin{equation*}
				S^{\bbH}_{\disc}(f)=I_{\disc}^{\bbH}(f)-\sum_{\underline{\bbH}'\in \mathcal{E}_{ell}(\bbH)\setminus \{ \bbH\}} \iota(\bbH, \underline{\bbH}') S^{\bbH'}_{\disc}(f^{\underline{\bbH}'})
			\end{equation*}
			and we can assume, by induction on the dimension of the group for example, that the stable distributions $S^{\bbH'}_{\disc}$ appearing in the right hand side have already been defined. This can be turn into a definition of $S^{\bbH}_{\disc}$ and, from this point of view, the content of the hypothesis is that the distribution $S^{\bbH}_{\disc}$ is always stable, that is factors through $SI(\bbH(\bbA_{\bbK}))$.
			
			\item Hypothesis \ref{hyp2} can be reduced to the more usual stabilization of the trace formula in the pinned case (i.e. when $\cH={}^L H$), by projecting afterwards onto a given central character. In the pinned case, the above hypothesis is the exact analog of \cite[Global Theorem 2]{Artstab} as revisited in \cite[\S 3.2]{Art13}.
			
			\item The assumption that $H$ is quasi-split with anisotropic center is for convenience only. More precisely, the anisotropic center condition simplifies the constant $\iota(\bbH,\underline{\bbH}')$. Moreover, we eventually want to apply the stabilization only to the semi-simple quasi-split group $G$ and its endoscopic groups and they all satisfy this assumption.
			
		\end{itemize}
	\end{remark}
	
	For our purpose, it will be more convenient to use a variant of the stabilization after projection to a given family of Satake parameters $c\in \cC(\bbH)$. This variant, that we formulate as an hypothesis below, can be shown to follow from Hypothesis \ref{hyp2} and the fundamental lemma for spherical Hecke algebras see \cite[Lemma 3.3.1]{Art13}. In the statement below, we make implicit use of the map $\cC(\bbH')\to \cC(\bbH)$, for any $\underline{\bbH}'\in \cE_{ell}(\bbH)$, that is induced from the $L$-embedding $\cH'\subset \cH$.
	
	\begin{hyp}\label{hyp3}
		For every $L$-pair $\bbH=(H,\cH)$ over $\bbK$, with $H$ quasi-split and having anisotropic center, and every $c\in \cC(\bbH)$, there exists a stable distribution $S^{\bbH}_{\disc,c}: SI(\bbH(\bbA_{\bbK}))\to \overline{\bbQ}_{\ell}$ such that the following identities
		\begin{equation}\label{stabilization2}
			\displaystyle I_{\disc,c}^{\bbH}(f)=\sum_{\underline{\bbH}'\in \mathcal{E}_{ell}(\bbH)} \iota(\bbH, \underline{\bbH}') \sum_{c'\in \cC(\bbH')\mapsto c} S^{\bbH'}_{\disc,c'}(f^{\underline{\bbH}'})
		\end{equation}
		are satisfied for every $L$-pair $\bbH=(H,\cH)$ as before, $c\in \cC(\bbH)$ and $f\in C_c^\infty(\bbH(\bbA_{\bbK}))$.
	\end{hyp}
	
	Finally, we will also need to assume the stabilization of the twisted trace formula for $(R_{\bbL/\bbK}G_{\bbL},\theta)$. This is the content of the next hypothesis, where for convenience we directly state its variant after projection onto a family of Satake parameters. Note that in the case of number fields the analog of this hypothesis is again a theorem thanks to the work of M\oe{}glin-Waldspurger \cite{MWstab1} \cite{MWstab2}.
	
	\begin{hyp}\label{hyp4}
		Assume Hypothesis \ref{hyp3}. Then, we have the following identity
		\begin{equation}\label{stabilization3}
			\displaystyle I_{\disc,c}^{R_{\bbL/\bbK}G_{\bbL},\theta}(f)=\sum_{\underline{\bbH}\in \mathcal{E}_{ell}(R_{\bbL/\bbK}G_{\bbL},\theta)} \iota((R_{\bbL/\bbK}G_{\bbL},\theta), \underline{\bbH}) \sum_{c'\in \cC(\bbH)\mapsto c} S^{\bbH}_{\disc,c'}(f^{\underline{\bbH}})
		\end{equation}
		for every $c\in \cC(R_{\bbL/\bbK}G_{\bbL})^\theta$ and $f\in C_c^\infty(G(\bbA_{\bbL}))$.
	\end{hyp}
	
	\subsection{Proof of Theorem \ref{hyp rho traces}}\label{proof theo rho traces}

	We will use the following lemma, which is a direct consequence of Proposition \ref{prop Idiscc} using induction.
	
	\begin{lemma}\label{lem Sdiscc}
		For every $L$-pair $\bbH=(H,\cH)$, with $H$ quasi-split and having anisotropic center, and every $c\in \cC(\bbH)$, we have $S^{\bbH}_{\rm disc,c}=0$ unless there exists an $L$-parameter $\phi: \Gamma_{\bbK}\to \mathcal{H}$ such that $c_v$ is conjugate to $\phi(\Frob_v)$ for almost all $v$.
	\end{lemma}
	
	For every $L$-parameter $\rho': \Gamma_{\bbK}\to {}^LG(\overline{\bbQ}_{\ell})$, we denote by $c_{\rho'}\in \cC(G)$ the family of Satake parameters given by $c_{\rho',v}=\rho'(\Frob_v)$ for almost all $v$. Let $\rho: \Gamma_{\bbK}\to {}^LG(\overline{\bbQ}_{\ell})$ be our $L$-parameter with Zariski dense image. By Proposition \ref{prop ext BHKT}, for every  $\rho': \Gamma_{\bbK}\to {}^LG(\overline{\bbQ}_{\ell})$, we have $c_{\rho'}=c_{\rho}$ if and only if $\rho'$ is conjugate to $\rho$. This implies that
	$$\displaystyle I_{\disc,c_{\rho}}^{G_{\bbK}}(f)=I_{\rho}(f) \mbox{ for every } f\in C_c^\infty(G(\bbA_{\bbK})),$$
	but also, by Lemma \ref{lem Sdiscc}, that for every $\underline{\bbH}\in \cE_{ell}(G_{\bbK})$ and $c\in \cC(\bbH)$ mapping to $c_{\rho}$, we have either $S^{\bbH}_{\rm disc,c}=0$ or $\underline{\bbH}=G_{\bbK}$ and $c=c_{\rho}$. Therefore, Hypothesis \ref{hyp3} applied to $c=c_\rho$ leads to the identity
	\begin{equation}\label{eq Irhostable}
		\displaystyle I_{\rho}(f)=S^{G_{\bbK}}_{\rm disc,c_{\rho}}(f) \mbox{ for every } f\in C_c^\infty(G(\bbA_{\bbK})).
	\end{equation}
	
	Note that we may identify $c_{\rho_{\bbL}}$ with an element of $\cC(R_{\bbL/\bbK}G_{\bbL})=\cC(G_{\bbL})$. Since $\rho_{\bbL}$ still has Zariski dense image, the above reasoning also shows
	$$\displaystyle I_{\disc,c_{\rho_{\bbL}},\theta}^{R_{\bbL/\bbK}G_{\bbL}}(f)=I_{\rho_{\mathbb{L}},\theta}(f) \mbox{ for every } f\in C_c^\infty(G(\bbA_{\bbL})).$$
	On the other hand, by Lemma \ref{lem Sdiscc} and Lemma \ref{lem endobc}, we have $S^{\bbH}_{\disc,c'}=0$ for every $\underline{\bbH}\in \mathcal{E}_{ell}(R_{\bbL/\bbK}G_{\bbL},\theta)$ and $c'\in \cC(\bbH)$ that maps to $c_{\rho_{\bbL}}$ unless $\underline{\bbH}$ is the principal endoscopic group $(G_{\bbK},{}^L G_{\bbK}, \theta, b_{\bbL/\bbK})$ and $c'=c_{\rho'}$ where $\rho': \Gamma_{\bbK}\to {}^L G(\overline{\bbQ_\ell})$ is a $L$-parameter with $b_{\bbL/\bbK}\circ \rho'$ equivalent to $b_{\bbL/\bbK}\circ \rho$.
	
	\begin{lemma}
		Up to equivalence, the $L$-parameters $\rho': \Gamma_{\bbK}\to {}^L G(\overline{\bbQ_\ell})$ with $b_{\bbL/\bbK}\circ \rho'$ equivalent to $b_{\bbL/\bbK}\circ \rho$ are exactly the twists $\widehat{\chi}\rho$ where $\widehat{\chi}:\Gal(\bbL/\bbK)\to Z(\widehat{G})^{\Gamma_{\bbL}}$ is a $1$-cocyle. Moreover, for two such $1$-cocycles $\widehat{\chi}$, $\widehat{\chi}'$, the $L$-parameters $\widehat{\chi}\rho$ and $\widehat{\chi}'\rho$ are equivalent if and only if the images of $\widehat{\chi}$ and $\widehat{\chi}'$ in $H^1(\bbL/\bbK,Z(\widehat{G})^{\Gamma_{\bbL}})$ are the same.
	\end{lemma}
	
	\begin{proof}
		First, we note that the $L$-parameters $b_{\bbL/\bbK}\circ \rho'$ and $b_{\bbL/\bbK}\circ \rho$ are equivalent if and only if the restrictions $\rho'_{\bbL}$ and $\rho_{\bbL}$ of $\rho$ and $\rho'$ to $\Gamma_{\bbL}$ respectively are equivalent as $L$-parameters $\Gamma_{\bbL}\to {}^L G_{\bbL}(\overline{\bbQ_{\ell}})$. Since $\rho_{\bbL}$ still has Zariski dense image, its centralizer in $\widehat{G}$ is reduced to the $\Gamma_{\bbL}$-fixed part of the center $Z(\widehat{G})^{\Gamma_{\bbL}}$. This readily implies that $\rho'_{\bbL}$ and $\rho_{\bbL}$ are equivalent if and only if $\rho'$ and $\rho$ differ, up to equivalence, by a $1$-cocycle $\Gal(\bbL/\bbK)\to Z(\widehat{G})^{\Gamma_{\bbL}}$. Similarly, $\widehat{\chi}\rho$ and $\widehat{\chi}'\rho$ are equivalent if and only if they are conjugate by an element of $Z(\widehat{G})^{\Gamma_{\bbL}}$, i.\ e.\ if and only if $\widehat{\chi}$ and $\widehat{\chi}'$ have the same image in $H^1(\bbL/\bbK,Z(\widehat{G})^{\Gamma_{\bbL}})$.
	\end{proof}
	
	By the above lemma, together with identity \eqref{eq Irhostable} (applied to the various twists $\widehat{\chi}\rho$ instead of $\rho$), Hypothesis \ref{hyp4} applied to $c=c_{\rho_{\bbL}}$ reads
	\begin{equation}\label{Irho1}
		\displaystyle I_{\rho_{\mathbb{L}},\theta}(f)=\iota((R_{\bbL/\bbK} G_{\bbL},\theta), \underline{G}_{\bbK})\sum_{\hat{\chi}\in H^1(\bbL/\bbK,Z(\widehat{G})^{\Gamma_{\bbL}})} I_{\hat{\chi}\rho}(f')
	\end{equation}
	for every $f\in C_c^\infty(G(\bbA_{\bbL}))$ with base-change transfer $f'\in C_c^\infty(G(\bbA_{\bbK}))$. We henceforth fix $f$ and $f'$.
	
	By a result of Langlands (see e.g. \cite{LLap15}), there is a natural bijection
	$$\displaystyle H^1(\bbK,Z(\widehat{G}))/\ker^1(\bbK,Z(\widehat{G}) )\simeq \Hom_{sm}(G(\bbK)\backslash G(\bbA_{\bbK}),\overline{\bbQ}_{\ell}^\times),\;\; \widehat{\chi}\mapsto \chi,$$
	where the right hand side denotes the group of smooth characters $G(\bbA_{\bbK})\to\overline{\bbQ}_{\ell}^\times$ that are trivial on $G(\bbK)$. By Proposition \ref{prop ext BHKT}, we have $\ker^1(\bbK,Z(\widehat{G}))=1$. Moreover, the natural morphism $H^1(\bbK,Z(\widehat{G}))\to H^1(\bbL,Z(\widehat{G}))$ corresponds by the isomorphism above (applied both to $\bbK$ and $\bbL$) to the pullback of characters $G(\bbK)\backslash G(\bbA_{\bbK})\to \overline{\bbQ}_{\ell}^\times$ to $G(\bbL)\backslash G(\bbA_{\bbL})$ by the norm mapping. Thus, the above bijection restricts to a bijection between $H^1(\bbL/\bbK,Z(\widehat{G})^{\Gamma_{\bbL}})$ and the group of automorphic characters $\chi: G(\bbK)\backslash G(\bbA_{\bbK})\to \overline{\bbQ}_{\ell}^\times$ that are trivial on the image of the norm mapping. Compatibility of global Lafforgue's parameterization with twisting then yields
	\begin{equation*}
		\displaystyle I_{\hat{\chi}\rho}(f')=I_{\rho}(\chi f')
	\end{equation*}
	for every $\widehat{\chi}\in H^1(\bbL/\bbK,Z(\widehat{G})^{\Gamma_{\bbL}})$ with corresponding character $\chi: G(\bbK)\backslash G(\bbA_{\bbK})\to \overline{\bbQ}_{\ell}$. However, since $f'$ is a base-change transfer, the stable orbital integrals of $f'$ vanish outside of the image of the norm mapping and therefore $f'$ and $\chi f'$ have the same stable orbital integrals. Thus, $I_{\rho}$ being a stable distribution (by \eqref{eq Irhostable}), for every $\hat{\chi}\in H^1(\bbL/\bbK,Z(\widehat{G})^{\Gamma_{\bbL}})$ we have
	\begin{equation}\label{Irho2}
		\displaystyle I_{\hat{\chi}\rho}(f')=I_{\rho}(f'\chi)=I_{\rho}(f').
	\end{equation}
	
	It remains to compute $\iota((R_{\bbL/\bbK} G_{\bbL},\theta),\underline{G}_{\bbK})$. Again by Proposition \ref{prop ext BHKT}, and using that the $\Gamma_{\bbK}$-module $Z(\widehat{R_{\bbL/\bbK}G_{\bbL}})$ is induced from $Z(\widehat{G})$ seen as a $\Gamma_{\bbL}$-module, we have
	$$\displaystyle \ker^1(\bbK,Z(\widehat{G}))=0=\ker^1(\bbL,Z(\widehat{G}))=\ker^1(\bbK,Z(\widehat{R_{\bbL/\bbK}G_{\bbL}}))$$
	so that
	\begin{equation*}
		\displaystyle \tau(R_{\bbL/\bbK}G_{\bbL})\tau(G_{\bbK})^{-1}=\lvert Z(\widehat{G})^{\Gamma_{\bbL}}\rvert\times \lvert Z(\widehat{G})^{\Gamma_{\bbK}}\rvert^{-1}.
	\end{equation*}
	On the other hand, we readily compute that the group of outer automorphisms of the principal endoscopic group $\underline{G}_{\bbK}$ is
	\begin{equation*}
		\displaystyle \Out(\underline{G}_{\bbK})=\left(Z(\widehat{R_{\bbL/\bbK}G_{\bbL}})/Z(\widehat{G})\right)^{\Gamma_{\bbK}}.
	\end{equation*}
	By the exact sequence
	\begin{multline*}
		\displaystyle 1\to Z(\widehat{G})^{\Gamma_{\bbK}}\to Z(\widehat{G})^{\Gamma_{\bbL}}=Z(\widehat{R_{\bbL/\bbK}G_{\bbL}})^{\Gamma_{\bbK}} \\ \to \left(Z(\widehat{R_{\bbL/\bbK}G_{\bbL}})/Z(\widehat{G})\right)^{\Gamma_{\bbK}}\to H^1(\bbL/\bbK,Z(\widehat{G})^{\Gamma_{\bbL}})\to 1,
	\end{multline*}
	it follows that
	$$\displaystyle \iota((R_{\bbL/\bbK} G_{\bbL},\theta),\underline{G}_{\bbK})=\left\lvert H^1(\bbL/\bbK,Z(\widehat{G})^{\Gamma_{\bbL}})\right\rvert^{-1}.$$
	Combining this with \eqref{Irho1} and \eqref{Irho2}, we deduce $I_{\rho_{\mathbb{L}},\theta}(f)=I_{\rho}(f')$; that is the identity of Theorem \ref{hyp rho traces}.
	
	\subsection{A weakened assumption on stabilization}\label{Sect weak Stabilization}
	
	In this last subsection, we assume that $p$ does not divide the order of the absolute Weyl group of $G$.	We now want to explain  important weakenings of the Hypotheses \eqref{hyp2} and \eqref{hyp4}, pertaining to the stabilization of (twisted) trace formulas, that are actually sufficient to deduce Theorem \ref{hyp rho traces} in full generality (and therefore also sufficient as working hypotheses for the rest of the paper). More specifically, we will explain that we actually just need to assume the stabilization for test functions that are supported in the (strongly) regular semi-simple locus (or, more accurately, away from the nonregular semi-simple locus). This seems to us to be a rather important reduction since most of the extra difficulties appearing for function fields (as opposed to number fields) in the stabilization process have to do with the nonregular semi-simple contributions (including the unipotent ones) on the geometric side of the trace formula.
	
	Let us now explain more formally the Hypotheses involved. Recall that for $H$ a reductive group (defined over $\bbK$), we denote by $H_{sreg}$ the open subset of strongly regular semisimple elements. Let $H_{sing}=H\setminus H_{sreg}$ be its closed complement. Then, we say that a test function $f\in C_c^\infty(H(\bbA_{\bbK}))$ is {\em strongly regular} if it vanishes on $H_{sing}(\bbA_{\bbK})$. Let $C_c^\infty(H(\bbA_{\bbK}))_{sreg}$ be the subspace of strongly regular test functions. Concretely, we will only need to consider functions $f\in C_c^\infty(H(\bbA_{\bbK})_{sreg}$ of the following form: $f=\prod_v f_v$, where $f_v\in C_c^\infty(H(\bbK_v))$ and for at least one place $v_0$, $f_{v_0}$ is supported in $H_{sreg}(\bbK_{v_0})$. Finally, we let $SI(H(\bbA_{\bbK}))_{sreg}$ be the image of $C_c^\infty(H(\bbA_{\bbK}))_{sreg}$ in $SI(H(\bbA_{\bbK}))$.
	
	We define also, for any $L$-pair $\bbH=(H,\cH)$ over $\bbK$, a subspace $C_c^\infty(\bbH(\bbA_{\bbK}))_{sreg}$ and its image $SI(\bbH(\bbA_{\bbK}))_{sreg}\subset SI(\bbH(\bbA_{\bbK}))$. Finally, we say that a test function $f\in C_c^\infty(G(\bbA_{\bbL}))$ is {\em $\theta$-strongly regular} if it vanishes on $(R_{\bbL/\bbK}G_{\bbL})_{\theta-sing}(\bbA_{\bbK})$ where 
	$$(R_{\bbL/\bbK}G_{\bbL})_{\theta-sing}=R_{\bbL/\bbK}G_{\bbL}\setminus (R_{\bbL/\bbK}G_{\bbL})_{\theta-sreg}$$ is the complement of the (open) locus of $\theta$-strongly regular semisimple points, and we denote by $C_c^\infty(G(\bbA_{\bbL}))_{\theta-sreg}$ the subspace of $\theta$-strongly regular functions.
	
	We then make the following assumptions regarding the stabilization of (twisted) trace formulas (where this time we state everything before projection to any family of Satake parameters):
	
	\begin{hyp}\label{hyp2bis}
		For every $L$-pair $\bbH=(H,\cH)$ over $\bbK$ with $H$ quasi-split and having anisotropic center, there exists a stable distribution $S^{\bbH}_{\disc}: SI(\bbH(\bbA_{\bbK}))_{sreg}\to \overline{\bbQ}_{\ell}$ such that the following identities
		\begin{equation*}
			\displaystyle I_{\disc}^{\bbH}(f)=\sum_{\underline{\bbH}'\in \mathcal{E}_{ell}(\bbH)} \iota(\bbH, \underline{\bbH}') S^{\bbH'}_{\disc}(f^{\underline{\bbH}'})
		\end{equation*}
		are satisfied for every $L$-pair $\bbH=(H,\cH)$ as before and $f\in C_c^\infty(\bbH(\bbA_{\bbK}))_{sreg}$.
	\end{hyp}

\begin{hyp}\label{hyp4bis}
	Assume Hypothesis \ref{hyp3}. Then, we have the following identity
	\begin{equation*}
		\displaystyle I_{\disc,c}^{R_{\bbL/\bbK}G_{\bbL},\theta}(f)=\sum_{\underline{\bbH}\in \mathcal{E}_{ell}(R_{\bbL/\bbK}G_{\bbL},\theta)} \iota((R_{\bbL/\bbK}G_{\bbL},\theta), \underline{\bbH}) S^{\bbH}_{\disc}(f^{\underline{\bbH}})
	\end{equation*}
	for every $f\in C_c^\infty(G(\bbA_{\bbL}))_{\theta-sreg}$.
\end{hyp}

Note that in the context of both of these two hypotheses, the endoscopic transfers $f^{\underline{\bbH}}$ can also be chosen to be ($\theta$-)strongly regular. Now we can show that, by the same proof as the one given in Subsection \ref{proof theo rho traces}, assuming Hypotheses \eqref{hyp2bis}, \eqref{hyp4bis} (instead of Hypotheses \eqref{hyp2} and \eqref{hyp4}) as well as Hypothesis \eqref{hyp transfer}, the conclusion \eqref{iden rho traces} of Theorem \ref{hyp rho traces} holds whenever $f\in C_c^\infty(G(\bbA_{\bbL}))_{\theta-sreg}$. Therefore, all we need to see is the following:

\begin{prop}\label{prop rhotrace weak}
	Assume that $p$ does not divide the order of the absolute Weyl group of $G$, and that the identity \eqref{iden rho traces} of Theorem \ref{hyp rho traces} holds for every pair of test functions $f=\prod_v f_v$, $f'=\prod_v f'_v$ such that $f'$ is a base-change transfer of $f$ and $f'\in C_c^\infty(G(\bbA_{\bbL}))_{\theta-reg}$. Then, \eqref{iden rho traces} holds for every pair of test functions $f=\prod_v f_v$, $f'=\prod_v f'_v$ such that $f'$ is a base-change transfer of $f$.
\end{prop}

\begin{proof}
First, let us pick a place $v_0$ of $\bbK$ such that:
\begin{itemize}
	\item The group $G$ is unramified at $v_0$.  In particular we can fix a Borel subgroup $B=TU\subset G_{v_0}$;
	
	\item The residue field of $\bbK_{v_0}$ is of cardinality $q>5$;
	
	\item $\rho$ is unramified at $v_0$;
	
	\item The centralizer $Z_{\widehat{G}}(\rho(\Frob_{v_0}))$ is connected modulo $Z(\widehat{G})^{\Frob_{v_0}}$ and, for every place $w_0\mid v_0$ of $\bbL$, the centralizer $Z_{\widehat{G}}(\rho(\Frob_{w_0}))$ is connected modulo $Z(\widehat{G})^{\Frob_{w_0}}$. 
\end{itemize}
For any splitting field $\bbK'/\bbK$ of $G$, the subset of elements $\gamma\in \widehat{G}\rtimes \Gal(\bbK'/\bbK)$ whose centralizer is contained in a maximal torus of $\widehat{G}$ is open (and nonempty) in the Zariski topology. Moreover, for every such element $\gamma$, the centralizer $Z_{\widehat{G}}(\gamma)$ is connected modulo its intersection $Z_{\widehat{G}}(\gamma)\cap Z(\widehat{G})$ with the center of $\widehat{G}$. Thus, it follows from the fact that $\rho$ has Zariski dense image and the Cebotarev density theorem that there exist infinitely many places $v_0$ satisfying the four conditions above.

Moreover, since $\rho$ has Zariski dense image,  in particular it is irreducible (that is its image is not contained in any proper parabolic subgroup).  Then by \cite[Lemma 11.4]{ST21} the elements $\rho(\Frob_{v_0})$ and $\rho(\Frob_{w_0})$, $w_0\mid v_0$, of ${}^L G$ are {\em pure}, that is their eigenvalues in the adjoint representation of ${}^L G$ are all Weil numbers of weight $0$. In particular, using our fixed isomorphism $\iota: \overline{\bbQ_{\ell}}\simeq \bbC$, $\rho(\Frob_{v_0})$ and $\rho(\Frob_{w_0})$ are the Satake parameters of unitary unramified principal series $I(\chi')$ and $I(\chi_{w_0})$ of $G(\bbK_{v_0})$ and $G(\bbL_{w_0})$ respectively. Moreover, by the spectral characterization of $R$-groups for principal series (see \cite[Section 5]{Sol25}) and the condition imposed on the centralizers of $\rho(\Frob_{v_0})$, $\rho(\Frob_{w_0})$, we know that these principal series are irreducible. Set
$$\displaystyle I(\chi)=\bigotimes_{w_0\mid v_0} I(\chi_{w_0}),$$
a unitary irreducible principal series of $G(\bbL_{v_0})$.

Let $f^{v_0}=\prod_{v\neq v_0} f_v\in C_c^\infty(G(\bbA_{\bbL}^{v_0}))$ and $f^{',v_0}=\prod_{v\neq v_0} f'_v\in C_c^\infty(G(\bbA_{\bbK}^{v_0}))$ be arbitrary test functions such that, for each $v\neq v_0$, $f'_v$ is a base-change transfer of $f_v$. Then, we need to prove that for every $f_{v_0}\in C_c^\infty(G(\bbL_{v_0}))$, $f'_{v_0}\in  C_c^\infty(G(\bbK_{v_0}))$, where $f'_{v_0}$ is a base-change transfer of $f_{v_0}$, we have
\begin{equation}\label{eq1 rhotrace weak}
\displaystyle I_{\rho_{\mathbb{L}},\theta}(f^{v_0}f_{v_0})=I_\rho(f^{',v_0}f'_{v_0}).
\end{equation}
By assumption, this identity is verified every time $f_{v_0}$ is supported in
$$G(\bbL_{v_0})_{\theta-sreg}:=(R_{\bbL/\bbK} G_{\bbL})_{\theta-sreg}(\bbK_{v_0})\subset G(\bbL_{v_0}).$$
We will consider both sides as distributions in $f_{v_0}$ and $f'_{v_0}$ respectively. The crucial point is now the following fact:
\begin{num}
\item\label{eq2 rhotrace weak} For every automorphic representation $\pi'\subset \CA_{0,\rho}(G)$ (resp. $\pi\subset \CA_{0,\rho_{\mathbb{L}}}(G_{\mathbb{L}})$), its local component at $v_0$ is $\pi'_{v_0}\simeq I(\chi')$ (resp. $\pi_{v_0}\simeq I(\chi)$).
\end{num}
Indeed, by the local-global compatibility of \cite{GLa}, we know that the Genestier-Lafforgue parameter of $\pi'_{v_0}$ (resp. $\pi_{w_0}$ for $w_0\mid v_0$) is $\rho\mid_{W_{\bbK_{v_0}}}: W_{\bbK_{v_0}}\to {}^L G$ (resp. $\rho\mid_{W_{\bbL_{w_0}}}: W_{\bbL_{w_0}}\to {}^L G$). Since this local parameter is pure and unramified, by \cite[Theorem 1.4(i)]{GHS}, the unique irreducible representation of $G(\bbK_{v_0})$ with Genestier-Lafforgue parameter $\rho\mid_{W_{\bbK_{v_0}}}$ (resp. $\rho\mid_{W_{\bbL_{w_0}}}$) is $I(\chi')$ (resp. $I(\chi_{w_0})$). Here we note that \cite[Theorem 1.4(i)]{GHS} applies in the case at hand to any supercuspidal representation of a Levi subgroup of $G(\bbK_{v_0})$ or $G(\bbL_{w_0})$ since we are assuming that $q>5$ and $p\nmid \lvert W(G)\rvert$, which by a result of Fintzen \cite{Fin} implies that every such supercuspidal representation is compactly induced.

By \eqref{eq2 rhotrace weak}, the local distributions
\begin{equation*}
\displaystyle f_{v_0}\mapsto I_{\rho_{\mathbb{L}},\theta}(f^{v_0}f_{v_0}) \mbox{ and } f'_{v_0}\mapsto I_\rho(f^{',v_0}f'_{v_0})
\end{equation*}
are proportional to
\begin{equation*}
\displaystyle f_{v_0}\mapsto \Tr(I(\chi,f_{v_0})\circ I(\chi,\theta)) \mbox{ and } f'_{v_0}\mapsto \Tr(I(\chi',f'_{v_0}))
\end{equation*}
respectively. Here, we have denoted by $I(\chi,\theta)$ the natural extension of the automorphism $\theta$ to the principal series $I(\chi)$ (note that, by construction, $\chi\circ \theta=\chi$). Thus, we are reduced to proving the following: if $C_1,C_2\in \bbC$ are such that
$$\displaystyle C_1 \Tr(I(\chi,f_{v_0})\circ I(\chi,\theta))=C_2 \Tr(I(\chi',f'_{v_0}))$$
every time $f'_{v_0}\in C_c^\infty(G(\bbK_{v_0}))$ is a base-change transfer of $f_{v_0}\in C_c^\infty(G(\bbL_{v_0}))$ with $\Supp(f_{v_0})\subset G(\bbL_{v_0})_{\theta-sreg}$, then the identity holds without this last assumption. That fact is clear once we know that both distributions are represented by locally integrable functions which are moreover stable (i.e. constant on strongly regular ($\theta$-)stable conjugacy classes) by the characterization of base-change transfer given in Proposition \ref{Prop char bc GIT}. That (twisted) characters of principal series are given by locally integrable stable functions is well-known, see the proof of Proposition \ref{ptw}.
\end{proof}
	
	\section{Construction of types}\label{types}
	
	\subsection{Statement of the main result}
	
	In this section, we let $F$ be a local field of residual characteristic $p$ and $G$ a connected semisimple group over $F$, that is unramified (in particular, quasi-split). Let $B=TU$ be a Borel subgroup of $G$, with $U$ its unipotent radical and $T$ a maximal torus. We fix a non-degenerate character $\psi^F: U(F) \rightarrow \mathbb{C}^{\times}$. Let $U_{der}$ be the derived subgroup of $U$. Then, $U/U_{der}$ is a product of additive groups and, for every finite separable extension $K/F$, we let $\psi^K: U(K)\to \mathbb{C}^\times$ be the composition of $\psi^F$ with the trace map $\Tr_{K/F}: U(K)/U_{der}(K)=(U/U_{der})(K)\to (U/U_{der})(F)=U(F)/U_{der}(F)$.
	
	For every integer $n\geq 0$, we let $T(F)_n$ be the $n$-th congruence subgroup of $T(F)$, that is
	\begin{equation}
		\displaystyle T(F)_n=\Ker(\bbT(\mathcal{O}_F)\to \bbT(\mathcal{O}_F/\mathfrak{m}_F^n)),
	\end{equation}
	where $\bbT$ denotes the unique $\mathcal{O}_F$-torus with generic fiber $T$ (which exists since $T$ is unramified).

	Recall that $p$ is said to be {\em very good} for $G$ if it does not divide the coefficients of the decomposition of the highest root in simple roots nor the order of the center of its simply-connected cover. In particular, if $p$ does not divide the order of the Weyl group of $G$, $p$ is very good for $G$.

	\begin{theorem}\label{MWKV} 
		Assume that $p$ is very good for $G$ and let $E/F$ be a finite Galois tamely ramified extension. Then, there exist sequences $\left(\cP_{n}^K, \chi_{n}^K\right)_{n\geq 1}$, indexed by subextensions $F\subset K\subset E$, with the following properties:
		\begin{enumerate}
			\item For every $n$ and $K$, $\cP_{n}^K$ is a compact-open subgroup of $G(K)$ and $\chi_{n}^K: \cP_{n}^K \to\mathbb{C}^{\times}$ is a character. We denote by $f_{n}^K \in C_c^{\infty}(G(K))$ the function supported in $\cP_{n}^K$ and equal to $\operatorname{vol}\left(\cP_{n}^K\right)^{-1} (\chi_{n}^K)^{-1}$ there.
			\item For any cyclic subextension $L / K$ of $E/F$ with $\operatorname{Gal}(L / K)=\langle \theta\rangle$ and $[L:K]$ prime to $p$ and any $n\geq 1$, $f_{n}^K$ is a base-change transfer of $f_{n}^L$; that is
			$$\displaystyle
			SO_\gamma\left(f_{n}^K\right)= \begin{cases}0 & \text { if } \gamma \text { is not a norm, } \\ SO_{\delta, \theta}\left(f_{n}^L\right) & \text { if } \gamma\sim N(\delta),\end{cases}
			$$
			for every strongly regular element $\gamma \in G(K)$. (See Subsection \ref{subsection rhotrace} for the relevant notation.)
        \end{enumerate}
Moreover, there exists an integer $n_0\geq 1$, that only depends on the Dynkin type of $G$ (that is, on the absolute root system of $G$)\footnote{This bound can be made explicit, e.g. we can take $n_0=4$ times the sum of the coefficients of the highest root.}, such that:
        \begin{enumerate}
        \setcounter{enumi}{2}
			\item  For every $F\subset K\subset E$ and every smooth irreducible representation $\pi$ of $G(K)$, the sequence $(\operatorname{dim} \pi^{\cP_{n}^K, \chi_{n}^K})_{n\geq n_0}$ is increasing and for $n$ large enough, we have
			\begin{equation}\label{eq KtypeandWhitt}
				\operatorname{dim} \pi^{\cP_{n}^K, \chi_{n}^K}=\operatorname{dim} \Hom_{U(K)}(\pi,\psi^K)
			\end{equation}
			i.e. $\operatorname{dim} \pi^{\cP_{n}^K, \chi_{n}^K}=1$ for $n$ large enough if $\pi$ is generic with respect to $\psi^K$ and 0 otherwise.
			
			\item For every $n\geq n_0$ and every character $\chi: T(K)/T(K)_n\to \mathbb{C}^\times$, we have
			\begin{equation*}
				\displaystyle \dim I_{B(K)}^{G(K)}(\chi)^{\cP_n^K,\chi_n^K}=1,
			\end{equation*}
			where $I_{B(K)}^{G(K)}(\chi)$ denotes the (normalized) parabolic induction of $\chi$.
		\end{enumerate}
	\end{theorem}
	
	In the next subsection, we explain a general mechanism for the construction of sequences of types $\left(\cP_{n}^L, \chi_{n}^L\right)_{n\geq 1}$ as in the theorem above (that is pairs satisfying point $(i)$ above) using Moy-Prasad filtrations. In Subsection \ref{ssection proof of ii}, we show that these sequences always satisfy property $(ii)$ of the theorem. Finally, in Subsection \ref{ssection proof iii} we specialize this construction to a special case and show that it verifies properties $(iii)$ and $(iv)$ under the assumption that $p$ is very good for $G$. 
	
	Let us emphasize that the above theorem is inspired by work of Rodier \cite{Rod74} later generalized by M\oe{}glin-Waldspurger \cite{MW87}. Indeed, in {\em loc. cit.} the authors construct a sequence of $K$-types $\left(\cP_{n}^F, \chi_{n}^F\right)_{n\geq 1}$ of $G(F)$ when $F$ is a local non-Archimedean field of {\em characteristic zero} satisfying the property \eqref{eq KtypeandWhitt} (for $n$ large enough). The constructions given in these references use in a crucial way the exponential map, which cannot be applied directly here and is in a sense the reason for our appeal to Moy-Prasad filtrations; the Moy-Prasad isomorphisms serve as an adequate substitute for the exponential mapping.
	
	Here is some notation that will be useful in the subsequent subsections:
	\begin{itemize}
		\item We let $A$ be the maximal split subtorus of $T$ and $\overline{B}=T\overline{U}$ the Borel subgroup opposite to $B$ (with respect to $T$). 
		
		\item Lie algebras will be denoted by the corresponding small gothic letter e.g. we write $\mathfrak{g}$, $\mathfrak{t}$, $\mathfrak{a}$, $\mathfrak{u}$ for the Lie algebras of $G$, $T$, $A$ and $U$ respectively.
		
		\item We denote by $\Phi$  the set of roots of $A$ in $G$, $\Phi^+\subset \Phi$  the subset of positive roots with respect to $B$, $\Phi^-$ the subset of negative roots and $\Delta\subset \Phi^+$ the corresponding simple roots.
		
		\item For every $\alpha\in \Phi$, we denote by $\mathfrak{u}_\alpha\subset \mathfrak{g}$ the corresponding root subspace. In particular, we have the decomposition $\mathfrak{u}=\bigoplus_{\alpha\in \Phi^+} \mathfrak{u}_\alpha$.
		
		\item We write $\mathcal{O}_F$ for the ring of integers in $F$, $\mathfrak{m}_F\subset \mathcal{O}_F$ for its maximal ideal, $\varpi_F\in \mathfrak{m}_F$ for a uniformizer and $k_F=\mathcal{O}_F/\mathfrak{m}_F$ for its residue field. Similar notation will be applied to finite extensions of $F$.

	\end{itemize}
	
	\subsection{Review of Moy-Prasad filtrations}
	
	Let $\mathcal{B}(G,F)$ be the Bruhat-Tits building of $G/F$. To every $x\in \mathcal{B}(G,F)$, we denote by $\mathbb{G}^F_x$ the affine smooth model of $G$ over $\mathcal{O}_F$ such that $G(F)_{x,0}:=\mathbb{G}^F_x(\mathcal{O}_F)$ is the parahoric subgroup corresponding to the facet containing $x$. This admits a filtration by the compact-open subgroups
	$$\displaystyle G(F)_{x,n}:= \mathbb{G}_x^F(\mathfrak{m}_F^n):=\Ker(\mathbb{G}_x^F(\mathcal{O}_F)\to \mathbb{G}_x^F(\mathcal{O}_F/\mathfrak{m}_F^n)),\;\; n\in \bbN.$$
	(More generally, for every group scheme $X$ over $\mathcal{O}_F$, we will henceforth write $X(\mathfrak{m}_F^n)$ for the kernel of the morphism $X(\mathcal{O}_F)\to X(\mathcal{O}_F/\mathfrak{m}_F^n)$.) The sequence $(G(F)_{x,n})_{n\geq 0}$ satisfies the following basic property \cite[Proposition 13.2.5 (4)]{KP}:
	\begin{num}
		\item\label{eq1 MPfiltration} for every $n,m\geq 0$, $[G(F)_{x,n},G(F)_{x,m}]\subseteq G(F)_{x,n+m}$ where $[.,.]$ denotes the commutator map $[g,h]=ghg^{-1}h^{-1}$.	
	\end{num}

	Let $\mathcal{A}\subset \mathcal{B}(G,F)$ be the apartment corresponding to $T$. Assume that $x\in \mathcal{A}$. Let $\mathbb{T}$, $\mathbb{A}$ be the schematic closures of $T$, $A$ respectively in $\mathbb{G}^F_x$. These are $\mathcal{O}_F$-tori, $\mathbb{A}$ being a split $\mathcal{O}_F$-torus, that do not depend on the choice of $x$. Let $\mathbb{U}_{x}^F$, $\overline{\mathbb{U}}_x^F$ be the schematic closures of $U$ and $\overline{U}$ in $\mathbb{G}_x^F$. These are smooth $\mathcal{O}_F$-models of $U$ and $\overline{U}$ respectively (this follows e.g. from \cite[Proposition 8.3.14 (1)-(2)]{KP}). Set, for every $n\geq 0$,
	$$\displaystyle T(F)_{n}=\mathbb{T}(\mathfrak{m}_F^n),\; U(F)_{x,n}=\mathbb{U}_{x}^F(\mathfrak{m}_F^n),\; \overline{U}(F)_{x,n}=\overline{\mathbb{U}}_x^F(\mathfrak{m}_F^n).$$
	Then, for every $n\geq 1$ the group law induces a homeomorphism (see \cite[Theorem 13.2.5(3)]{KP})
	\begin{equation}\label{eq2 MPfiltration}
		\displaystyle  \overline{U}(F)_{x,n} \times T(F)_{n}\times U(F)_{x,n}\simeq G(F)_{x,n}.
	\end{equation}
	
	Let $\mathfrak{G}_x^F:=\Lie(\mathbb{G}^F_x)$, $\mathfrak{T}:=\Lie(\mathbb{T})$, $\mathfrak{U}_x^F=\Lie(\mathbb{U}_x^F)$ be the Lie algebras of the group schemes $\mathbb{G}^F_x$, $\mathbb{T}$ and define, for $n\in \bbN$,
	$$\displaystyle \mathfrak{g}(F)_{x,n}:=\mathfrak{G}_x^F(\mathfrak{m}_F^n)=\varpi_F^n \mathfrak{g}(F)_{x,0},\;\; \mathfrak{t}(F)_{n}:=\mathfrak{T}(\mathfrak{m}_F^n),\;\; \mathfrak{u}(F)_{x,n}:=\mathfrak{U}_x^F(\mathfrak{m}_F^n).$$
	For every $\alpha$, let $\mathfrak{U}_{\alpha,x}^F$ be the schematic closure of $\mathfrak{u}_\alpha$ in $\mathfrak{G}_x^F$. This is also the $\alpha$-eigenspace for the adjoint action of $\mathbb{A}$ on $\mathfrak{G}_x^F$. Set
	$$\displaystyle \mathfrak{u}_\alpha(F)_{x,n}:=\mathfrak{U}_{\alpha,x}^F(\mathfrak{m}_F^n),\;\; \alpha\in \Phi,n\geq 0.$$
	Then, we have
	\begin{equation}\label{eq3 MPfiltration}
		\displaystyle \mathfrak{g}(F)_{x,n}=\mathfrak{t}(F)_{n}\oplus \bigoplus_{\alpha\in \Phi} \mathfrak{u}_\alpha(F)_{x,n}.
	\end{equation}
	
	For every $n\geq 1$, there exists a natural isomorphism, called the {\em Moy-Prasad isomorphism} \cite[Theorem 13.5.1]{KP}
	$$\displaystyle MP_{x,n}^F: G(F)_{x,n}/G(F)_{x,2n}\simeq \mathfrak{g}(F)_{x,n}/\mathfrak{g}(F)_{x,2n},$$
	which verifies the following properties:
	\begin{num}
		\item\label{eq4 MPfiltration} $MP_{x,n}^F$ sends $T(F)_{n}$ to $\mathfrak{t}(F)_{n}/\mathfrak{t}(F)_{2n}$ and $U(F)_{x,n}$ to $\mathfrak{u}(F)_{x,n}/\mathfrak{u}(F)_{x,2n}$ for every $\alpha\in \Phi$ (this follows from the functoriality of Moy-Prasad isomorphisms \cite[Proposition A.5.19(3)]{KP} and the fact that the group schemes $\mathbb{T}$, $\mathbb{U}_{x}^F$ are smooth);
	\end{num}
	
	\begin{num}
		\item\label{eq5 MPfiltration} It sends commutators to Lie brackets: for $n,m\geq 1$, $g\in G(F)_{x,n}$, $h\in G(F)_{x,m}$, we have
		$$MP_{x,n+m}^F([g,h])=[MP_{x,n}^F(g),MP_{x,m}^F(h)] \mbox{ modulo } \mathfrak{g}(F)_{x,n+m+\min(m,n)}.$$
		(see \cite[Proposition 1.6.2]{Adl98}.)
	\end{num}

	Let now $K/F$ be a tamely ramified finite Galois extension with Galois group $\Gamma=\Gal(K/F)$. Then, we have a natural embedding $\mathcal{B}(G,F)\subset \mathcal{B}(G,K)$. In particular, for every $x\in \mathcal{B}(G,F)$, we may consider the Moy-Prasad subgroups $G(K)_{x,n}\subset G(K)$ ($n\in \bbN$). For our purpose, it will be more convenient to index those with respect to $F$ rather than $K$. More precisely, if $\mathbb{G}_x^K$ denotes the parahoric integral model of $G$ over $\mathcal{O}_K$ associated to $x$, then we set
	$$\displaystyle G(K)_{x,n}:=\Ker(\mathbb{G}_x^K(\mathcal{O}_K)\to \mathbb{G}_x^K(\mathcal{O}_K/\mathfrak{m}_F^n \mathcal{O}_K)).$$
	This corresponds to rescaling the depth $n$ by the ramification index of $K/F$. We apply a similar convention to the Lie subalgebras $\mathfrak{g}(K)_{x,n}=\varpi_F^n \mathfrak{g}(K)_{x,0}$.
	
	Let $x\in \mathcal{B}(G,F)$. We will need the following properties relating the Moy-Prasad filtrations over $F$ and $K$: for each $n\in \bbN$,
	\begin{num}
		\item\label{eq7 MPfiltration} $G(K)_{x,n}$, $\mathfrak{g}(K)_{x,n}$ are $\Gamma$-stable and $G(K)_{x,n}^\Gamma=G(F)_{x,n}$, $\mathfrak{g}(K)_{x,n}^\Gamma=\mathfrak{g}(F)_{x,n}$ \cite[Proposition 12.9.4]{KP};
		
		\item\label{eq8 MPfiltration} $(G(K)_{x,n}/G(K)_{x,2n})^\Gamma=G(F)_{x,n}/G(F)_{x,2n}$ and $(\mathfrak{g}(K)_{x,n}/\mathfrak{g}(K)_{x,2n})^\Gamma=\mathfrak{g}(F)_{x,n}/\mathfrak{g}(F)_{x,2n}$ \cite[Theorem 13.8.5]{KP};
		
		\item\label{eq9 MPfiltration} the isomorphism $MP_{x,n}^K$ is $\Gamma$-equivariant (this follows from \cite[Theorem 13.5.1]{KP});
		
		\item\label{eq10 MPfiltration} $\Tr_{K/F} \mathfrak{g}(K)_{x,n}=\mathfrak{g}(F)_{x,n}$ (where $\Tr_{K/F}: \mathfrak{g}(K)\to \mathfrak{g}(F)$ denotes the trace map). Indeed, by \eqref{eq7 MPfiltration} we have $\Tr_{K/F} \mathfrak{g}(K)_{x,n}\subset \mathfrak{g}(K)^\Gamma_{x,n}=\mathfrak{g}(F)_{x,n}$ and $\mathcal{O}_K\otimes_{\mathcal{O}_F} \mathfrak{g}(F)_{x,n}\subset \mathfrak{g}(K)_{x,n}$, so that, since $\Tr_{K/F}(\mathcal{O}_K)=\mathcal{O}_F$, the trace map $\Tr_{K/F}: \mathfrak{g}(K)_{x,n}\to \mathfrak{g}(F)_{x,n}$ is surjective.
	\end{num}
	
	\subsection{Proof of property (ii)}\label{ssection proof of ii}
	
	Let $x\in \mathcal{A}\subset \mathcal{B}(G,F)$ and denote, for every tamely ramified extension $K/F$, by $(\cP^K_n)_{n\geq 1}$ the sequence of Moy-Prasad subgroups $(G(K)_{x,n})_{n\geq 1}$. We fix, for every $n\geq 1$, a character $\chi_n^F: \mathfrak{g}(F)_{x,n}/\mathfrak{g}(F)_{x,2n}\to \bbC^\times$ that we identify, by composition with $MP^F_{x,n}$, with a character of $\cP^F_n$ trivial on $\cP^F_{2n}$. More generally, for every tamely ramified extension $K/F$, we denote by $\chi_n^K$ the composition of $\chi_n^F$ with the trace map $\Tr_{K/F}: \mathfrak{g}(K)_{x,n}\to \mathfrak{g}(F)_{x,n}$ that we identify, again by means of $MP^K_{x,n}$, with a character of $\cP^K_n/\cP^K_{2n}$.
	
	As in Proposition \ref{MWKV}, we let $f_{n}^K$ be the function on $G(K)$ that is equal to $vol(\cP_{n}^K)^{-1}(\chi_{n}^K)^{-1}$ on $\cP_{n}^K$ and to zero outside of $\cP_{n}^K$. 
	
	The proposition below shows that point (ii) of Proposition \ref{MWKV} is satisfied for the sequences of $K$-types $(\cP_n^K,\chi_n^K)_{n\geq 1}$:
	
	\begin{prop}\label{prop bc transfer types}
		Let $F\subset K\subset L$ be finite extensions with $L/K$ cyclic of degree prime to $p$. Let $\theta$ be a generator of $\Gal(L/K)$. Then, for every $n\geq 1$, $f_{n}^K$ is a base-change transfer of $f_{n}^L$ i.e. for every strongly regular element $\gamma\in G(K)$, we have
		$$\displaystyle
		SO_\gamma\left(f_{n}^K\right)= \begin{cases}0 & \text { if } \gamma \text { is not a norm, } \\ SO_{\delta, \theta}\left(f_{n}^L\right) & \text { if } \gamma\sim N(\delta),\end{cases}
		$$
		where $SO_{\delta,\theta}$ and $SO_\gamma$ stand for stable orbital integrals normalized as in Subsection \ref{subsection rhotrace}.
	\end{prop}
	
	\begin{proof}
		This is essentially a consequence of the following lemma.
		
		\begin{lemma}
			The conjugation action induces a homeomorphism
			\begin{align}\label{slice isom}
				\displaystyle \cP_{n}^K\theta\times^{\cP_{n}^K} \cP_{n}^L\simeq \cP_{n}^L\theta, \\
				\nonumber (h\theta,g)\mapsto g^{-1}h\theta g,
			\end{align}
			where $\cP_{n}^K\theta\times^{\cP_{n}^K} \cP_{n}^L$ denotes the quotient of $\cP_{n}^K\theta\times \cP_{n}^L$ by the (free) $\cP_{n}^K$-action given by $h\cdot (h'\theta,g)=(hh' \theta h^{-1},hg)$.
			
			Moreover, the isomorphism \eqref{slice isom} sends the quotient of the normalized\footnote{By which we mean the measure of total mass one.} invariant measure on $\cP_{n}^K\theta\times \cP_{n}^L$ by the normalized Haar measure on $\cP_n^K$ to the normalized Haar measure on $\cP_{n}^L\theta$.
		\end{lemma}
		
		\begin{proof}
			Since both sides are compact, to show that \eqref{slice isom} is an isomorphism it suffices to show that it is a bijection that is:
			\begin{itemize}
				\item Every element of $\cP_{n}^L\theta$ is $\cP_{n}^L$-conjugate to an element of $\cP_{n}^K\theta$;
				
				\item For every $k\in \cP_{n}^L$ and $k'\in \cP_{n}^K$ if $kk'\theta k^{-1}\in \cP_{n}^K\theta$ then $k\in \cP_{n}^K$.
			\end{itemize}
			The proofs of both points are similar. We start with the first point. By induction and compactness, it suffices to show that for every $\ell\geq 1$, every element of $\cP_{n}^K\cP_{n\ell}^L\theta$ is $\cP_{n}^L$-conjugate to an element in $\cP_{n}^K\cP_{n(\ell+1)}^L\theta$. Let $k\in \cP_{n}^K$ and $k'\in \cP_{n\ell}^L$. Because the degree of $L/K$ is prime to $p$, and by \eqref{eq7 MPfiltration}, we have the decomposition
			\begin{equation}\label{eq1 prop bc transfer types}
				\displaystyle \mathfrak{g}(L)_{x,\ell n}=\mathfrak{g}(L)_{x,\ell n}^\theta\oplus (1-\theta) \mathfrak{g}(L)_{x,\ell n}=\mathfrak{g}(K)_{x,\ell n}\oplus (1-\theta) \mathfrak{g}(L)_{x,\ell n}.
			\end{equation}
			From there and \eqref{eq8 MPfiltration}, and \eqref{eq9 MPfiltration}, we deduce the existence of $j\in \cP_{\ell n}^L$ such that $jk'\theta j^{-1}\in \cP_{\ell n}^K \cP_{(\ell+1)n}^L\theta$. Since $jk=kj$ modulo $\cP_{(\ell+1)n}^L$ (by \eqref{eq1 MPfiltration}), this implies $jkk'\theta j^{-1}\in \cP_{n}^K\cP_{n(\ell+1)}^L\theta$ and this shows the first point.
			
			The second point can be dealt with in a similar manner. Namely, it suffices to show that for $k\in \cP_{n}^L$ and $k'\in \cP_{n}^K$ such that $kk'\theta k^{-1}\in \cP_{n}^K\theta$ we have $k\in \cP_{\ell n}^L \cP_{n}^K$ for every $\ell\geq 1$. For this, we proceed by induction: up to translating $k$ and conjugating $k'\theta$ by an element of $\cP_n^K$, we may assume that $k\in \cP_{\ell n}^L$ for some $\ell\geq 1$ and we aim to deduce that $k\in \cP_{\ell n}^K \cP_{(\ell+1)n}^L$. Since $kk'=k'k \mod \cP_{(\ell+1)n}^L$ (still by \eqref{eq1 MPfiltration}), the assumption implies $k \theta(k)^{-1}\in \cP_{\ell n}^K \cP_{(\ell+1)n}^L$. Hence, by \eqref{eq1 prop bc transfer types}, $\theta(k\cP_{(\ell+1)n}^L)=k\cP_{(\ell+1)n}^L$ that is, by \eqref{eq8 MPfiltration}, $k\in \cP_{\ell n}^K \cP_{(\ell+1)n}^L$ as claimed.

			We now turn to the claim regarding compatibility of the isomorphism with normalized Haar measures. The map \eqref{slice isom} is clearly locally analytic and it therefore suffices to show that its Jacobian is equal to $1$ at each point. Since it is equivariant for the obvious $\cP_n^L$-actions on both sides, we just need to compute this Jacobian at a point $(h\theta,1)\in \cP_{n}^K\theta\times \cP_{n}^L$. We can identify the tangent space at this point with the quotient $(\mathfrak{g}(K)\oplus \mathfrak{g}(L))/\mathfrak{g}(K)$ where $\mathfrak{g}(K)$ is embedded in the direct sum $\mathfrak{g}(K)\oplus \mathfrak{g}(L)$ by means of the map $Z\mapsto (Z-\Ad(h\theta)(Z),Z)$. The differential of \eqref{slice isom} at $(h\theta,1)$ is then given by
			\begin{equation}\label{eq2 slice isom}
				\displaystyle (\mathfrak{g}(K)\oplus \mathfrak{g}(L))/\mathfrak{g}(K)\to \mathfrak{g}(L),
			\end{equation}
			$$\displaystyle (X,Y)\mapsto X+\Ad(h\theta)(Y)-Y.$$
			What we need to check is that this differential is an isomorphism preserving the measures on both sides. More precisely, the normalized Haar measures on $\cP_n^K$ and $\cP_n^L$ induce Haar measures on $\mathfrak{g}(K)$ and $\mathfrak{g}(L)$ respectively hence on the quotient $(\mathfrak{g}(K)\oplus \mathfrak{g}(L))/\mathfrak{g}(K)$, and we aim to show that it corresponds to the measure on $\mathfrak{g}(L)$ by the linear map \eqref{eq2 slice isom}. For this, it obviously suffices to show the existence of $\mathcal{O}_K-$ and $\mathcal{O}_L-$lattices $\Lambda^K\subset \mathfrak{g}(K)$, $\Lambda^L\subset \mathfrak{g}(L)$ such that \eqref{eq2 slice isom} reduces to an isomorphism $(\Lambda^K\oplus \Lambda^L)/\Lambda^K\simeq \Lambda^L$ i.e. such that we have a short exact sequence
			$$\displaystyle 0\to \Lambda^K\to \Lambda^K\oplus \Lambda^L\to \Lambda^L\to 0$$
			where the first map is $Z\mapsto (Z-\Ad(h\theta)(Z),Z)$ and the second one $(X,Y)\mapsto X+\Ad(h\theta)(Y)-Y$. We claim that this works for $\Lambda^K=\mathfrak{g}(K)_{x,n}$ et $\Lambda^L=\mathfrak{g}(L)_{x,n}$. That it is indeed the case is just a variation on the two above bullet points at the Lie algebra level and is left to the reader.
		\end{proof}
		
		Consider now the following diagram
		$$\displaystyle \xymatrix{G(K)\times G(L) \ar[r]^{\;\;\;\;\;\;\; p} \ar[d]_N & G(L) \ar[d]^{\pi_\theta} \\ G(K) \ar[r]_{\pi} & \mathcal{Q}(K)}$$
		where $p$ is given by the conjugation action $(\gamma,g)\mapsto g^{-1}\gamma \theta(g)$, $N$ is the first projection followed by the norm map $\gamma\mapsto \gamma^d$ ($d=[L:M]$), $\mathcal{Q}$ denotes the GIT quotient $G\sslash \Ad(G)$ for the conjugation action that we identify (by mean of the norm mapping) with the GIT quotient $R_{L/K}G_L\sslash \Ad_\theta(R_{L/K} G_L)$ for the action by $\theta$-conjugation of $R_{L/K}G_L$ on itself and $\pi_\theta$, $\pi$ stand for the natural projections. Then, the statement of the proposition is equivalent to the equality
		\begin{equation}\label{eq3}
			\displaystyle \pi_{*}f_{n}^K=\pi_{\theta,*}f_{n}^L
		\end{equation}
		between the pushforward to the regular locus in $\mathcal{Q}(K)$ of the products of $f_{n}^K$ and $f_{n}^L$ with Haar measures. Let $g_n$ be the product of the normalized Haar measure on $\cP_{n}^K\times \cP_{n}^L$ with the restriction of the character $\chi_{n}^L$ to the first factor. Note that, by construction, $\chi_{n}^L\mid \cP_{n}^K=(\chi_{n}^K)^d=\chi_{n}^K\circ N$. Then, we have $N_* g_n=f_{n}^K$ whereas, by the lemma above, $p_*g_n=f_{n}^L$. Hence, \eqref{eq3} simply follows from the commutativity of the diagram.
	\end{proof}
	
	\subsection{Approximation of Whittaker coinvariants and properties (iii)-(iv)}\label{ssection proof iii}
	
	Our goal is now to show how we can choose the point $x\in \mathcal{A}$ and the characters $\chi_{n}^F$ such that point $(iii)$ and $(iv)$ of Theorem \ref{MWKV} are satisfied. 
	
	Let $\overline{B}=T\overline{U}$ be the Borel subgroup opposite to $B$ with respect to $T$. Let $\overline{\mathfrak{b}}$ and $\mathfrak{u}_{der}$ be the Lie algebras of $\overline{B}$ and $U_{der}$. Then, we have a natural isomorphism $U/U_{der}\simeq \mathfrak{u}/\mathfrak{u}_{der}$ so that we can identify $\psi^F$ with a character of $\mathfrak{u}(F)$ trivial on $\mathfrak{u}_{der}(F)$; the hypothesis that $\psi^F$ is non-degenerate translates to: for every $\alpha\in \Delta$, the restriction of $\psi^F$ to $\mathfrak{u}_{\alpha}(F)$ is non-trivial. We further extend $\psi^F$ to a character of $\mathfrak{g}(F)$ that is trivial on  $\overline{\mathfrak{b}}(F)$, through the decomposition
	$$\displaystyle \mathfrak{g}=\overline{\mathfrak{b}}\oplus \mathfrak{u}.$$
	For every subextension $F\subset K\subset E$, we set $\psi^K:=\psi^F\circ \Tr_{K/F}: \mathfrak{g}(K)\to \bbC^\times$.
	
	We choose $x\in \mathcal{A}$ such that:
	\begin{num}
		\item $x$ is hyperspecial as a point of the building $\mathcal{B}(G,E)$ of $G/E$;
		\item\label{cond2 x} For every $\alpha\in \Delta$, $\psi^E$ is trivial on $\mathfrak{u}_\alpha(E)_{x,0}$ but nontrivial on $\varpi_E^{-1}\mathfrak{u}_\alpha(E)_{x,0}$. In particular, $\psi^E$ is trivial on $\mathfrak{g}(E)_{x,0}$.
	\end{num}
	This is always possible. Indeed, let $x_0\in \mathcal{A}\subset \mathcal{B}(G,F)$ be a hyperspecial point, which exists since $G$ is an unramified group, and $A_{ad}$ be the image of $A$ in the adjoint quotient $G_{ad}$ of $G$. For every $\alpha\in \Delta$, let $n_\alpha$ be the largest integer such that $\psi^E$ is trivial on $\varpi_E^{-n_\alpha}\mathfrak{u}_\alpha(E)_{x_0,0}$. Let $a\in A_{ad}(E)$ be such that $\alpha(a)=\varpi_E^{n_\alpha}$ for every $\alpha\in \Delta$. Then, the point $x=a\cdot x_0\in \mathcal{A}$ works. Indeed, $x_0$ is still hyperspecial in $\mathcal{B}(G,E)$ hence so is $x$. Moreover, \eqref{cond2 x} follows from the identities $\mathfrak{u}_\alpha(E)_{x,0}=\varpi_E^{n_\alpha}\mathfrak{u}_\alpha(E)_{x_0,0}$, $\alpha\in \Delta$.
	
	Since $G$ is unramified, up to enlarging the tamely ramified extension $E$, we may assume that $G_E$ is split. Let $\Phi_E$ be the set of roots of $T_E$ in $G_E$, $\Phi_E^+$ be the subset of positive roots with respect to $B_E$ and $\Delta_E$ the corresponding set of simple roots. Similarly to \eqref{eq3 MPfiltration}, we have a decomposition
	\begin{equation*}
		\displaystyle \mathfrak{g}(E)_{x,0}=\mathfrak{t}(E)_{x,0}\oplus \bigoplus_{\beta\in \Phi_E} \mathfrak{u}_\beta(E)_{x,0}.
	\end{equation*}
	Then, \eqref{cond2 x} implies the following seemingly stronger property:
	\begin{equation}\label{cond2 bis}
		\displaystyle \mbox{For every } \beta\in \Delta_E, \; \psi^E \mbox{ is trivial on } \mathfrak{u}_\beta(E)_{x,0} \mbox{ but non-trivial on } \varpi_E^{-1}\mathfrak{u}_\beta(E)_{x,0}.
	\end{equation}
	Indeed, there is a natural surjection $p:\Delta_E\to \Delta$, given by restriction to $A_E$, and the fibers are $\Gamma=\Gal(E/F)$-orbits for the natural action. Moreover, for each $\alpha\in \Delta$, we have
	$$\displaystyle \mathfrak{u}_\alpha(E)_{x,0}=\bigoplus_{\beta\in p^{-1}(\alpha)} \mathfrak{u}_\beta(E)_{x,0}$$
	and $\sigma(\mathfrak{u}_\beta(E)_{x,0})=\mathfrak{u}_{\sigma \beta}(E)_{x,0}$ for $\sigma\in \Gamma$ and $\beta\in\Delta_E$. In particular, by \eqref{cond2 x}, $\psi^E$ is trivial on $\mathfrak{u}_\beta(E)_{x,0}$ for every $\beta\in \Delta_E$ and, for each $\alpha\in \Delta$, there exists at least one $\beta\in p^{-1}(\alpha)$ such that $\psi^E$ is non-trivial on $\varpi_E^{-1}\mathfrak{u}_\beta(E)_{x,0}$. Using that $\psi^E$ is $\Gamma$-invariant, and what precedes, we deduce that $\psi^E$ is non-trivial on $\varpi_E^{-1}\mathfrak{u}_\beta(E)_{x,0}$ for every $\beta\in \Delta_E$.
	
	For $n\geq 0$, we let $\chi_{n}^F$ be the character of $\mathfrak{g}(F)_{x,n}=\varpi_F^n \mathfrak{g}(F)_{x,0}$ given by $\chi_{n}^F(X):=\psi^F(\varpi_F^{-2n}X)$. As $\Tr_{E/F}(\mathfrak{g}(F)_{x,0})=\mathfrak{g}(F)_{x,0}$, by \eqref{cond2 x}, $\chi_{n}^F$ factors through a character of $\mathfrak{g}(F)_{x,n}/\mathfrak{g}(F)_{x,2n}$ that we identify, as before, with a character of $\cP_n^F=G(F)_{x,n}$. More generally, for every subsextension $F\subset K\subset E$ we define, as in the previous subsection, a character $\chi_n^K$ of $\mathfrak{g}(K)_{x,n}/\mathfrak{g}(K)_{x,2n}$, that we identify with a character of $\cP_n^K=G(K)_{x,n}$ by mean of $MP_{x,n}^K$, by $\chi_n^K=\chi_n^F\circ \Tr_{K/F}$. Thus, $\chi_n^K(X)=\psi^K(\varpi_F^{-2n}X)$, $X\in \mathfrak{g}(K)_{x,n}$. 
	
	\begin{lemma}\label{lem regular nilp}
		Let $n\geq m\geq 0$ be integers, $F\subset K\subset E$ be a subextension and $X\in \mathfrak{u}(K)_{x,m}$. Then, the character
		$$\displaystyle \mathfrak{g}(K)_{x,n}\ni Y\mapsto \chi_n^K([X,Y])$$
		is trivial if and only if $X\in \mathfrak{u}(K)_{x,n}$.
	\end{lemma}
	
	\begin{proof}
		Since $\mathfrak{g}(K)_{x,n}=\Tr_{E/K}(\mathfrak{g}(E)_{x,n})$, $\mathfrak{u}(K)_{x,n}=\mathfrak{u}(K)\cap \mathfrak{u}(E)_{x,n}$ and
		$$\chi_n^K([X,\Tr_{E/K}(Y)])=\chi_n^E([X,Y]) \mbox{ for } Y\in \mathfrak{g}(E),$$
		we are reduced to proving the lemma for $K=E$.
		
		Moreover, as we have $\mathfrak{u}(E)_{x,m}=\varpi_F^m\mathfrak{u}(E)_{x,0}$, $\mathfrak{g}(E)_{x,n}=\varpi_F^n\mathfrak{g}(E)_{x,0}$ and $\chi_n^E=\psi^E(\varpi_F^{-2n})$, it suffices to show that for $X\in \mathfrak{u}(E)_{x,0}$, the character
		$$\displaystyle \mathfrak{g}(E)_{x,0}\ni Y\mapsto \psi^E(\varpi_E^{-1}[X,Y])$$
		is trivial if and only if $X\in \varpi_E\mathfrak{u}(E)_{x,0}$.
		
		Let $G_x$, $T_x$, $\mathfrak{u}_{\beta,x}$ ($\beta\in \Phi_E$) be the special fibers of the smooth models $\mathbb{G}_x^E$, $\mathbb{T}_{\mathcal{O}_E}$ and $\mathfrak{U}_{\beta,x}^E$ respectively. Since $x\in \mathcal{A}$ is hyperspecial in $\mathcal{B}(G,E)$ and $G_E$ is split, $G_x$ is split reductive, $T_x$ is a split maximal torus of it, the set of roots of $T_x$ in $G_x$ is again $\Phi_E\subset X^*(T_E)=X^*(T_x)$ and for each $\beta\in \Phi_E$ the corresponding root subspace is $\mathfrak{u}_{\beta,x}$. 
		
		In particular, we have the decomposition
		$$\displaystyle \mathfrak{g}_x=\mathfrak{t}_x\oplus \bigoplus_{\beta\in \Phi_E} \mathfrak{u}_{\beta,x}$$
		where $\mathfrak{g}_x:=\Lie(G_x)$, $\mathfrak{t}_x:=\Lie(T_x)$.
		
		As $p$ is very good for $G$, there exists a $G_x$-invariant symmetric non-degnerate pairing $\langle .,.\rangle$ on $\mathfrak{g}_x$ (see \cite[Lemma 4.2.3]{Rich}). In particular, this induces a perfect pairing
		$$\displaystyle \langle .,.\rangle: \mathfrak{g}_x(k_E)\times \mathfrak{g}_x(k_E) \to k_E$$
		on $\mathfrak{g}_x(k_E)=\mathfrak{g}(E)_{x,0}/\varpi_E \mathfrak{g}(E)_{x,0}$.
		Let $\kappa: k_E\to \bbC^\times$ be a non-trivial character. As $\psi^E$ is trivial on $\mathfrak{g}(E)_{x,0}$, there exists a unique element $\overline{N}\in \mathfrak{g}_x(k_E)$ such that
		\begin{equation*}
			\displaystyle \psi^E(\varpi_E^{-1}Z)=\kappa(\langle \overline{N},\overline{Z}\rangle), \mbox{ for every } Z\in \mathfrak{g}(E)_{x,0},
		\end{equation*}
		where $\overline{Z}$ denotes the image of $Z$ in $\mathfrak{g}_x(k_E)$. In particular, we have
		$$\displaystyle \psi^E(\varpi_E^{-1}[X,Y])=\kappa(\langle \overline{N}, [\overline{X},\overline{Y}]\rangle )=\kappa(-\langle [\overline{X},\overline{N}],\overline{Y}\rangle),\;\;\; X,Y\in \mathfrak{g}(E)_{x,0},$$
		and the lemma is equivalent to the injectivity of the map
		\begin{equation}\label{map regular nilp}
			\displaystyle \overline{X}\in \mathfrak{u}_x(k_E):=\bigoplus_{\alpha\in \Phi^+} \mathfrak{u}_{\alpha,x}(k_E)\mapsto [\overline{X},\overline{N}].
		\end{equation}

		By \eqref{cond2 bis}, $\langle \overline{N},.\rangle$ is non-trivial on $\mathfrak{u}_{\beta,x}(k_E)=\mathfrak{u}_\beta(E)_{x,0}/\varpi_E\mathfrak{u}_\beta(E)_{x,0}$, for $\beta\in \Delta_E$, whereas it is trivial on $\mathfrak{u}_{\beta,x}(k_E)$ for $\beta\in \Phi_E\setminus \Delta_E$ as well as on $\mathfrak{t}_x(k_E)$. Thus, by $T_x$-invariance of $\langle .,.\rangle$, $\overline{N}$ is a sum of basis elements of $\mathfrak{u}_{-\beta,x}(k_E)$ for $\beta\in \Delta_E$. By \cite[Lemma 3.1.2]{Rich}, and the assumption that $p$ is very good for $G$, this implies that the map \eqref{map regular nilp} is injective.
	\end{proof}
	
	Let $2\check{\rho}\in X_*(A)$ be the sum of all the positive coroots of $A$ with respect to $B$. Let $F\subset K\subset E$ be a subextension. For every $n\geq 1$, we let
	\begin{itemize}
		\item $a_n:=(2\check{\rho})(\varpi_F)^{-n}\in A(F)$;
		
		\item $G^K_n:=a_n\cP^K_na_n^{-1}=a_nG(K)_{x,n}a_n^{-1}$;
		
		\item $\psi^K_n: G^K_n\to \mathbb{C}^\times$ be the character defined by $\psi^K_n(g_n)=\chi_n^K(a_n^{-1}g_na_n)$ for $g_n\in G_n^K$;
		
		\item $\varphi^K_n\in C_c^\infty(G(K))$ be the test function that is equal to $\vol(G^K_n)^{-1}\psi^K_n$ on $G^K_n$ and zero on the complement.
	\end{itemize}

	\begin{prop}\label{prop Rodier sequence}
		Assume that $p$ is very good for $G$. Then, for every $n\geq 1$, we have
		\begin{enumerate}[(i)]
			\item The Iwahori decomposition
			$$\displaystyle G^K_n=\overline{B}^K_n U^K_n$$
			where $\overline{B}^K_n:= G^K_n\cap \overline{B}(K)$, $U^K_n:= G^K_n\cap U(K)$;
			
			\item $\psi^K_n\mid_{\overline{B}^K_n}=1$ and $\psi^K_n\mid_{U^K_n}=\psi^K\mid_{U^K_n}$;
			
			\item $\overline{B}^K_{n+1}\subset \overline{B}^K_n$ and $U^K_n\subset U^K_{n+1}$.
			
		\end{enumerate}
		Furthermore, there exists an integer $n_0\geq 1$ which only depends on the Dynkin type of $G$ (that is, on the absolute root system of $G$) such that for every $n\geq n_0$,
		\begin{enumerate}
			\setcounter{enumi}{3}
			\item $\varphi^K_n\star \varphi^K_{n+1}\star \varphi^K_n=[U^K_{n+1}:U^K_n]^{-1}\varphi^K_n$.
		\end{enumerate}
	\end{prop}
	
	\begin{proof}
		$(i)-(iii)$ follow readily from the decomposition \eqref{eq2 MPfiltration} and the fact that $\langle \alpha, 2\check{\rho}\rangle=2$ for every $\alpha\in \Delta$.
		
		Let us prove $(iv)$. By $(i)$ and $(iii)$, we have $G^K_{n+1}/(G^K_n\cap G^K_{n+1}) =U^K_{n+1}/U^K_n$. From there and since $\psi^K_{n+1}\mid_{U^K_n}=\psi^K_n$ (by $(ii)$), it follows that
		$$\displaystyle \varphi^K_n\star \varphi^K_{n+1}\star \varphi^K_n=[U^K_{n+1}:U^K_n]^{-1} \sum_{u\in U^K_{n+1}/U^K_n} \psi^K_{n+1}(u) \varphi^K_n\star \delta_u\star \varphi^K_n.$$
		Thus, it suffices to check that for $u\in U^K_{n+1}- U^K_n$,  $\varphi^K_n\star \delta_u\star \varphi^K_n=0$ or equivalently that there exists $g_n\in G^K_n\cap uG^K_nu^{-1}$ such that $\psi^K_n(g_n)\neq \psi^K_n(u^{-1}g_nu)$.
		
		Let $u\in U^K_{n+1}$ and assume that $\psi^K_n(g_n)= \psi^K_n(u^{-1}g_nu)$ for every $g_n\in  G^K_n\cap uG^K_nu^{-1}$. Set $v=a_n^{-1}ua_n$. Then, $v\in a_1U(K)_{x,n+1}a_1^{-1}$ and there exists $n'_0\geq 1$, only depending on the Dynkin type of $G$\footnote{e.g. we can take $n'_0=2$ times the sum of the coefficients of the highest root.} such that, provided $n\geq n'_0$, $v\in U(K)_{x,n-n'_0}$. In particular, since $G(K)_{x,n}$ is normal in $G(K)_{x,0}$ (by \eqref{eq1 MPfiltration}), we have $vG(K)_{x,n}v^{-1}=G(K)_{x,n}$ and the assumption on $u$ translates to: for every $g_n\in G(K)_{x,n}$, $\chi_n^K(v^{-1}g_nv)=\chi_n^K(g_n)$.
		
		Assume that $n\geq n_0:=2n'_0$ and let $X+\mathfrak{g}(K)_{x,2n-n_0}=MP_{x,n-n_0'}^K(v)$ be the image of $v$ by the Moy-Prasad isomorphism. By \eqref{eq4 MPfiltration}, we may assume that $X\in \mathfrak{u}(K)_{x,n-n_0'}$. By \eqref{eq5 MPfiltration}, we have (note that $n+2(n-n_0')\geq 2n$)
		$$\displaystyle MP_{x,n}^K(vg_nv^{-1})=MP_{x,n}^K(g_n)+[X,MP_{x,n}^K(g_n)] \mbox{ mod } \mathfrak{g}(K)_{x,2n}.$$
		Therefore the assumption translates to: $\chi_n^K([X,Y])=1$ for every $Y\in \mathfrak{g}(K)_{x,n}$. By Lemma \eqref{lem regular nilp} this implies $X\in \mathfrak{u}(K)_{x,n}$, that is $v\in U(K)_{x,n}$ or equivalently $v\in U^K_n$.
	\end{proof}
	
	We now prove that for these choices of $x$ and $\chi_n^F$ the sequence of $K$-types $(\cP_n^K,\chi_n^K)_{n\geq 1}$ satisfies point $(iii)$ of Theorem \ref{MWKV}. This is obviously equivalent to the following proposition.
	
	\begin{prop}
		Assume that $p$ is very good for $G$. Let $\pi$ be an irreducible representation of $G(K)$ and let $n_0\geq 1$ be an integer such that the point $(iv)$ of Proposition \ref{prop Rodier sequence} is satisfied. Then, the sequence $(\dim \pi^{G^K_n,\psi^K_n})_{n\geq n_0}$ is increasing and stabilizes to $\dim \Hom_{U(K)}(\pi, \psi^K)$, in particular for $n$ large enough (depending on $\pi$) we have
		$$\displaystyle \dim \pi^{G^K_n,\psi^K_n}=\dim \Hom_{U(K)}(\pi, \psi^K).$$
	\end{prop}
	
	\begin{proof}
		As explained in \cite[Appendix A]{Hansen}, it follows from Proposition \ref{prop Rodier sequence} that we have
		$$\displaystyle c-ind_{U(K)}^{G(K)}(\psi^K)\simeq \varinjlim_{n} c-ind_{G^K_n}^{G(K)}(\psi^K_n),$$
		where $c-ind$ denotes compact induction, and the transition morphisms 
        \[ c-ind_{G^K_n}^{G(K)}(\psi^K_n)\to c-ind_{G^K_{n+1}}^{G(K)}(\psi^K_{n+1}) \]
        are given by $f\mapsto \varphi^K_{n+1}\star f$. Moreover, Proposition \ref{prop Rodier sequence}(iv) implies that these transition morphisms are split injections whenever $n\geq n_0$. It follows that 
		\begin{multline}\label{eq whitt functional}
			\displaystyle \Hom_{G(K)}(c-ind_{U(K)}^{G(K)}(\psi^K),\pi)\simeq \varprojlim_{n\geq n_0} \Hom_{G(K)}(c-ind_{G^K_n}^{G(K)}(\psi^K_n),\pi) \\
            =\varprojlim_{n} \pi^{G^K_n,\psi^K_n}
		\end{multline}
		with transition morphisms split surjections. In particular, $\dim(\pi^{G^K_n,\psi^K_n})$ increases for $n\geq n_0$ and equals $\dim \Hom_{G(K)}(c-ind_{U(K)}^{G(K)}(\psi^K),\pi)$ for $n$ large enough. By duality and Frobenius reciprocity, this last space is isomorphic to the space $\Hom_{U(K)}(\tilde{\pi},(\psi^K)^{-1})$, that is the space of $(\psi^K)^{-1}$-Whittaker functionals on the contragredient $\tilde{\pi}$ of $\pi$. As is well-known the vector spaces $\Hom_{U(K)}(\pi,\psi^K)$ and $\Hom_{U(K)}(\tilde{\pi},(\psi^{K})^{-1})$ are naturally in duality (this can also be deduced directly from \eqref{eq whitt functional} since $\pi^{G^K_n,\psi^K_n}$ and $\tilde{\pi}^{G^K_n,(\psi^K_n)^{-1}}$ are in duality). The proposition follows.
	\end{proof}
	
	It remains to show the last property (point $(iv)$) of Theorem \ref{MWKV}. Let $n\geq n_0$ and $\chi: T(K)/T(K)_n\to \mathbb{C}^\times$ be a character. As is well-known, the principal series $I(\chi)=I_{\overline{B}(K)}^{G(K)}(\chi)$ admits a unique Whittaker functional up to a scalar, that is $\dim \Hom_{U(K)}(I(\chi),\psi^K)=1$. Thus, by the previous proposition, it suffices to check that the $K$-type $(\cP_n^K,\chi_n^K)$ appears at least once in $I(\chi)$, which is readily seen since by the decomposition \eqref{eq2 MPfiltration}, the character $\chi$ is trivial on the intersection $\overline{B}(K)\cap \cP_n^K$.
	
	\subsubsection{Relation to principal congruence subgroups}     
	\begin{lemma}\label{ncongruence}  Let $\pi$ be an irreducible representation of $G(F)$ and $n\geq 1$ be an integer. Choose a hyperspecial point $x_0$ so that $G(F)_{x_0,0} = G(\cO_F)$ corresponds to a Chevalley presentation of the group and $G(F)_{x_0,n}$ is the $n$-th principal congruence subgroup.  Suppose that the $K$-type $(\cP_n^F,\chi_n^F)$ occurs in $\pi$.  Then, $\pi$ also contains nonzero vectors fixed by $G(F)_{x_0,2n+1}$, that is:
		$$\pi^{G(F)_{x_0,2n+1}} \neq 0.$$
	\end{lemma}
	
	\begin{proof}  Let $x$ be the point of the building from which the sequence of $K$-types $(\cP_n^F,\chi_n^F)_{n\geq 1}$ was constructed. By its very definition, the kernel of the character $\chi_n^F$ contains $G(F)_{x,2n}$, in particular $\pi^{G(F)_{x,2n}} \neq 0$. Up to $G(F)$-conjugacy, we may assume that $x$ and $x_0$ belong to the closure of a common alcove of the building. Let $x'$ be a point in the interior of this alcove. Then, we have the inclusions
		$$\displaystyle G(F)_{x_0,2n+1}\subseteq G(F)_{x',2n}\subseteq G(F)_{x,2n}	$$
		and the lemma follows.
	\end{proof}
	
	\section{Globalization}\label{sec_globalization}
	
	Let $q$ be a power of the prime $p$, and let $F$ be a local field of characteristic $p$ and residue field $\bbF_q$. We suppose given a quasi-split semisimple group $G_F$ over $F$, minimally split by a Galois extension $F_1 / F$; thus the $L$-group of $G_F$ is the semi-direct product ${}^L G_F = \widehat{G} \rtimes \Gal(F_1/F)$.
	
	Let $W$ denote the Weyl group of $\widehat{G}_F$. We assume that $p \nmid |W|$; in particular, $p > 2$. 
	\begin{proposition}\label{prop_globalisation_of_qs_parameter}
		Let $\tau : \Gamma_F \to {}^L G_F(\overline{\bbQ})$ be a continuous homomorphism, which is irreducible. Then, for any sufficiently large prime $l$, we can find the following data:
		\begin{enumerate} 
			\item A place $\lambda | l$ of $\overline{\bbQ}$.
			\item A global field $K$ of characteristic $p$, with field of constants $\mathbb{F}_q$; a Galois extension $K_1 / K$; a place $v$ of $K$ that remains prime in $K_1$, together with an identification $K_{1, v} / K_v \cong F_1 / F$; and a quasi-split semisimple group $G$ over $K$, split by $K_1$, and with a compatible identification $G_{K_{v}} = G_F$.			\item A continuous, finitely ramified $L$-homomorphism $\rho : \Gamma_K \to {}^L G (\overline{\mathbb{Q}}_\lambda)$ of Zariski dense image, such that $\rho|_{\Gamma_{K_v}} \sim \tau$. 
		\end{enumerate} 
	\end{proposition}
	\begin{proof}
		In the case that $F_1 = F$ (i.e.\ that $G_F$ is split), the existence of $\rho$ with the required properties follows from \cite[Proposition 9.5]{Boc19} (i.e.\ the Khare--Wintenberger argument). We need to explain how to generalise these arguments to the quasi-split case. The condition that $\tau$ is irreducible implies that it has finite image. Choose a prime $l > 2 \dim \widehat{G}_F$ and a place $\lambda | l$ of $\overline{\bbQ}$ such that $\tau$ takes values in ${}^L G_F(\overline{\mathbb{Z}}_\lambda)$, $\tau(\Gamma_F)$ has order prime to $l$, and $H^2(\Gamma_F, \widehat{\mathfrak{g}}_{\overline{\mathbb{F}}_l}) = 0$. Let $\overline{\tau} : \Gamma_F \to {}^L G_F(\overline{\mathbb{F}}_l)$ denote the reduction of $\tau$ modulo $\mathfrak{m}_{\overline{\bbZ}_l}$. Our first task is to globalise the residual representation $\overline{\tau}$.
		
		Let $k \leq \overline{\bbF}_\lambda$ be a finite extension such that $\overline{\tau}(\Gamma_F) \leq {}^L G_F(k)$, and choose data $E_1 / E$, $v$, and $G$ as in part (ii) of the statement of the Proposition. We claim that we can find the following data:
		\begin{itemize}
			\item A finite extension $K / E$, linearly disjoint from $E_1 / E$, in which $v$ splits. In particular, for each place $w | v$ of $K$, we can identify $K_w = F$ and $G_{K_w} = G_F$.
			\item A surjective $L$-homomorphism $\overline{\rho} : \Gamma_K \to {}^L G(k)$ such that, for each place $w |  v$ of $K$, $\overline{\rho}|_{\Gamma_{K_w}}$ and $\overline{\tau}$ are $\widehat{G}(k)$-conjugate. 
		\end{itemize}  
		This can be proved using the main theorem of \cite{Mor89}, using a variant of the strategy used in \cite{Mor90}. Indeed, let $X = \Spec E_1[ \{ x_g \}_{g \in {}^L G(k)}]$, and let ${}^L G(k)$ act on $X$ as follows: it acts $E$-linearly on $E_1$ through the map ${}^L G(k) \to \Gal(E_1 / E)$, and on the variables $x_g$ by permutation. Let $X' \subset X$ denote the non-empty ${}^L G(k)$-invariant open subscheme defined by the equations $x_g - x_h \neq 0$, for all $g, h \in {}^L G(k)$, $g \neq h$. Then ${}^L G(k)$ acts freely on the affine scheme $X'$; let $Y = X' / {}^L G(k)$. Then $Y$ is a geometrically integral affine $K$-scheme. If $M / E$ is a field extension, then a point $y \in Y(M)$ determines a Galois extension $M' / M$, containing $M \cdot E_1 / M$, together with an $L$-homomorphism $\Gal(M' / M) \to {}^L G(k)$. Conversely, if $M' / M$ is such an extension, then a choice of normal basis for $M'$ can be used to construct an $M$-point $y \in Y(M)$. In particular, we can write down a point $x_\tau \in Y(E_v)$ associated to $\overline{\tau}$. 
		
		The main theorem of \cite{Mor89} implies that we can find the following data:
		\begin{itemize}
			\item An extension $K / E$, split over $v$, and linearly disjoint from $E_1 / E$.
			\item A point $x \in Y(K)$, corresponding to a surjective $L$-homomorphism $\overline{\rho} : \Gamma_K \to {}^L G(k)$ such that for each place $w | v$ of $K$, $x$ approximates $x_\tau$ in $Y(K_w) = Y(E_v)$ to any desired degree of precision. 
		\end{itemize}
		Krasner's lemma shows that, if $x$ is $w$-adically close enough to $x_\tau$, then $\overline{\rho}$ is necessarily $\widehat{G}(k)$-conjugate to $\overline{\rho}$. 
		After possibly adjoining another extension to $K$, we can assume that $\overline{\rho}$ is unramified away from the set of places $w | v$ of $K$. Let $S$ denote this set, let $\mathcal{C}$ denote the category of complete Noetherian local $W(k)$-algebras, and let $\operatorname{Def} : \mathcal{C} \to \operatorname{Sets}$ denote the functor that associates, to each $A \in \mathcal{C}$, the set of $\widehat{G}(A)$-conjugacy classes of continuous $L$-homomorphisms $\rho_A : \Gamma_{K, S} \to {}^L G(A)$ lifting $\overline{\rho}$. Then $\operatorname{Def}$ is representable by an object $R \in \mathcal{C}$. 
		
		To complete the proof, it is enough (following the Khare--Wintenberger method) to show that $R$ is a finite $W(k)$-algebra of Krull dimension at least 1. Indeed, this necessarily implies the existence of an $L$-homomorphism $\rho : \Gamma_K \to {}^L G(\overline{\mathbb{Z}}_\lambda)$ lifting $\overline{\rho}$, and unramified outside $S$. We need to check that $\rho$ has Zariski dense image, and that for each $w \in S$, $\rho|_{\Gamma_{K_w}}$ and $\tau$ are $\widehat{G}(\overline{\mathbb{Q}}_\lambda)$-conjugate. The Zariski density follows from the Zariski density of $\rho(\Gamma_{K \cdot E_1})$ in $\widehat{G}(\overline{\mathbb{Q}}_\lambda)$, which follows from \cite[Proposition 6.7]{Boc19}. The local conjugacy follows from the fact that $H^1(\Gamma_{K_w}, \widehat{\mathfrak{g}}_k) = H^2(\Gamma_{K_w}, \widehat{\mathfrak{g}}_k) = 0$. 
		
		For the remainder of the proof, we therefore focus on explaining why $R$ has the claimed properties. We first claim that $R$ has a presentation 
        \[ R \cong W(k) \llbracket X_1, \dots, X_g \rrbracket / (f_1, \dots, f_g), \]
        for some $g \geq 0$. This is proved in the case $K = K_1$ as \cite[Proposition 5.11]{Boc19}, and the same proof (which is now a rather classical part of Galois deformation theory) applies in our case; we omit the details. This implies in particular that the Krull dimension of $R$ is $\geq 1$. 
		
		Let $S'$ denote set of places of $K_1 = K \cdot E_1$ lying above $v$, and let $R' \in \mathcal{C}$ denote the object representing the functor of deformations of $\overline{\rho}|_{G_{K_1}} : \Gamma_{K_1} \to \widehat{G}(k)$, as defined and studied in \cite[\S 5]{Boc19}. Then \cite[Proposition 5.14]{Boc19} implies that $R'$ is a finite flat $W(k)$-algebra. On the other hand, restriction of the universal deformation defines a morphism $R' \to R$ in $\mathcal{C}$. The proof will be complete if we can show that this morphism is surjective (as $R$ is then also a finite $W(k)$-algebra). Looking at Zariski tangent spaces, we see that it suffices to show that the restriction map $H^1(\Gamma_{K, S}, \widehat{\mathfrak{g}}_k) \to H^1(\Gamma_{K_1, S'}, \widehat{\mathfrak{g}}_k)$ is injective, or equivalently that $H^1(\Gal(K_1 / K), \widehat{\mathfrak{g}}_k^{\Gamma_{K_1}}) = 0$. This is true because $\widehat{\mathfrak{g}}_k^{\Gamma_{K_1}} = \widehat{\mathfrak{g}}_k^{\widehat{G}(k)} = 0$, by \cite[Lemma 5.1]{Boc19}. 
		
	\end{proof}
	We apply the above proposition to prove the following theorem. We remark that in our application in \S \ref{sec_descent_arguments}, the group $G_F$ will be chosen to be unramified, and so the extension $F_1 / F$ will be an  unramified (and therefore cyclic) extension of local fields. In particular, the extension $K' / K$ provided by Theorem \ref{potaut} will also be cyclic, and so the global group $G$ will split over a cyclic extension. 
	\begin{theorem}\label{potaut}
		Let $\tau : \Gamma_F \to {}^L G_F(\overline{\bbQ})$ be a continuous homomorphism, which is elliptic. Then we can find the following data:		\begin{enumerate}
			\item A prime $l$ and a place $\lambda | l$ of $\overline{\bbQ}$.
			\item A smooth, geometrically connected curve $X$ over $\bbF_q$ of genus $g > 1$. We set $K = \bbF_q(X)$.
			\item A finite set $S$ of places of $K$, a Galois extension $K' / K$ such that each place $v \in S$ remains prime in $K'$, and for each $v \in S$, an identification $K'_v / K_v \cong F_1 / F$.
			\item A finite, soluble (and therefore Galois) extension $L / K$, containing $K'$.
			\item A quasi-split semisimple group $G$ over $K$, split by $K'$, and for each $v \in S$, an identification $G_{K_v} = G_{F}$.
			\item An $L$-homomorphism $\rho : \Gamma_K \to {}^L G(\overline{\bbQ}_\lambda)$ of Zariski dense image. 
			\item A cuspidal automorphic representation $\Pi$ of $G(\bbA_L)$ over $\overline{\bbQ}$, everywhere unramified, such that $r_\lambda(\Pi) : \Gamma_L \to \widehat{G}(\overline{\bbQ}_\lambda)$ is conjugate to $\rho|_{\Gamma_L}$ (in particular, $r_\lambda(\Pi) $ has Zariski dense image and therefore depends, up to conjugacy, only on $\Pi$).
		\end{enumerate}
		Moreover, the following conditions are satisfied:
		\begin{enumerate} 
			\item For each place $v \not\in S$ of $K$, $\rho|_{\Gamma_{K_v}}$ is unramified.
			\item For each $v \in S$, there is exactly one place $v(+)$ of $L$ lying above $v$, $\rho|_{\Gamma_{K_v}}$ is $\widehat{G}(\overline{\mathbb{Q}}_\lambda)$-conjugate to $\tau$, and $\rho|_{\Gamma_{L_{v(+)}}}$ is trivial. 
		\end{enumerate}
	\end{theorem}
	\begin{proof}
		We claim that we can find the following data:
		\begin{enumerate} 
			\item A smooth, geometrically connected curve $X_0$ over $\bbF_q$ of genus $g > 1$, of function field $K_0 = \bbF_q(X_0)$, and a Galois extension $K_0' / K_0$. 
			\item A finite set $S_0$ of places of $K_0$, such that each $v_0 \in S_0$ remains prime in $K_0'$, and for each $v_0 \in S_0$, an isomorphism $K'_{0, v_0} / K_{0, v_0} \cong F_1 / F$.
			\item A quasi-split semisimple group $G$ over $K_0$, split by $K_0'$, and for each $v_0 \in S_0$, an identification $G_{K_{0, v_0}} = G_{F_v}$. 
			\item A finite, soluble extension $L_0 / K_0$, with constant field $\bbF_q$, and a finite set $S_0$ of places of $K_0$, such that for each $v_0 \in S_0$, there is a unique place $v_0(+)$ of $L_0$ lying above $v_0$, together with an isomorphism $\iota_{v_0} : K_{0, v_0} \to F$ of $\bbF_q$-algebras such that $\iota_{v_0}^{-1, \ast} \tau|_{\Gamma_{L_{v_0(+)}}}$ is trivial. 
			\item A prime $l$, a place $\lambda | l$ of $\overline{\mathbb{Q}}_\lambda$, and a continuous $L$-homomorphism $\rho_0 : \Gamma_{K_0} \to {}^L G(\overline{\bbQ}_\lambda)$ of Zariski dense image, unramified outside $S_0$, such that for each $v_0 \in S_0$, $\rho_0|_{\Gamma_{K_{0, v_0}}}$ is $\widehat{G}(\overline{\bbQ}_\lambda)$-conjugate to $\tau$.  
			\item Writing $\eta : G \to H$ for the adjoint quotient, a continuous, everywhere unramified representation $\widetilde{\rho}_0 : \Gamma_{L_0} \to \widehat{H}(\overline{\bbQ}_\lambda)$, lifting $\rho_0|_{\Gamma_{L_0}}$.
		\end{enumerate} 
		Indeed, by Proposition \ref{prop_globalisation_of_qs_parameter}, we can find a smooth, geometrically connected curve $X_a$ over $\bbF_q$, of function field $K_a = \bbF_q(X_a)$, a Galois extension $K_a' / K_a$, a quasi-split semisimple group $G$ over $K_a$, a place $v_a$ of $K_a$ that remains prime in $K'_a$, identifications $K_{a, v_a}' / K_{a, v_a} \cong F_1 / F$ and $G_{K_{a, v_a}} = G_F$, and a continuous $L$-homomorphism $\rho_a : \Gamma_{K_a} \to {}^L G(\overline{\bbQ}_\lambda)$ of Zariski dense image such that $\rho_a|_{\Gamma_{K_{a, v_a}}}$ is $\widehat{G}(\overline{\bbQ}_\lambda)$-conjugate to $\tau$.  There is no obstacle to assuming the genus of $X_a$ is at least $1$.  Let $M_a / K_a$ be a Galois extension containing $K_a'$ such that $\rho_a|_{\Gamma_{M_a}}$ is everywhere unramified, such that some place of $K_a$ splits in $M_a$, and such that if $w_a$ is a place of $M_a$ lying above $v_a$, then $\rho_a|_{\Gamma_{M_{a, w_a}}}$ is trivial. By  \cite[Theorem 1.4]{Con15}, there exists a representation $\widetilde{\rho}_a : \Gamma_{M_a} \to \widehat{H}(\overline{\bbQ}_{\lambda})$ lifting $\rho_a|_{\Gamma_{M_a}}$. After possibly enlarging $M_a$ further, we can assume that $\widetilde{\rho}_a$ is everywhere unramified. Now, choose a place $w_a$ of $M_a$ lying above $v_a$, and let $K_b = M_a^{D_{w_a / v_a}}$, and let $v_b$ denote the place of $K_b$ lying below $w_a$. Then $K_{b, v_b} = K_{a, v_a}$, $\rho_a|_{\Gamma_{K_b, v_b}}$ is conjugate to $\tau$, and $M_a / K_b$ is soluble. Let $K_c / K_b$ be a Galois extension in which $v_b$ splits, and such that for each place $w$ of $K_c$ not lying above $v_b$, $\rho_a|_{\Gamma_{K_{c, w}}}$ is unramified. Then $K_c / K_b$ is linearly disjoint from $M_a$. We take $K_0 = K_c$, $L_0 = M_a \cdot K_c$, $S_0$ to be the set of places of $L_0$ lying above $v_b$, and for $v_0 \in S_0$, $\rho_0 = \rho_a|_{\Gamma_{K_0}}$, and $\widetilde{\rho}_0 = \widetilde{\rho}_a|_{\Gamma_{L_0}}$. By construction, there is a place of $K_0$ that splits in $L_0$, so we see in particular that the field of constants of $L_0$ equals $\bbF_q$.

		Our aim will be to construct the following data:
		\begin{itemize}
			\item A Galois extension $K / K_0$,  split over $S_0$. 
			\item A cyclic extension $E_0 / K_0$, of prime degree $r \nmid [L_0 : K_0]$, such that each place $v \in S_0$ is inert (and hence unramified) in $E_0$.
			\item A cuspidal, everywhere unramified, automorphic representation $\widetilde{\Pi}_L$ of $H(\bbA_L)$, where $L = L_0 \cdot K \cdot E_0$, such that $r_\lambda(\widetilde{\Pi}_L) \sim \widetilde{\rho}_{0}|_{\Gamma_L}$. 
		\end{itemize}
		We explain why this suffices to prove the theorem. We will take $S$ to be the set of places of $K$ lying above a place of $S_0$, and $\rho = \rho_0|_{\Gamma_K}$. By construction, the extension $L / K$ will be soluble, with each place $v \in S$ having a unique place $v(+)$ lying above it, identifications $K'_v / K_v \cong F_1 / F$ and $G_{K_v} \cong G_F$, and  $\rho|_{\Gamma_{K_v}}$, $\tau$ will be $\widehat{G}(\overline{\mathbb{Q}}_\lambda)$-conjugate. Also by construction, $\rho$ will have Zariski dense image, and be unramified at each place $v\not\in S$ of $K$. We need to explain why $\rho|_{\Gamma_L}$ is associated to a cuspidal, everywhere unramified automorphic representation $\Pi_L$ of $G(\bbA_L)$. The proof of \cite[Theorem 11.1]{Boc19} shows that if we choose $f \in \widetilde{\Pi}_L - \{ 0 \}$ to be the (unique up to scalar) $H(\widehat{\cO}_L)$-invariant function, then $f \circ \eta : G(L) \backslash G(\bbA_L) / G(\widehat{\cO}_L) \to \overline{\bbQ}$ is non-zero, and generates a cuspidal, everywhere unramified automorphic representation $\Pi_L$ with the correct associated Galois representation.   Finally, our assumption on the genus of $X_a$ guarantees that $K$ and $L$ are function fields of genus at least $1$.
		
		Constructing the above data is proving the potential automorphy of $\widetilde{\rho}_{0}$ over an extension $L_0 \cdot K \cdot E_0 / L_0$, with $K_0$ split at $S_0$. If $L_0 = E_0 = K_0$ and we omit the splitness condition, this is the content of \cite[Theorem 11.1]{Boc19}, and our arguments here will be a variation on the proof of that theorem. We note however that the reduction there from the case of a general split, adjoint group $H$ to the case where $H$ is simple (this case being treated in \cite[Theorem 11.4]{Boc19}) contains a gap. The argument we give here will therefore both treat this variation, and also suffice to repair this gap. See Remark \ref{rmk_mistake} below for further explanation.
		
		For the rest of the proof, we therefore focus on constructing data $K, E_0, \widetilde{\Pi}_L$ as above. By \cite[Theorem 6.5]{Boc19}, we can extend $\widetilde{\rho}_{0} = \widetilde{\rho}_{0, \lambda}$ to a compatible system $(S, (\widetilde{\rho}_{0,\mu})_\mu)$ of representations of Zariski dense image, in an essentially unique way.  Since $H$ is split, we can write $H \cong H_1 \times \dots \times H_s$ for some $s \geq 1$, where each $H_i$ is adjoint, split, and absolutely simple. Let $h_i$ denote the Coxeter number of $\widehat{H}_i$. We observe that if $h_i = 2$, then $H_i \cong \PGL_2$ (by inspection of the classification). After re-ordering the factors, we can suppose that $h_1 \geq \dots \geq h_a > 2$, and that $h_{a+1} = \dots = h_s = 2$, for some $a \geq 1$. By Zsigmondy's theorem, we can find primes $t_1, \dots, t_a$, not necessarily distinct, such that $t_i | q^{h_i} - 1 $ but $t_i \nmid q^k - 1$ for any $k = 1, \dots, h_i - 1$. (Here we note that $q \neq 2$ because of our assumption that $p \nmid |W|$, and that $h_i > 2$; so we are not in the exceptional cases of Zsigmondy's theorem.) In particular, the order of $q \text{ mod }t_i \in (\bbZ / t_i \bbZ)^\times$ is $h_i$.

		Let $Y$ be another smooth, geometrically connected curve over $\bbF_q$ and let $F = \bbF_q(Y)$. After possibly enlarging $F$, we claim that we can find a compatible system $(\emptyset, (R_\mu)_\mu)$ of representations $R_\mu : \Gamma_F \to \widehat{H}(\overline{\bbQ}_\mu)$ with the following properties:
		\begin{enumerate}
			\item Each representation $R_\mu$ is everywhere unramified and of Zariski dense image.
			\item There exists a prime $l_1 > \max(t_1, \dots, t_a, 2 \dim \hat{H})$ and a place $\lambda_1 | l_1$ of $\overline{\bbQ}$ such that $\overline{R}_{\lambda_1} = \widehat{H}(k)$, for some finite extension $k / \bbF_{l_1}$ containing a $t_i^\text{th}$ root of unity for each $i = 1, \dots, a$.
			\item The order of $q \text{ mod }l_1 \in (\bbZ / l_1 \bbZ)^\times$ is prime to $h_i$, for each $i = 1, \dots, a$.
		\end{enumerate}
		We justify the claim. We first choose a sufficiently large prime $l_1$ such that $l_1 \equiv 3 \text{ mod }4$, $q \text{ mod }l_1$ is a square, and $l_1 \not\equiv 1\text{ mod }r$ for each prime $r | h_1 \dots h_a$. (Note that $q$ is prime to $h_1 \dots h_a$, since each $h_i$ divides the order of the Weyl group of $G$). Then the order of $q \text{ mod }l_1$ is prime to $h_i$, for each $i = 1, \dots, a$. We then choose a finite extension $k / \bbF_{l_1}$ containing a $t_i^\text{th}$ root of unity for each $i = 1, \dots, a$, and argue as in \cite[Proposition 9.5]{Boc19}: applying the theorem of Moret-Bailly (in the form of \cite[Theorem 9.3]{Boc19}, we can find (after enlarging $F$) a surjective, everywhere unramified representation $\overline{R}_1 : \Gamma_F \to \widehat{H}(k)$. By the Khare--Wintenberger method, we can lift this to an everywhere unramified representation in characteristic 0, which can then be placed in a compatible system $(R_\mu)_\mu$ with $\overline{R}_{\lambda_1} = \overline{R}_1$. 		
		
		By \cite[Proposition 6.6]{Boc19}, we can assume, after possibly replacing $(\widetilde{\rho}_{0, \mu})_\mu$ and $(R_\mu)_\mu$ by equivalent compatible systems, that there is a number field $E / \bbQ$ and a set $\cL$ of rational primes of Dirichlet density 0 with the following properties:
		\begin{itemize}
			\item For each $m \not\in \cL$ and for each place $\mu | m$ of $\overline{\bbQ}$, $\widetilde{\rho}_{0, \mu}$ and $R_{\mu}$ take values in $H(E_\mu)$.
			\item If $m \not\in \cL$ splits in $E$, then $\widetilde{\rho}_{0,\mu}$ and $R_\mu$ have image equal to $\widehat{H}(\bbZ_m)$.
		\end{itemize}	
		Using the Chebotarev density theorem, we can find a prime $l_0 $ satisfying the following condition: 
		\begin{itemize}
			\item $l_0$ does not divide $q l_1$, and $l_0 > \max(3, \# W_G)$. 
			\item $l_0$ splits in $E$ and $l_0 \not\in \cL$.
		\end{itemize}
		We observe that the groups $\widehat{H}(\bbF_{l_0})$ and $\widehat{H}(k)$ are perfect and have no isomorphic non-trivial quotients, and that $l_0, l_1$ are very good characteristics for $\widehat{H}$. Fix a choice of place $\lambda_0$ of $\overline{\bbQ}$ of residue characteristic $l_0$.
		
		For each $i = 1, \dots, a$, choose an element $\alpha_i \in L_0^\times$ that is a $t_i^\text{th}$ power locally in $L_{0, v_0(+)}$ for each $v_0 \in S_0$, and that has valuation 1 at some place of $L_0$, and let $\phi_i : \Gamma_{L_0} \to \widehat{H}_i(\bbZ[\zeta_{t_i}])$ be the parameter associated to the splitting field of $Y^{t_i} - \alpha_i$ over $L_0$ as in the proof of \cite[Proposition 10.12]{Boc19}. Then $\phi_i$ becomes automorphic after base change to any finite separable extension $L_0' /L_0$, linearly disjoint from the extension of $L_0(\zeta_{l_1})$ cut out by $\phi_i^{ad}$, such that $\phi_i|_{\Gamma_{L_0'}}$ is everywhere unramified. Moreover, \cite[Proposition 10.10]{Boc19} shows that in this case $\overline{\phi}_{i, \lambda_1}|_{\Gamma_{L_0'}}$ is strongly $\widehat{H}$-irreducible and that $\overline{\phi}_{i,\lambda_1}(\Gamma_{L_0'})$ is an $\widehat{H}$-abundant subgroup of $\widehat{H}(\overline{\bbF}_{l_1})$. In fact, since $q \text{ mod }l_1$ has order prime to $h_i$,  we have $\overline{\phi}^{ad}_{i, \lambda_1}(\Gamma_{L_0'}) = \overline{\phi}^{ad}_{i, \lambda_1}(\Gamma_{L_0'(\zeta_{l_1})})$, so we see that $\overline{\phi}_{i,\lambda_1}(\Gamma_{L_0'(\zeta_{l_1})})$ is also $\widehat{H}$-abundant.

		We now apply the construction of \cite[\S 9]{Boc19}. Let $H_\pi = \prod_{i=1}^a H_i$, and let $\widetilde{\rho}_{\pi, \mu}$ (resp. $R_{\pi, \mu}$) denote the projection of $\widetilde{\rho}_{0, \mu}$ (resp. $R_\mu$) to $H_\pi$, $\phi_\pi = \prod_{i=1}^a \phi_i$. Let us introduce the following data:
		\begin{itemize}
			\item The group $P = \widehat{H}_\pi(\bbF_{l_0}) \times \widehat{H}_\pi(k)$.
			\item A geometrically connected curve $Z_0$ over $\bbF_q$ with $\bbF_q(Z_0) = L_0$ and such that $\phi_\pi$ factors through $\pi_1(Z_0, \overline{\eta}_{Z_0})$.
			\item The homomorphism $\varphi : \pi_1(Z_0, \overline{\eta}_{Z_0}) \to P$ given by $\overline{\widetilde{\rho}}_{\pi,\lambda_0} \times \overline{\phi}_{\pi, \lambda_1}$ (note that our choices imply that $\overline{\phi}_{\pi, \lambda_1}$ indeed takes values in $\widehat{H}_\pi(k)$).
			\item The homomorphism $\psi : \pi_1(Y, \overline{\eta}_Y) \to P$ given by $\overline{R}_{\pi, \lambda_0} \times \overline{R}_{\pi, \lambda_1}$. (Then $\psi|_{\pi_1(Y_{\overline{\bbF}_q}, \overline{\eta}_Y)}$ is surjective, because $\psi$ is surjective and $P$ is perfect.)
		\end{itemize} 
		We set $Z_{\psi, \varphi} = \underline{\operatorname{Isom}}_{Y_{L_0}, P}(X_{\varphi}, X_\psi)$, where $X_{\varphi}$, $X_\psi$ are the $P$-torsors over $Y_{L_0}$ associated to the homomorphisms $\varphi$ and $\psi$, respectively. Then \cite[Lemma 9.1]{Boc19} shows that $Z_{\psi, \varphi}$ is finite \'etale over $Y_{L_0}$ and geometrically connected over $L_0$ (note there is a typo in \emph{loc. cit.} -- $Z_{\psi, \varphi}$ is not finite \'etale over $L_0$).
		
		Fix a choice of prime $r \nmid [L_0 : K_0] t_1 \dots t_a h_1 \dots h_a (l_1 - 1)$ and a cyclic extension $E_0 / K_0$ of degree $r$ in which each place of $S_0$ is inert, and some place not in $S_0$ splits. Then for each $v_0 \in S_0$, there is a unique place $v_0(++)$ of $L_0 \cdot E_0$ lying above $v_0$. Choosing $r$ sufficiently large, we can assume that for each place $v_0 \in S_0$, the set $Z_{\psi, \varphi}((L_0 \cdot E_0)_{v_0(++)})$ is non-empty. Indeed, since $\varphi$ is unramified at $v_0(++)$, this follows from the Weil bounds and Hensel's lemma. 
		
		We claim that we can find a Galois extension $K / K_0$ and a point 
        \[ z \in \operatorname{Res}_{L_0 \cdot E_0/ K_0} (Z_{\psi, \varphi})(K) = Z_{\psi, \varphi}(L_0 \cdot E_0 \otimes_{K_0} K) \]
        with the following properties:
		\begin{itemize}
			\item $S_0$ splits in $K$. In particular, $K / K_0$ is linearly disjoint from $L_0 \cdot E_0 / K_0$.
			\item The extension $L = L_0 \cdot E_0 \cdot K / L_0$ is linearly disjoint from the extension of $L_0(\zeta_{l_1})$ cut out by $\overline{\widetilde{\rho}}_{\pi, \lambda_0} \times \overline{\phi}^{ad}_{\pi, \lambda_1}$. 
			\item The representation $\phi_\pi|_{\Gamma_L}$ is everywhere unramified.
			\item The image $y \in Y_{K_0}(L) = Y(L)$ of the point $z$ does not lie in $Y(\overline{\bbF}_q)$. 
		\end{itemize}
		This follows from \cite[Th\'eor\`eme 1.3]{Mor89}, in a very similar manner to the proof of  \cite[Proposition 9.2]{Boc19}. The necessary local conditions at places $v_0 \in S_0$ may be satisfied by choice of $E_0$. The linear disjointness may be obtained by choosing $K / K_0$ to be split at a suitable set of auxiliary places. 
		
		We now explain how to deduce the existence of $\widetilde{\Pi}_L$. The point $y$ determines an embedding $j : F \to L$, and the above conditions imply that we can apply \cite[Corollary 8.21]{Boc19} to conclude that $j^\ast R_{\pi, \lambda_1}$ is automorphic, associated to a cuspidal, everywhere unramified automorphic representation of $H_\pi(\bbA_L)$. Since $j^\ast R_{\pi, \lambda_0}$ is a member of the same Zariski dense compatible system, it is also automorphic, associated to the same automorphic representation. We can then apply \cite[Corollary 8.21]{Boc19} to conclude that $\widetilde{\rho}_{\pi, \lambda_0}|_{\Gamma_L}$ is automorphic, associated to a cuspidal, everywhere unramified automorphic representation $\widetilde{\Pi}_{L, \pi}$ of $H_\pi(\bbA_L)$. To conclude the proof, we note that $H \cong H_\pi \times \PGL_2^{s-a}$, and that the automorphy of the $\widehat{\PGL}_2^{s-a}$ factor of $\widetilde{\rho}_{0, \lambda_0}|_{\Gamma_L}$ follows from known results (see e.g.\ \cite{Dri83}). We can take $\widetilde{\Pi}_L$ to be the external tensor product of $\widetilde{\Pi}_{L, \pi}$ with the cuspidal, everywhere unramified automorphic representation of $\PGL_2^{s-a}(\bbA_L)$ so obtained. 
	\end{proof}
	\begin{remark}\label{rmk_mistake} We expand on the gap in the proof of \cite[Theorem 11.1]{Boc19}, described in the proof of Theorem \ref{potaut}, which was pointed out to us by Dmitri Whitmore. First, \cite[Theorem 11.4]{Boc19} shows the potential automorphy of a Zariski dense representation $\rho : \Gamma_K \to \widehat{G}(\overline{\bbQ}_\lambda)$ when $G$ is simple, split, and adjoint, over some Galois extension $K' / K$. \cite[Theorem 11.1]{Boc19} aims to treat the case where $G$ is split and semisimple (but not necessarily simple or adjoint). The proof first treats the same result when $G$ is assumed only to be semisimple, split and adjoint (i.e. a product of groups of the type appearing in \cite[Theorem 11.4]{Boc19}). It is asserted there that this case follows from the case where $G$ is simple.
		
		If one writes $G = \prod_{i=1}^s G_i$ as a product of its simple factors, with a corresponding decomposition $\rho = \prod_i \rho_i$, then \cite[Theorem 11.4]{Boc19} shows that one can find Galois extensions $K_1, \dots, K_s / K$ such that $\rho_i|_{\Gamma_{K_i}}$ is automorphic. However, it is not clear that one can choose a common extension $K' / K$ with the property that each $\rho_i|_{\Gamma_{K'}}$ is automorphic, which is what is needed in order to deduce that $\rho|_{\Gamma_{K'}}$ is automorphic. Moreover, the nature of the `universally automorphic' parameters constructed in the proof of \cite[Theorem 11.4]{Boc19} means that this is a subtle issue, as one needs to carefully control the field of constants $\bbF_q$, in order to ensure that, for some auxiliary prime $t$, $q \text{ mod }t$ has order $h$, the order of the Coxeter number of $G$ in the case that $G$ is simple. In the proof given here, we use an alternative construction of `universally automorphic' parameters using Zsigmondy's theorem, in which the field of constants $\bbF_q$ is kept fixed.  Specialising this argument to the situation of \cite[Theorem 11.1]{Boc19}, one recovers the theorem stated there. 
	\end{remark}
	
	\section{Multiplicity one}\label{sec_mult_one}
	
	The following theorem follows from results recently announced by Gaitsgory and Raskin:
	\begin{theorem}\label{multone}  Let $L$ be a global field of positive characteristic $p$ ; let $G$ be a split semisimple group over $L$; and let $\Pi$ be a cuspidal automorphic representation of $G$ that is everywhere unramified. Let $\lambda$ be a prime-to-$p$ place of $\overline{\mathbb{Q}}$, and suppose that the image of $r_\lambda(\Gamma_L)$ is Zariski dense in $\widehat{G}(\overline{\mathbb{Q}}_\lambda)$. Suppose finally that $p$ is prime to the order of the Weyl group of $G$.
    
		Then $\Pi$ has automorphic multiplicity one; in other words, if $\Pi' \leq \mathcal{A}_0(G)$ is an automorphic representation such that $\Pi' \isoarrow \Pi$ as representations of $G(\mathbb{A}_L)$, then $\Pi' = \Pi$.
	\end{theorem}
	\begin{proof}
	    Assuming that $p$ is a very good characteristic for $G$ (a consequence of our assumption that $p$ does not divide the order of the Weyl group), it follows from \cite[Corollary 1.5.8]{GR25} that the (generalised) eigenspace for the excursion operators of \cite{Lafa} in $\mathcal{A}_0(G)^{G(\widehat{\mathcal{O}}_L)}$, associated to the eigenvalues of $r_\lambda(\Pi)$, has dimension exactly 1. Since $r_\Lambda(\Gamma_L)$ has Zariski dense image, it follows from \cite[Proposition 6.4]{Boc19} that this eigenspace coincides with the eigenspace in $\mathcal{A}_0(G)^{G(\widehat{\mathcal{O}}_L)}$ for the unramified Hecke eigenvalues of $\Pi$. In particular, the multiplicity of $\Pi$ is 1. 
	\end{proof}

	\begin{remark}  The main result of \cite{GR24} is a categorical version of the geometric Langlands conjecture, which includes an existence statement:  any irreducible $\widehat{G}$-local system over a smooth projective curve $X$ over a field of characteristic zero is the image under their functor of a Hecke eigensheaf on the stack $Bun_G$ of principal $G$-bundles on $X$.   The categorical  geometric Langlands correspondence also includes the assertion that the Hecke eigensheaf is uniquely determined by this correspondence.  The multiplicity at most one statement of \cite[Corollary 1.5.8]{GR25} follows from the geometric version by a specialization argument.
		
		Gaitsgory and Raskin have explained that an appropriate version of a conjecture of Takeshi Saito on singular support would imply the corresponding existence theorem in the arithmetic setting; in other words, it would imply a complete global Langlands correspondence for irreducible everywhere unramified homomorphisms from $\Gamma_K$ to  $\widehat{G}(\overline{\bbQ}_\lambda)$.  The methods of the present paper would then allow us to deduce the global Langlands correspondence for irreducible 
		homomorphisms  $\phi: \Gamma_K \ra \widehat{G}(\overline{\bbQ}_\lambda)$ such that the image under $\phi$ of any inertia group is {\it finite}, at least for very good $p$.
	\end{remark}
	
	\begin{remark}\label{coeff}  It is announced in \cite[Theorem B]{Ras} that Gaitsgory, Lafforgue, and Raskin have proved that, if $\Pi$ is everywhere unramified and tempered then it is globally generic, and not just locally generic everywhere, which is the property we use in our inductive argument.
	\end{remark}
	
	\section{A descent argument, following results of Konno and Varma}\label{secKV}

	The descent argument in the following sections is inspired by a result, due initially to Konno \cite{Ko} and generalized and completed by S. Varma \cite{Va17}, on the preservation of genericity under cyclic stable descent.   These results are not known for local fields of positive characteristic.  We believe (but have not verified) that the arguments in \cite{Ko,Va17} remain valid when the characteristic $p$ is sufficiently large relative to the group $G$, but it is clear upon inspection that some of them definitely break down for small $p$.  As motivation for what follows, we assume the result to hold for our group $G(F)$ and derive the existence of cyclic descent under this hypothesis.  The argument is uniform and does not depend on a case-by-case analysis, unlike the substitute developed in the subsequent sections.  
	
	We admit Theorem 6.15 of \cite{Va17}. The complete statement of that theorem is based on \cite[Hypothesis 6.14]{Va17}, which is a statement about abstract finite sets (packets) of representations.  We reformulate the Theorem and its underlying hypotheses in terms of  global data on the model of those constructed in Theorem \ref{potaut}.   We thus introduce groups $H$ and $H_1$  as in \cite[\S 6.3]{Va17}, with $H$ taking the place of Varma's $G$ and $H_1$ twisted endoscopic for $H$ with respect to an automorphism
	$\theta$ of $H$ of finite order.  These are defined in terms of a global cyclic extension $\bbL/\bbK$ of function fields with $\Gamma := \Gal(\bbL/\bbK) = \langle \theta \rangle$ of prime order $\ell$; the case $\ell = p$ is not excluded.   We assume there is a finite cyclic extension $F/F_1$ of local fields of degree $\ell$ and sets $S$ and $S_1$ of places of $\bbL$ and $\bbK$, respectively, such that 
	\begin{enumerate}
		\item For all $v_1 \in S_1$ there is a unique $v$ dividing $v_1$, and $v \in S$;
		\item For all $v_1 \in S_1$,  $\bbK_{v_1} \isoarrow F_1$, $\bbL_v \isoarrow F$, and the decomposition group of $\bbL_v/\bbK_{v_1}$ is generated by $\theta \in \Gamma$;
		\item $\bbL/\bbK$ is unramified outside $S$;
		\item  There is an everywhere locally generic cuspidal automorphic representation $\Sigma$ of $G(\bbL)$, unramified outside $S$, and an isomorphism
		$$b: \Sigma \isoarrow \Sigma^\theta$$
		that acts as the identity locally everywhere on the Whittaker model (see Remark \ref{whit} below).
		\item  The representation $\Sigma$ has automorphic multiplicity $1$.
		\item The global parameter $r_\lambda(\Sigma):  \Gamma_{\bbL} \ra \widehat{G}(\overline{\bbQ}_\lambda)$ has Zariski dense image.
		\item  There is an irreducible $\theta$-invariant generic representation $\pi_F$ of $G(F)$ such that, for all $v \in S$, $\Sigma_v \isoarrow \pi_F$.
	\end{enumerate}
	We let $H = \prod_{v \in S} G(\bbL_v)$, $H_1 = \prod_{v \in S_1} G(\bbK_{v_1})$.  Then  $\theta$-defines an automorphism of $H$ for which $H_1$ is a twisted endoscopic  group.    We write $\Sigma \isoarrow \Sigma^S \otimes \Sigma_S$, where $\Sigma^S$ is an irreducible representation of $G(\ad^S) := \prod'_{w \notin S} G(\bbL_w)$
	and $\Sigma_S$ is an irreducible representation of $H$; thus $\Sigma_S \isoarrow \otimes_{v \in S} \pi_F$.  Let $K_S \subset H$ be an open compact subgroup such that
	$\Sigma_S^{K_S} \neq 0$.  
	Let $\Pi$ denote the set of irreducible representations $\pi$ of $H$ such that $\Sigma_\pi := \Sigma^S \otimes \pi$ is a cuspidal automorphic representation of $G(\ad_{\bbL})$ and $\pi^{K_S} \neq 0$; thus $\Pi$ is a finite set.  
	\begin{remark}\label{whit}  Whittaker models and genericity are defined here with respect to the data in the statement of \cite[Hypothesis 6.14]{Va17}.  Since the discussion in the present section is entirely conditional on the validity of the results of Konno and Varma for local fields of positive characteristic, we are not more specific here, and we refer the reader to Varma's paper for precise definitions.
	\end{remark}

	In what follows, we need to make the following hypothesis.
	
	\begin{hyp}\label{selecta}  Let $\pi$ be an irreducible representation of $H_? \in \{H, H_1\}$.  There is an open compact subgroup $K_? \subset H_?$, a finite collection $\Pi_?$ of irreducible representations of $H_?$, containing $\pi$, such that $\pi_i^{K_?} \neq 0$ for all $\pi_i \in \Pi_?$, and a set of coefficients $0 < m(\pi_i) \in \mathbb{N}$, for $\pi_i \in \Pi_?$, such that $\sum_{\pi_i \in \Pi_?} m(\pi_i)\pi_i$ is a stable distribution.
	\end{hyp}
	
	Over $p$-adic fields this is proved by Arthur in \cite{Art96}.  A weaker result, where $m(\pi_i) \in \mathbb{Z}$ rather than $\mathbb{N}$, follows from the main theorem of \cite{Fu24}, also for $p$-adic fields, provided the Fargues-Scholze parameter of $\pi$ is elliptic (which is not known in general even for supercuspidal representations).  Fu proves stability on the elliptic regular locus and uses \cite{Art96} to prove stability of the distribution.  Since the present section is only included for heuristic reasons, and since we do not state any theorems, it doesn't hurt to introduce a new hypothesis here.

	\begin{ideal}\label{KoVa}  Assume \cite[Theorem 6.15]{Va17} is valid in our present situation.  Assume $\Sigma_S$ is the unique  generic element of $\Pi$.  Assume the open compact subgroup $K_S \subset H$ can be chosen so that the set $\Pi$ is stable;  i.e., there are coefficients $m(\pi) \in \mathbb{N}$ for $\pi \in \Pi$, as in Hypothesis \ref{selecta}, such that
		$\sum_{\pi \in \Pi} m(\pi)\pi$ is a stable distribution.
		Then
		\begin{enumerate}
			\item There is an open compact subgroup $K_1 \subset H_1$ and a finite non-empty set $\Pi_1$ of representations of $H_1$, with $\pi_1^{K_1} \neq 0$ for all $\pi_1 \in \Pi_1$, satisfying conditions (a) and (d) of Hypothesis 6.14 of \cite{Va17}.  In particular the representations in $\Pi_1$ satisfy the endoscopic character identities with respect to $\Pi$.
			\item  There is an irreducible representation $\Sigma_1^{S_1}$ of $G(\ad_\bbK^{S_1})$ that is everywhere locally unramified and generic, whose stable base change is isomorphic to $\Sigma_S$;
			\item  For any $\pi_1 \in \Pi_1$, the representation $\Sigma_1^{S_1}\otimes \pi_1$ is cuspidal automorphic for $G(\ad_{\bbK})$.
			\item  There is an irreducible generic representation $\pi_{F_1}$ of $G(F_1)$ and a representation $\Sigma_{S_1} \in \Pi_1$ such that
			$\Sigma_{S_1} \isoarrow \otimes_{v_1 \in S_1} \pi_{F_1}$; moreover, $\Sigma_{S_1}$ is the unique generic member of $\Pi_1$.
			\item  The representation $\Sigma_1^{S_1} \otimes \Sigma_{S_1}$ has automorphic multiplicity $1$.
		\end{enumerate}
	\end{ideal}
	\subsubsection*{Sketch of proof}  
	In order to apply Theorem \ref{hyp rho traces}  we let $\rho_\bbL = r_\lambda(\Sigma)$ and let $\rho = (\rho_v)$ denote the set of collections of local parameters of $G(\bbK_v)$ as $v$ runs over places of $\bbK$, such that, for almost all $v$, with $w~|~v$ a place of $\bbL$, the restriction of $\rho_v$ to the decomposition group of $w$ is $\rho_{\bbL,w}$.  We choose $f^S$ and $f^{\prime,S_1}$ to be spherical matching functions away from $S$ so that, as $f_S$ varies over $K_S$-invariant test functions on  $H$, the distribution 
	$$D(f_S) = I_{\rho_\bbL,\theta}(f^S\otimes f_S)$$
	is non-trivial.  Then if $f'_{S_1}$ is a base change transfer of $f_S$, it follows from Theorem \ref{hyp rho traces}  that 
	\begin{equation}\label{tranf}
		D(f_S) = D_1(f'_{S_1}) := I_\rho(f^{\prime,S_1}\otimes f'_{S_1}).
	\end{equation}

	By hypothesis the restriction of $D$ to $K_S$-biinvariant functions is stable.    Let $\Pi_1$ be the set of representations $\pi_1$ of $H_1$ such that $\Sigma_1^{S_1} \otimes \pi_1$ is automorphic; for $\pi_1 \in \Pi_1$ let $m(\pi_1) \in \mathbb{N}$ denote the automorphic multiplicity of $\Sigma_1^{S_1}\otimes \pi_1$.   We let $K_1 \subset H_1$ be an open compact subgroup as in Hypothesis \ref{selecta} so that $\pi_j^{K_1} \neq 0$ for every $\pi_j \in \Pi_1$.  
	
	At this point, we would like to use \cite[Theorem 6.15]{Va17} to conclude.  More precisely, under  Hypothesis \ref{selecta}, we would like to show that \eqref{tranf} implies that the  sets $\Pi$ and $\Pi_1$ satisfy condition (i) of \cite[Hypothesis 6.14]{Va17}, namely that
	$$ \sum_{\pi_1 \in \Pi_1} m(\pi_1) \pi_1$$
	is stable.  This should follow from \eqref{tranf}, which in turn would suffice to establish the key condition (iii)(iv) of \cite[Hypothesis 6.14]{Va17}.   However,  we don't see how to complete the argument; therefore we leave this as a heuristic motivation for the argument of the next section.

	\begin{corollary}  Assume Ideal Theorem \ref{KoVa} is a thoerem and \cite[Theorem 6.15]{Va17} is valid in the situation of Theorem \ref{potaut}.  Then there is an everywhere locally generic automorphic representation $\Pi_K$ of
		$G(\mathbb{A}_K)$ with automorphic multiplicity one whose Galois parameter is equivalent to the parameter $\rho$  of Theorem \ref{potaut} (iv).  In particular, every irreducible local parameter $\tau$ is the image of a generic and (necessariily) supercuspidal representation of $G(\cF)$ under the Genestier--Lafforgue parametrization.
	\end{corollary}
	\begin{proof}  This is an induction argument, starting with the everywhere unramified cuspidal automorphic representation $\Pi$ of $G(\mathbb{A}_L)$ constructed in Theorem \ref{potaut}.  We choose a place $v \in S_0$ as in the statement of that theorem and a place $v_\bbL$ of $\bbL$ dividing $v$.  We induct on the order of the finite soluble group $\Gamma = \Gal(L/K)$.  Theorem \ref{potaut} asserts that the statement is true with $K$ replaced by $L$.  Let 
		$$L \supset \bbL \supset \bbK \supset K,$$ 
		with $\bbL$ and $\bbK$ Galois over $K$ and $\bbL/\bbK$ cyclic of prime order.  We assume the statement is true with $K$ replaced by $\bbL$.  
		
		The induction step is to show that it is true with $K$ replaced by $\bbK$. But this follows immediately from Theorem \ref{KoVa}.  We omit the details, because the Corollary is conditional on the validity of Varma's theorem in equal characteristic.
	\end{proof}

	\section{Descent arguments}\label{sec_descent_arguments}
	
	We now give arguments for descent that do not rely on \cite{Va17}. This will lead to a proof of our main theorem (Theorem \ref{introthm_descent}, restated as Theorem \ref{existenceofpi}) below). We suppose that the following conditions are satisfied. 
    \begin{itemize}
        \item $K$ is a global field of positive characteristic $p$, and $L / K$ is a soluble extension.
        \item $G$ is a quasi-split semisimple group over $K$, split by a cyclic subextension $K' / K$ of $L / K$, and with absolute Weyl group of order prime to $p$.
        \item $S_0$ is a finite set of places of $K$, with the property that for each $v \in S_0$, there is a unique place $v(+)$ of $L$ lying above $v$, such that 
        \[ \Gal(L_{v(+)} / K_v) \cong \Gal(L / K). \]
        \item $\lambda$ is a place of $\overline{\mathbb{Q}}$ of residue characteristic prime to $p$, and $\rho : \Gamma_K \to {}^L G(\overline{\mathbb{Q}}_\lambda)$ is a continuous representation, of Zariski dense image, unramified outside $S_0$.
        \item $F$ is a local field of positive characteristic, $\tau : W_F \to {}^L G(\overline{\mathbb{Q}})$ is an irreducible $L$-homomorphism, and for each $v \in S_0$, $F \cong K_v$ is an isomorphism making $\tau$ and $\rho|_{W_{K_v}}$ conjugate under the action of $\widehat{G}(\overline{\mathbb{Q}}_\lambda)$.
        \item $\Pi_L$ is an everywhere unramified cuspidal automorphic representation of $G(\mathbb{A}_L)$ such that $r_\lambda(\Pi_L)$ and $\rho|_{\Gamma_L}$ are $\widehat{G}(\overline{\mathbb{Q}}_\lambda)$-conjugate.
    \end{itemize}

	\subsection{Wild descent}
	
	We let $D = \Gal(L_{v(+)}/K_v) \isoarrow \Gal(L/K)$ and view it as the common decomposition group of all $v \in S_0$.  Let $I \subset D$ be the (common) inertia group and $P \subset I$ the (common) wild inertia group. By assumption, the restriction of $\tau$ to $\Gamma_{K_v}$ factors through $D$ for every $v\in S_0$. Let $P\subset D_p\subset D$ be the inverse image of the $p$-Sylow subgroup $\tau(W_F)_p$ of $\tau(W_F)$, a normal subgroup by Lemma \ref{irrparam}. Let $L^P \subset L$ denote the fixed field of $D_p$. It is a Galois extension of $K$ with Galois group isomorphic to $D/D_p$. Note that the splitting extension $K'/K$ is of degree prime to $p$ (by our assumption that $p$ does not divide the order of the absolute Weyl group of $G$), hence $K'\subset L^P$ and $G$ is split over $L^P$. Let $\rho_P$ denote the restriction of $\rho$ to $\Gamma_{L^P}$.

	\begin{prop}\label{wild}  
There is a unique cuspidal automorphic representation $\Pi_{L^P}$ of $G(\ad_{L^P})$ such that the global parameter $r_{\lambda}(\Pi_{L^P})$ of $\Pi_{L^P}$ is equivalent to $\rho_P$ and such that, for $w$ a prime of $L^P$ whose restriction to $K$ is not in $S_0$, the component $\Pi_{L^P,w}$ is unramified. Moreover, $\Pi_{L^P}$ satisfies the following properties:
		\begin{itemize}
\item $\Pi_{L^P}$ is locally generic everywhere.
			\item For $v \in S_0$ with unique divisor $v_P$ in $L^P$, the component $\Pi_{L^P,v_P}$ is a principal series representation with a (local) Whittaker model.
			\item $\Pi_{L^P}$ has automorphic multiplicity $1$.
			\end{itemize}
	\end{prop}
	
	The proof will be presented after a digression on principal series. Let $F/F_1$ be a cyclic extension of local fields of degree $\ell$, and let $\theta$ be a generator of $\Gal(F/F_1)$. Let $G$ be a quasi-split reductive group over $F_1$.  Let $B = T U\subset G$ be a Borel subgroup with maximal torus $T$ and unipotent radical $U$. Let $\chi_1$ be a character of $T(F_1)$, $\chi$ the composition of $\chi_1$ with the norm $N_{F/F_1}:  T(F) \ra T(F_1)$.  Let $I(\chi_1)$ (resp. $I(\chi)$) be the (normalized) principal series representation $I_{B(F_1)}^{G(F_1)} \chi_1$
	(resp. $I_{B(F)}^{G(F)} \chi$). The Galois automorphism $\theta$ extends naturally to $G(F)$ and it preserves both $B(F)$ and the character $\chi$.  In particular it induces an automorphism of the induced representation $I(\chi)$ that we shall denote also by $\theta$ for simplicity.
	
	Let $\varphi$ be a test function on $G(F)$, and $\varphi_1$ be a base-change transfer of $\varphi$ to $G(F_1)$. (See Subsection \ref{subsection rhotrace}.)
	
	\begin{proposition}\label{ptw}
	We have the following relation between the twisted trace for $G(F)$ and the trace for $G(F_1)$:
		\begin{equation}\label{BCPS} 
			\Tr(I(\varphi_1,\chi_1)) = \Tr(I(\varphi,\chi)\circ \theta).
		\end{equation}
	\end{proposition}
	
	\begin{proof}  The proof is based on comparing the explicit formulas for the characters on the two sides of \eqref{BCPS}, as determined in \cite[Cor. 7.3.9]{HL}.  We thank Jean-Loup Waldspurger for the reference to \cite{HL} and for the following argument. 
		
		Since we are working with principal series, the argument of \cite{HL} can be much simplified in the following way. By an easy computation, the operator $I(\varphi,\chi)\circ \theta$ is acting on $I(\chi)$ through the kernel function
		\begin{equation*}
			\displaystyle K_{\varphi,\chi,\theta}(x,y)=\int_{T(F)\times U(F)} \varphi(x^{-1}tu \theta(y)) \chi(t) \delta_{B(F)}(t)^{1/2} du dt,\;\; x,y\in H,
		\end{equation*}
		where $\delta_{B(F)}$ stands for the modular character of $B(F)$. In other words, for every $e\in I(\chi)$, we have
		\begin{equation*}
			\displaystyle \left((I(\varphi,\chi)\circ \theta)e \right)(x)=\int_{B(F)\backslash G(F)} K_{\varphi,\chi,\theta}(x,y)e(y)dy.
		\end{equation*}
		It follows that the twisted character $\Tr(I(\varphi,\chi)\circ \theta)$ is given by the formula
		\[\begin{aligned}
			\displaystyle \Tr(I(\varphi,\chi)\circ \theta) & =\int_{B(F)\backslash G(F)} K_{\varphi,\chi,\theta}(x,x) dx \\
			& = \int_{B(F)\backslash G(F)} \int_{T(F)\times U(F)}\varphi(x^{-1} tu \theta(x)) \chi(t) \delta_{B(F)}(t)^{1/2} dudt dx.
		\end{aligned}\]
		The map $v\in U(F)\mapsto u=t^{-1}v^{-1}t\theta(v)\in U(F)$ is an isomorphism for $t$ in general position with Jacobian equal to
		\[\begin{aligned}
			\displaystyle & \left\lvert \det(\Ad(t^{-1})-\theta)\mid_{\Lie(U(F))}\right\rvert_F=\delta_{B(F)}(t)^{-1} \left\lvert \det(1-\Ad(t)\theta)\mid_{\Lie(U(F))}\right\rvert_F \\
			& =\delta_{B(F_1)}(t_1)^{-1} \left\lvert \det(1-\Ad(t_1))\mid_{\Lie(U(F_1))}\right\rvert_{F_1}=\delta_{B(F)}(t)^{-1/2} D^{G(F_1)}(t_1)^{1/2},
		\end{aligned}\]
		where $t_1=N_{F/F_1}(t)$ and $D^{G(F_1)}(t_1)=\left\lvert \det(1-\Ad(t_1))\mid_{\Lie(G)/\Lie(T)}\right\rvert_{F_1}$. By a change of variable, the above identity gives
		\begin{align}\label{eq1 ptw}
			\displaystyle \Tr(I(\varphi,\chi)\circ \theta) & =\int_{B(F)\backslash G(F)} \int_{T(F)\times U(F)}\varphi((vx)^{-1} t \theta(vx)) \chi(t) D^{H_1}(t_1)^{1/2} dudt dx \\
			\nonumber & =\int_{T(F)\backslash G(F)} \int_{T(F)} \varphi(x^{-1}t\theta(x)) \chi_1(t_1) D^{G(F_1)}(t_1)^{1/2} dt dx.
		\end{align}
	
	\begin{lemma}
	Let $\Ker(N)$ and $\mathrm{Im}(N)$ be the kernel and the image of the norm morphism $N_{F/F_1}: T(F)\to T(F_1)$ respectively. We have an action of $T(F)$ on $\Ker(N)$ by $t\cdot z=zt\theta(t)^{-1}$. Let $\Ker(N)/T(F)$ be the finite set of orbits for this action\footnote{This set is in natural bijection with $H^1(F/F_1,T)$.}. Let $da$ be the quotient of Haar measures $\frac{dt}{dt_1}$ on $T(F)$ and $T(F_1)$. Then, we have the following integration formula
	$$\displaystyle \int_{T(F)} f(t)dt=\int_{\mathrm{Im}(N)} \sum_{z\in \Ker(N)/T(F)} \int_{T(F)/T(F_1)} f(a^{-1}z\tilde{t}_1 a) da dt_1,\;\; f\in C_c(T(F)),$$
	where for $t_1\in \mathrm{Im}(N)$, we have denoted by $\tilde{t}_1\in T(F)$ any element such that $N_{F/F_1}(\tilde{t}_1)=t_1$.
	\end{lemma}

\begin{proof}
This is just saying that when we disintegrate the Haar measure $dt$ along $N_{F/F_1}: T(F)\to T(F_1)$ with respect to $dt_1$ we get on the fiber $\Ker(N)$ the unique Haar measure whose restriction to $T(F)/T(F_1)\subset \Ker(N)$ is $da$. Since we are working with Haar measures, it suffices to check this in a neighborhood of $1\in T(F)$. It then reduces to the following readily checked linear-algebraic statement: the two $F_1$-linear isomorphisms $\bigwedge^{top} \mathfrak{t}(F)\simeq \bigwedge^{top} \mathfrak{t}(F)/\mathfrak{t}(F_1)\otimes \bigwedge^{top} \mathfrak{t}(F_1)$ induced by the short exact sequences
$$\displaystyle 0\to \mathfrak{t}(F_1)\to \mathfrak{t}(F)\to \mathfrak{t}(F)/\mathfrak{t}(F_1)\to 0$$
and
$$\displaystyle 0 \to \mathfrak{t}(F)/\mathfrak{t}(F_1)\xrightarrow{1-\theta} \mathfrak{t}(F)\xrightarrow{\Tr_{F/F_1}} \mathfrak{t}(F_1)\to 0$$
coincide.
\end{proof}

By the above lemma, \eqref{eq1 ptw} can be rewritten as
\begin{equation*}
\displaystyle \Tr(I(\varphi,\chi)\circ \theta)=\int_{\mathrm{Im}(N)} \sum_{z\in \Ker(N)/T(F)} O_{z\tilde{t}_1,\theta}(\varphi) \chi_1(t_1) D^{G(F_1)}(t_1)^{1/2} dt_1,
\end{equation*}	
where, for $t\in T(F)$, $O_{t,\theta}(\varphi)=\int_{T(F_1)\backslash G(F)} \varphi(x^{-1}t\theta(x)) dx$ denotes the $\theta$-twisted orbital integral of $\varphi$ at $t$. For $t\in T(F)$ that is $\theta$-regular semisimple, the $\theta$-stable conjugacy class of $t$ is the union of the $\theta$-conjugacy classes of $zt_1$ for $z\in \Ker(N)/T(F)$. Therefore, we have
$$\displaystyle SO_{t,\theta}(\varphi)=\sum_{z\in \Ker(N)/T(F)} O_{zt,\theta}(\varphi)$$
and
\begin{equation*}
\Tr(I(\varphi,\chi)\circ \theta)=\int_{\mathrm{Im}(N)}  SO_{\tilde{t}_1,\theta}(\varphi) \chi_1(t_1) D^{G(F_1)}(t_1)^{1/2} dt_1
\end{equation*} 
On the other hand, the stable orbital integral $SO_{t_1}(\varphi_1)=O_{t_1}(\varphi_1)$ at $t_1\in T(F_1)$ reduces to the usual orbital integral. Therefore, by definition of the base-change transfer, we obtain the formula
		\[\begin{aligned}
			\displaystyle \Tr(I(\varphi,\chi)\circ \theta)=\int_{T(F_1)} O_{t_1}(\varphi_1) \chi_1(t_1) D^{G(F_1)}(t_1)^{1/2} dt_1
		\end{aligned}\]
Finally, the character of the induced $I(\chi_1)$ is given in a similar way by right hand side of the above formula and the proposition follows.
	\end{proof}

	Since the equation \eqref{BCPS} only computes the trace of functions $\varphi_1$ whose (regular semisimple) stable orbital integrals are supported on the norms from $G(F)$, we will also need the following lemma.
	
	Let $G_{sreg}(F_1)_{Nm}$ be the subset of strongly regular semisimple elements $\gamma\in G(F_1)$ that are norms of elements $\delta\in G(F)$. It is open and closed in the strongly regular semisimple locus $G_{sreg}(F_1)$. We denote similarly by $T(F_1)_{Nm}\subset T(F_1)$ the image of the norm $N_{F/F_1}: T(F)\to T(F_1)$.

	\begin{proposition}\label{principaldescent}
		Assume that $p$ does not divide the order of the (absolute) Weyl group of $G$. Then, the following hold:
		\begin{enumerate}[(i)]
			\item There exists a neighborhood $\mathcal{U}\subset G(F_1)$ of $1$ such that $\mathcal{U}\cap G_{sreg}(F_1)\subseteq G_{sreg}(F_1)_{Nm}$.
			
			\item Let $\chi: T(F_1)\to \bbC^\times$ be a smooth character and assume that for every character 
            \[ \eta: T(F_1)/T(F_1)_{Nm}\to \mathbb{C}^\times \]
            the (normalized) induction $I_{B(F_1)}^{G(F_1)}(\chi\eta)$ is irreducible. Let $(\pi_i)_{i\in I}$ be a family of smooth irreducible representations of $G(F_1)$  and $\epsilon\in \bbC^\times$ be a root of unity such that for every test function $\varphi\in C_c^\infty(G(F_1))$ with $\Supp(\varphi)\cap G_{sreg}(F_1)\subset G_{sreg}(F_1)_{Nm}$ we have
			$$\displaystyle \sum_{i\in I} \Tr \pi_i(\varphi)=\epsilon \Tr I(\chi,\varphi).$$
			Then, $I=\{i_0\}$ is a singleton and there exists a character 
            \[ \eta: T(F_1)/T(F_1)_{Nm}\to C^\times \]
            such that $\pi_{i_0}\simeq I_{B(F_1)}^{G(F_1)}(\chi\eta)$. In particular, $\epsilon=1$.
		\end{enumerate}
		
	\end{proposition}
	
	\begin{proof}
		$(i)$ Let $\mathcal{T}$ be a set of representatives for the $G(F_1)$-conjugacy classes of maximal tori of $G$ defined over $F_1$. Since $p$ does not divide the order of the Weyl group, $\mathcal{T}$ is finite. Indeed, since all maximal tori are conjugate over a separable closure $F_1^{sep}$ of $F_1$, there is a natural bijection
		$$\displaystyle \mathcal{T}\simeq \Ker(H^1(F_1,N_G(T))\to H^1(F_1,G)),$$
		where $N_G(T)$ denotes the normalizer of $T$ in $G$. On the other hand, the fiber of the map $H^1(F_1,N_G(T))\to H^1(F_1,W)$, where $W=N_G(T)/T$, above $\alpha\in H^1(F_1,W(G))$ can be identified with $H^1(F_1,T_\alpha)$, where $T_\alpha$ denotes the torus obtained by twisting the Galois action on $T(F_1^{sep})$ by $\alpha$. In particular, it is always finite. Moreover, since $p$ is prime to the order of $W$, $H^1(F_1,W)$ is finite, hence so is $\mathcal{T}$.
		
		Now, every regular semisimple element of $G(F_1)$ is conjugate to an element of $S(F_1)$ for some $S\in \mathcal{T}$ and the result follows from noting that the norm map $N_{F/F_1}: S(F)\to S(F_1)$ has open image.
		
		$(ii)$ Let $A\subset T$ be the maximal split subtorus and $\Delta$ be the set of simple roots of $A$ in $G$ with respect to $B$. Let
		$$H_T:T(F_1)\to \Hom(X^*(A),\bbQ)=\Hom(X^*(T),\bbQ)$$
		be the map defined by $\langle \chi^{alg},H_T(t)\rangle= \log \lvert \chi^{alg}(t)\rvert$ for $t\in T(F_1)$ and $\chi^{alg}\in X^*(T)$. Set
		$$T^+=\{t\in T(F_1)\mid \langle \alpha, H_T(t)\rangle<0,\forall\, \alpha\in \Delta \}$$
		and $G^+$ be the open subset of elements $g\in G(F_1)$ that are conjugate to an element of $T^+$. For $\varphi\in C_c^\infty(G(F_1))$, we define $\varphi^{(B)}\in C_c^\infty(T(F_1))$ by
		$$\displaystyle \varphi^{(B)}(t)=\delta_B(t)^{1/2}\int_{K\times N(F_1)} \varphi(k^{-1}tuk) dudk,\;\; t\in T(F_1),$$
		where $K$ is any hyperspecial subgroup of $G(F_1)$ (so that we have an Iwasawa decomposition $G(F_1)=B(F_1)K$). Then, by Casselman's asymptotic formula \cite[\S 4]{Ca95}, when $\Supp(\varphi)\subset G^+$ we have 
		$$\Tr \pi(\varphi)=\Tr \pi_U(\mathbf{1}_{T^+}\varphi^{(B)})$$
		for every irreducible representation $\pi$ where $\pi_U$ stands for the (normalized) Jacquet module of $\pi$ with respect to $B=TU$ and $\mathbf{1}_{T^+}$ for the characteristic function of $T^+$. The map $\varphi\mapsto \mathbf{1}_{T^+}\varphi^{(B)}$ induces a surjection $C_c^\infty(G(F_1))\to C_c^\infty(T^+)$. Therefore, the assumption entails that for every $\varphi_T\in C_c^\infty(T(F_1))$ with $\Supp(\varphi_T)\subset T^+\cap T(F_1)_{Nm}$ we have
		$$\displaystyle \sum_{i\in I} \Tr (\pi_i)_U(\varphi_T)=\epsilon \Tr I(\chi)_U(\varphi_T).$$
		All the Jacquet modules above are $T(F_1)$-representations of finite length and the previous identity translates to an equality in the Grothendieck group of finite length representations of $T(F_1)_{Nm}$ as:
		$$\displaystyle \sum_{i\in I} [(\pi_i)_U\mid_{T(F_1)_{Nm}}]=\epsilon [I(\chi)_U\mid_{T(F_1)_{Nm}}]=\sum_{w\in W(G_{F_1})} [(w\chi)\mid_{T(F_1)_{Nm}}],$$
		where $W(G_{F_1})=N_{G(F_1)}(T) / T(F_1)$ is the (relative) Weyl group. (The last equality is a simple consequence of the geometric lemma.) Inducing to $T(F_1)$, this implies in particular that for some $i_0\in I$, $(\pi_{i_0})_U$'s contains a character of the form $\chi\eta$, $\eta: T(F_1)/T(F_1)_{Nm}\to \mathbb{C}^\times$, in its semisimplification. Then, there exists a surjection $(\pi_{i_0})_U\twoheadrightarrow \chi\eta$ hence, by Frobenius reciprocity and since both $\pi_{i_0}$ and $I(\eta\chi)$ are irreducible (the latter by assumption), $\pi_{i_0}\simeq I(\eta\chi)$. Now, since the characters of $I(\eta\chi)$ and $I(\chi)$ coincide on $G_{rs}(F_1)_{Nm}$, we deduce that
		$$\displaystyle \sum_{i\in I\setminus\{i_0\}} \Tr \pi_i(\varphi)=(\epsilon-1) \Tr I(\chi,\varphi)$$
		for every $\varphi\in C_c^\infty(G(F_1))$ with $\Supp(\varphi)\cap G_{sreg}(F_1)\subset G_{sreg}(F_1)_{Nm}$. Applying this to $\varphi$ the characteristic function of a small enough compact-open subgroup (which is possible thanks to $(i)$), we obtain that $I\setminus \{i_0 \}$ must be empty and $\epsilon=1$.
	\end{proof}

	\begin{proof}[Proof of Proposition \ref{wild}]

		If $v \in S_0$, then $v$ remains prime in $L$; we will write abusively $v$ also for the extension of this place to any intermediate field. By construction, $\rho|_{\Gamma_{L^P_v}}$ is valued in a conjugate of $\widehat{T}$. We choose a character $\chi_0 : T(L^P_v) \to \overline{\bbQ}_\lambda^\times$ such that the Langlands dual character $\widehat{\chi}_0 : \Gamma_{L^P_v} \to \widehat{T}(\overline{\bbQ}_\lambda)$ is $\widehat{G}$-conjugate to $\rho|_{\Gamma_{L^P_v}}$.
		
		We will argue by induction.  Let $L^P \subset \bbK \subset \bbL \subset L$ with $\bbL/\bbK$ cyclic, $\rho_\bbL = \rho|_{\Gamma_\bbL}$, and $\rho_\bbK = \rho|_{\Gamma_\bbK}$. By Lemma \ref{irrparam}, $\rho_\bbL$ and $\rho_\bbK$ are conjugated to parameters $\widehat{\chi}:\Gamma_{\bbL_v}\to \widehat{T}(\overline{\bbQ}_\lambda)$ and $\widehat{\chi}_1:\Gamma_{\bbK_v}\to  \widehat{T}(\overline{\bbQ}_\lambda)$ respectively. We let $\chi: T(\bbL_v)\to \overline{\bbQ}_\lambda^\times$ and $\chi_1: T(\bbK_v)\to \overline{\bbQ}_\lambda^\times$ be the corresponding characters.
		
		Let $S$ denote the set of places of $\bbL$ lying above $S_0$, let $S_1$ denote the set of places of $\bbK$ lying above $S_0$, and let $U^S = \prod_{v \not\in S} G(\cO_{\bbL_v})$, and define $U^{S_1}$ similarly.
		
		We will show by induction on $[L : \bbL]$ that:  
		\begin{enumerate}
			\item The space $\CA_{0, \rho_\bbL}^{U^S}$ is an irreducible $\prod_{v \in S} G(\bbL_v) = G_S$-module. In particular, there is a unique cuspidal automorphic representation $\Pi_\bbL$ of $G(\bbA_\bbL)$ such that $r_\lambda(\Pi_\bbL) \sim \rho_\bbL$ and $\Pi_\bbL^{U^S} \neq 0$, which necessarily appears in the space of cusp forms with multiplicity 1.
			\item $\Pi_\bbL$ is locally generic everywhere.
			\item  For $v$ as above, its local component 
			$\Pi_{\bbL,v}$ is isomorphic to the principal series representation $I(\chi_v)$.  
		\end{enumerate}
		First, we claim that this holds for $L = \bbL$. Indeed, by Theorem \ref{multone} it suffices to know that for every cuspidal representation $\Pi$ of $G(\bbA_{\bbL})$ with global parameter $\rho_L$ and every $v\in S_0$, $\Pi_{v(+)}\simeq I(\mathbf{1})$, the unramified principal series induced from the trivial character; so that in particular $\Pi$ is unramified everywhere. Since $\rho\mid_{\Gamma_{L_{v(+)}}}$ is trivial by assumption, by the local-global compatibility of \cite{GLa}, $\Pi_{v(+)}$ has trivial Genestier-Lafforgue parameter. Hence, the result follows from \cite[Theorem 1.4(i)]{GHS}, since 
		\begin{itemize}
		\item by the assumption that $p$ does not divide the order of the Weyl group every supercuspidal representation of a Levi subgroup of $G(\bbL_{v(+)})$ is compactly induced \cite{Fin}, and
		\item  the representation $\rho\mid_{\Gamma_{L_{v(+)}}}$ is pure because it is the restriction of an irreducible parameter for $K_v$.
		\end{itemize}
		
		We now assume that the result holds over $\bbL$ and we want to prove that the analogous conditions hold over $\bbK$, with $I(\chi_v)$ replaced by $I(\chi_{v, 1})$. By Theorem \ref{hyp rho traces} and the base-change fundamental lemma \cite{K86}, we have for any matching test functions $f_S \in C_c^\infty( G_S)$, $f'_{S_1} \in C_c^\infty( G_{S_1})$, an identity 
		\[ \tr ( R(f_S) \times \theta \mid \CA_{0, \rho_\bbL}^{U^S} ) = \tr ( R(f'_{S_1}) \mid \CA_{0, \rho_\bbK}^{U^{S_1}} ). \]
		By the induction hypothesis, there exists a root of unity $\epsilon$ such that the left-hand side equals $\epsilon \tr ( f_S \times \theta \mid I(\chi_S) )$ where $\chi_S=\bigotimes_{v\in S} \chi_v$ and $\epsilon$ measures the difference between the global action of $\theta$ on $\CA_{0, \rho_\bbL}^{U^S}$ (defined as composition by the corresponding automorphism of $G(\bbL)\backslash G(\bbA_{\bbL})$) and its action on $I(\chi_S)$ (also coming from composition, realizing $I(\chi_S)$ as the space of sections of a line bundle on $B(\bbL_S)\backslash G(\bbL_S)$). By Proposition \ref{ptw}, this equals $\epsilon \tr( f'_{S_1} \mid I(\chi_1) )$. By Hypothesis \ref{hyp image bc}, we can choose $f'_{S_1}$ to be any function on $C_c^\infty(G_{S_1})$ with the property that $\operatorname{Supp}(f'_{S_1}) \cap G_{sreg}(\bbK_{S_1}) \subset G_{sreg}(\bbK_{S_1})_{Nm}$. We can therefore apply Lemma \ref{principaldescent}(ii) (verifying the irreducibility hypothesis as in the proof of Lemma \ref{irrprin}) to conclude that $\CA_{0, \rho_\bbK}^{U^{S_1}}$ is irreducible and in fact isomorphic to $\otimes_{v \in S_1} I(\chi_{1, v} \eta_v)$, for some characters $\eta_v : T(\bbK_v) / T(\bbK_v)_{Nm} \to \bbC^\times$. 
		
		To conclude the proof, it remains to check that we can in fact take $\eta_{v}  = 1$ for each $v \in S_1$. Since $I(\chi_{1, v} \eta_v)$ is irreducible, it suffices to show that the characters $\chi_{1, v}$ and $\chi_{1, v} \eta_v$ are conjugate under the Weyl group of $G$ (or equivalently, that their Langlands duals are conjugate under the Weyl group of $\widehat{G}$). Let $\Pi_\bbK$ denote the unique automorphic representation contributing to $\CA_{0, \rho_\bbK}^{U^{S_1}}$. Since $\rho_\bbK \sim r_\lambda(\Pi_\bbK)$, and by local-global compatibility, we find that the Langlands dual characters $\widehat{\chi}_{1, v} \widehat{\eta}_v : \Gamma_{\bbK_{v}} \to \widehat{T}(\overline{\bbQ}_\lambda)$ and $\widehat{\chi}_{1, v} : \Gamma_{\bbK_{v}} \to \widehat{T}(\overline{\bbQ}_\lambda)$ are $\widehat{G}(\overline{\bbQ}_\lambda)$-conjugate. Lemma \ref{lem_simultaneous_conjugacy} shows that they are in fact conjugate under the Weyl group of $\widehat{G}$; and this completes the proof. 
	\end{proof}
	\begin{lemma}\label{lem_simultaneous_conjugacy} 
		Let $n \geq 1$, and let $s_1, \dots, s_n, t_1, \dots, t_n \in \widehat{T}(\overline{\bbQ}_\lambda)$ be such that the tuples $(s_1, \dots, s_n)$ and $(t_1, \dots, t_n) \in \widehat{G}(\overline{\bbQ}_\lambda)^n$ are $\widehat{G}(\overline{\bbQ}_\lambda)$-conjugate. Then they are in fact conjugate under the Weyl group of $\widehat{G}$.
	\end{lemma}
	\begin{proof}
		We use induction on $n$. The case $n = 1$ is well-known. In general, we can first act by an element of the Weyl group to assume that in fact $s_i = t_i$ for each $i = 1, \dots, n-1$. Take $g \in \widehat{G}$ such that $\Ad(g)(s_i) = t_i$ for each $i = 1, \dots, n$. Let $C = \operatorname{Cent}_{\widehat{G}}(s_1, \dots, s_{n-1})$; then $g \in C$, $C^\circ$ is a reductive subgroup of $\widehat{G}$ with maximal torus $\widehat{T}$, and the elements $s_1, \dots, s_{n-1}$ lie in the centre of $C^\circ$. Moreover, every connected component of $C$ contains an element of the form $\dot{w} \in \operatorname{Norm}_{\widehat{G}}(\widehat{T})$, where the image $w$ of $\dot{w}$ in the Weyl group of $\widehat{G}$ fixes $s_1, \dots, s_{n-1}$ \cite[II.4.1]{Spr70}. Therefore we can find $\dot{w}$ of this form and $c \in C^\circ$ such that $g = c \dot{w}$. In particular, we have $\Ad(c)(w(s_n)) = t_n$, so we can find an element $\dot{w}' \in \operatorname{Norm}_{C^\circ}(\widehat{T})$ such that $\Ad(\dot{w}' \dot{w})(s_n) = t_n$. Then the element $w' w $ in the Weyl group of $\widehat{G}$ takes $(s_1, \dots, s_n)$ to $(t_1, \dots, t_n)$. 
	\end{proof}
    	\subsection{Base change}   Henceforward we admit the hypotheses on the twisted stable trace formula, as stated in \S  \ref{tracehypo}.  		
	Before we begin the proof of the descent we draw some conclusions from these hypotheses and from the construction of test functions in \S \ref{types}.
	
	Let $\bbL/\bbK$ be a cyclic subextension of $L/K$ of degree {\em prime to $p$}. Set $\rho_{\bbK}=\rho\mid_{\Gamma_{\bbK}}$, $\rho_{\bbL}=\rho\mid_{\Gamma_{\bbL}}$ and let $S_{0,\bbK}$ (resp. $S_{0,\bbL}$) be the finite set of places of $\bbK$ (resp. $\bbL$) that divide a place in $S_0$.
	
	For each place $v\in S_0$, we introduce sequences of $K$-types $(\cP^{M_v}_n,\chi_n^{M_v})_{n\geq n_0}$ for each subextension $K_v\subset M_v\subset L_v$ satisfying the conclusion of Theorem \ref{MWKV}, where $n_0\geq 1$ is an integer satisfying points $(iii)$ and $(iv)$ of that theorem. For each $n\geq n_0$ we set
	$$\displaystyle \mathcal{P}_{n,S_0}^{\mathbb{K}}=\prod_{v\in S_0} \cP_{n}^{\mathbb{K}_v}, \; \chi_{n,S_0}^{\mathbb{K}}=\prod_{v\in S_0} \chi_n^{\mathbb{K}_v},\; G(\mathbb{O}_{\mathbb{K}}^{S_0})=\prod_{v\notin S_0} G(\mathcal{O}_{\mathbb{K}_v})$$
	where $\mathcal{O}_{\mathbb{K}_v}$ denotes the ring of integers in $\mathbb{K}\otimes_K K_v=\prod_{w\mid v} \mathbb{K}_w$. Then, $\mathcal{P}_{n,S_0}^{\mathbb{K}}\times G(\mathbb{O}_{\mathbb{K}}^{S_0})$ is a compact-open subgroup of $G(\bbA_{\bbK})$ and $\chi_{n,S_0}^{\mathbb{K}}$ is a character on it. We denote by
	$$\mathcal{A}_{0,\rho_{\mathbb{K}}}(G)^{\mathcal{P}_{n,S_0}^{\mathbb{K}}\times G(\mathbb{O}_{\mathbb{K}}^{S_0}), \chi_{n,S_0}^{\mathbb{K}}}$$
	the subspace of $\mathcal{A}_{0,\rho_{\mathbb{K}}}(G)$ (see Subsection \ref{subsection rhotrace}) on which $\mathcal{P}_{n,S_0}^{\mathbb{K}}\times G(\mathbb{O}_{\mathbb{K}}^{S_0})$ acts by the character $\chi_{n,S_0}^{\mathbb{K}}$. Similar notation apply to $\bbL$ instead of $\bbK$.
	
	\begin{prop}\label{BCc}
		Assume that $\bbL/\bbK$ is of degree prime to $p$ and
		\begin{equation*}
			\displaystyle \dim \mathcal{A}_{0,\rho_{\mathbb{L}}}(G_{\bbL})^{\mathcal{P}_{n,S_0}^{\mathbb{L}}\times G(\mathbb{O}_{\mathbb{L}}^{S_0}), \chi_{n,S_0}^{\mathbb{L}}}=1.
		\end{equation*}
		Then,
		\begin{equation}\label{eq1BCc}
			\displaystyle \dim \mathcal{A}_{0,\rho_{\mathbb{K}}}(G_{\bbK})^{\mathcal{P}_{n,S_0}^{\mathbb{K}}\times G(\mathbb{O}_{\mathbb{K}}^{S_0}), \chi_{n,S_0}^{\mathbb{K}}}=1.
		\end{equation}
		In particular, $\mathcal{A}_{0,\rho_{\mathbb{K}}}(G)$ contains a unique cuspidal automorphic representation that is locally generic at every place in $S_0$.
	\end{prop}
	
	\begin{proof}   
		Let $f=\prod_v f_v\in C_c^\infty(G(\bbA_{\bbL}))$	and $f'=\prod_v f'_v\in C_c^\infty(G(\bbA_{\bbK}))$ be defined by
		\begin{itemize}
			\item For $v\notin S_0$, $f_v=\vol(G(\mathcal{O}_{\bbL_v}))^{-1}\mathbf{1}_{G(\mathcal{O}_{\bbL_v})}$, $f'_v=\vol(G(\mathcal{O}_{\bbK_v}))^{-1}\mathbf{1}_{G(\mathcal{O}_{\bbK_v})}$;
			\item For $v\in S_0$, $f_v=\vol(\cP^{\bbL_v}_n)^{-1} (\chi_n^{\bbL_v})^{-1}$ (where we extend $\chi_n^{\bbL_v}$ by zero to $G(\bbL_v)$), $f'_v=\vol(\cP^{\bbK_v}_n)^{-1} (\chi_n^{\bbK_v})^{-1}$.
		\end{itemize}
		Then, by the base-change fundamental lemma \cite{K86} as well as Theorem \ref{MWKV}, $f$ and $f'$ satisfy the assumption of Theorem \ref{hyp rho traces} i.e. for every $v$, $f'_v$ is a base-change transfer of $f_v$. Therefore, Theorem \ref{hyp rho traces} yields the identity
		\begin{equation*}
			\displaystyle \Tr(f\times \theta\mid \cA_{0,\rho_{\bbL}}(G_{\bbL}))=\Tr(f'\mid \cA_{0,\rho_{\bbK}}(G_{\bbK})).
		\end{equation*}
		Note that the right-hand side is $\dim \mathcal{A}_{0,\rho_{\mathbb{K}}}(G_{\bbK})^{\mathcal{P}_{n,S_0}^{\mathbb{K}}\times G(\mathbb{O}_{\mathbb{K}}^{S_0}), \chi_{n,S_0}^{\mathbb{K}}}$ whereas the left hand side is the trace of $\theta$ acting on $\mathcal{A}_{0,\rho_{\mathbb{L}}}(G_{\bbL})^{\mathcal{P}_{n,S_0}^{\mathbb{L}}\times G(\mathbb{O}_{\mathbb{L}}^{S_0}), \chi_{n,S_0}^{\mathbb{L}}}$. (Here we are using that the characters $\chi_n^\bbL$ are invariant under $\theta$ by construction, and thus the test functions $f_v$ are $\theta$-invariant as well.) As the latter is of dimension one by assumption, we see that the LHS is just a root of unity and \eqref{eq1BCc} follows.
		
		Since the (global) $K$-type $\mathcal{P}_{n,S_0}^{\mathbb{K}}\times G(\mathbb{O}_{\mathbb{K}}^{S_0}), \chi_{n,S_0}^{\mathbb{K}}$ appears with multiplicity one in $\mathcal{A}_{0,\rho_{\mathbb{K}}}(G_{\bbK})$, it appears in a unique cuspidal automorphic representation $\Pi_{\bbK}\subset \mathcal{A}_{0,\rho_{\mathbb{K}}}(G_{\bbK})$. Moreover, for every $v\in S_0$, the local representation $\Pi_{\bbK,v}$ is generic by Theorem \ref{MWKV} (iv).
		
	\end{proof}
    We now return to the situation of Proposition \ref{wild} -- thus we have a cuspidal automorphic representation $\Pi_{L^P}$ of $G(\mathbb{A}_{L^P})$ with global parameter $\rho_P$ and unramified at every place outside $S_0$. Moreover, $\Pi_{L^P}$ has automorphic multiplicity one.  Since $\rho_P$ extends to a Galois parameter over $K$, $\Pi_{L^P}$ is invariant under $\Gal(L^P/K)$, and the remainder of the proof consists in descending by cyclic subquotients of this soluble group, which is of order prime to $p$.  We deduce the following:
	\begin{prop}\label{tametrace}  Let $K \subset \bbK \subset \bbL \subset L^P$ be a cyclic extension of prime order.   
		Suppose that for every $v \in S_0$, there is an integer $n_v$ and a test function $f_{n_v,\bbL_v}$ as in Proposition \ref{MWKV} such that the trace of 
		$\otimes_{v \in S_0} f_{n_v,\bbL_v}$  on $\CA_{0, \rho_{\mathbb{L}}}(G(\bbL))_{\rho}$ equals $1$.
		
		Then for every $v \in S_0$, there is a test function $f_{n_v,\bbK_v}$ as in Proposition \ref{MWKV} such that the trace of 
		$\otimes_{v \in S_0} f_{n_v,\bbK_v}$  on $\CA_{0, \rho_{\mathbb{K}}}$ equals $1$.
		
	\end{prop}
	\begin{proof}  Since $\bbL/\bbK$ is of order prime to $p$, we may apply Proposition \ref{BCc} to obtain an identity of traces.  	\end{proof}

    \subsection{Deduction of the main theorem}
    We can now deduce our main theorem. For the convenience of the reader, we restate our assumptions from scratch. 
	\begin{thm}\label{existenceofpi} Let $F$ be a local field of positive characteristic $p$, and let $G$ be an unramified semisimple group over $F$. Suppose that:
    \begin{enumerate}
        \item The order of the absolute Weyl group of $G$ is prime to $p$.
        \item The hypotheses of Theorem \ref{hyp rho traces} are satisfied.
    \end{enumerate}
    Let $\tau : W_F \to {}^L G(\mathbb{C})$ be an irreducible $L$-homomorphism. Then there exists a generic supercuspidal representation $\pi$ of $F$ such that $\CL^{ss}(\pi) = \tau$.
	\end{thm}
	\begin{proof} Theorem \ref{potaut} guarantees the existence of data $L / K$, $G$, $S_0$, $\rho$, $\Pi_L$ etc. as at the beginning of \S \ref{sec_descent_arguments}. By considering successive subextensions $K \leq \mathbb{K} \leq \mathbb{L} \leq L$, we show that there exists a cuspidal automorphic representation $\Pi_{\mathbb{K}}$ of $G(\mathbb{A}_{\mathbb{K}})$, locally generic everywhere, and occurring with automorphic multiplicity 1, such that $r_\lambda(\Pi_{\mathbb{K}})$ and $\rho|_{\Gamma_{\mathbb{K}}}$ are equivalent. This holds for $\mathbb{L} = L$ by construction, and by Theorem \ref{multone}. When $\mathbb{L} = L^P$, it holds by Proposition \ref{wild}. This reduces us to considering the case where $\mathbb{L} / \mathbb{K}$ is cyclic of prime order, prime to $p$. In this case, we can make repeated applications of Proposition \ref{BCc} to conclude the existence of $\Pi_{\mathbb{K}}$ in the case $K = \mathbb{K}$. The result follows on taking $\pi = \Pi_{K, v}$ for any $v \in S_0$ and applying local-global compatibility (i.e. Theorem \ref{thm_loc_correspondence}).
	\end{proof}
	
	\begin{remark}  In particular, our method confirms Shahidi's Generic Packet Conjecture \cite[Conjecture 9.4]{Sh90} under the hypotheses of the Theorem; this also holds for $p$-adic groups under the hypotheses of Theorem \ref{p-adic} below.  Shahidi's conjecture is known unconditionally in many cases, including for $p$-adic classical groups and (under appropriate assumptions) for $\operatorname{GSpin}$ groups  \cite{Va17}. We plan to revisit the general case in a future paper.	
	\end{remark}

	\section{Applications to $p$-adic groups}\label{sec_p-adic_apps}

    So far in this paper, we have considered the local Langlands correspondence over equicharacteristic local fields. We now draw consequences for mixed characteristic ($p$-adic) local fields, using the recent results of Li-Huerta \cite{LH24}. Let us therefore assume that $F$ is such a field, and let $G$ be an unramified semisimple group over $F$. We will show:
	\begin{thm}\label{p-adic}  Assume the hypotheses on the trace formula over function fields as in \S \ref{tracehypo}, and that the prime $p$ does not divide the order of the absolute Weyl group of $G$. Then for any $n \geq 1$, there exists $d$ (depending only on $n$, the Dynkin type of $G$, and the cardinality $q$ of the residue field of $F$) with the following property: suppose that the absolute ramification index $e_F$ of $F$ is at least $d$, and let 
    \[ \tau : W_F \to {}^L G(\overline{\mathbb{Q}})\]
    be an irreducible parameter, 
    trivial on the inertia subgroup $I_F^n \leq \Gamma_F$ of level $n$. Then there exists a generic supercuspidal representation $\pi$ of $G(F)$ such that $\mathcal{L}_F^{ss}(\pi) = \tau$. 
	\end{thm}
    Here $\CL^{ss}_F$ denotes the semisimple correspondence constructed in \cite{FS} (contrary to the usage defined in \S \ref{sec_review_of_Galois}, which takes $F$ to be an equicharacteristic local field). 
    \begin{proof}

        Set $m=\max(n,n_0)$ where $n_0$ is the integer appearing in the statement of Theorem \ref{MWKV}.

        Let $F^\flat = \mathbb{F}_q((t))$. We can choose $d \geq n$ such that, if $e_F \geq d$, then the Hecke algebras $\mathcal{H}(G(F), K^{2m+1})$, $\mathcal{H}(G(F^\flat), K^{\flat, 2m+1})$ (where $K^?$ denotes the principal congruence subgroup of level $?$) are related by the Deligne--Kazhdan correspondence for local fields \cite{D84, Ka86, G21}. Increasing $d$ further, we can assume that the conclusion of \cite[Theorem B]{LH24} is valid with $n = 2m+1$. Since $d \geq n$, $\tau$ is identified with an irreducible parameter $\tau^\flat : W_{F^\flat} \to {}^L G (\overline{\mathbb{Q}})$, trivial on $I_{F^\flat}^n$. 

        By Lemma \ref{irrparam}, the restriction of $\tau^\flat$ to the wild inertia subgroup (up to conjugacy) through the maximal torus $\widehat{T}$. Thus, there exists a tamely ramified extension $K^\flat/F^\flat$ such that the restriction of $\tau^\flat$ to $W_{K^\flat}$ factors through $\widehat{T}$ and is therefore Langlands dual to a character $\chi : T(K^\flat) \to \overline{\mathbb{Q}}^\times$. Since $G$ is unramified, so does $T$. Thus, by the preservation of depth by the local Langlands correspondence for tamely ramified tori \cite[Theorem 7.10]{Yu09}, and the assumption on the ramification of $\tau^\flat$, $\chi$ is trivial on $T(K^\flat)_n$.
        
        By Theorem \ref{introthm_descent}, there exists an irreducible, admissible, generic supercuspidal representation $\pi^\flat$ of $G(F^\flat)$ such that $\mathcal{F}^{ss}_{F^\flat}(\pi^\flat) = \tau^\flat$. Moreover, Theorem \ref{MWKV}(iv) and Proposition \ref{BCc} together show that $(\pi^\flat)^{\mathcal{P}^{F^\flat}_m, \chi_m^{F^\flat}} \neq 0$. Proposition \ref{ncongruence} then shows that $(\pi^\flat)^{K^{\flat, 2m+1}} \neq 0$, and therefore that $\pi^\flat$ is related, by the Deligne--Kazhdan correspondence, to an irreducible, admissible representation $\pi$ of $G(F)$ such that $\pi$, $\pi^\flat$ give isomorphic modules for the Hecke algebras of the respective principal congruence subgroups of level $2m+1$; and in particular, such that $\pi^{\mathcal{P}_m, \chi_m} \neq 0$. Since $m \geq n_0$, Theorem \ref{MWKV}(iii) shows that $\pi$ is generic. Invoking \cite[Theorem B]{LH24} (allowable because of our choice of $d$), we see that $\mathcal{L}^{ss}_F(\pi) = \tau$. To conclude the proof, we need only check that $\pi$ is supercuspidal. However, this follows from the irreducibility of $\tau$, and the compatibility of $\mathcal{L}^{ss}_F$ with parabolic induction. 
    \end{proof}

    \subsection{Stability, following Chenji Fu}

    Let us continue to assume that $F$ is a $p$-adic field, and that $G$ is an unramified semisimple group over $F$. The expected properties of the conjectural local correspondence $\mathcal{L}_F$ (see \cite{Ka16} and the introduction to this paper) include not only surjectivity, but also a collection of endoscopic character identities to be satisfied by representations both of $G(F)$, and of the pure inner forms of $G$. In particular, for irreducible parameters with a non-empty packet, one expects that the associated $L$-packet should satisfy a stable character identity, and that the associated $L$-packets on the pure inner forms of $G$ should be likewise non-empty. 

    In this section, we recall results of Chenji Fu \cite{Fu24}, and explain consequences for stability. In an appendix to this paper, David Hansen gives a result about $L$-packets on inner forms. \begin{thm}\label{stable} Let $\tau : W_F \to {}^L G (\overline{\mathbb{Q}})$ be an irreducible parameter, and let $\pi$ be an irreducible admissible representation of $G(F)$ such that $\mathcal{L}_F^{ss}(\pi) = \tau$. Then there is a stable distribution $\Theta_{\pi_0}$ on $G(F)$, which is a finite $\Z$-linear combination of  characters of irreducible, supercuspidal representations of $G(F)$, and satisfying the following conditions:
    \begin{enumerate}
        \item Every representation $\pi'$ appearing with non-zero coefficient in $\Theta_{\pi_0}$ satisfies $\mathcal{L}_F^{ss}(\pi') = \tau$.
        \item The representation $\pi$ appears with non-zero coefficient in $\Theta_{\pi_0}$. 
    \end{enumerate}
    	 \end{thm}
     \begin{proof}
          This is an immediate consequence of \cite[Theorem 1.0.3]{Fu24}.
     \end{proof}
     As Hansen explained to us, Fu's result has the following implication:
	 \begin{corollary}\label{genericstable} Let $\tau, \pi$ be as in Theorem \ref{p-adic}, and suppose moreover that $\operatorname{Cent}_{\widehat{G}}(\tau(W_F)) = Z_{\widehat{G}}^{W_F}$. Then the character of $\pi$ is stable. 
	 \end{corollary}
	 Thus Shahidi's Generic Packet Conjecture is satisfied in a strong form when the $L$-packet is expected to be a singleton, provided the representation we construct is the only generic representation with parameter $\tau$.  

     At present Fu's argument is only available for $p$-adic groups.  This is because it relies on the Lefschetz trace formula computations in the article \cite{HKW} of Hansen, Kaletha, and Weinstein, which are only applicable to $p$-adic groups; and because Arthur's proof that stability on the semisimple locus implies stability does not apply in equal characteristic.
	
	 We briefly sketch the proof of Corollary \ref{genericstable}.  What Fu actually proves in \cite{Fu24} is that there is an action, denoted $\star$, of $\cO(S_\tau)$ on the derived category of smooth representations of $G(F)$ with the property that
	 $\cO(S_\tau)\star\pi$ has a stable character $\Theta_{\pi_0}$ on the regular semisimple elements of $G(F)$.  Here $S_\tau$ is the centralizer of $\tau$ in $\hG$, $\cO(S_\tau)$ is its coordinate ring, and the action on $\pi$ is deduced from the spectral action, constructed in \cite{FS}, of the category of quasi-coherent sheaves on the moduli stack $LS_{\lG}$ of $L$-parameters for $G(F)$ on the derived category of $\ell$-adic sheaves on the moduli stack $Bun_G$ of $G$-bundles on the Fargues-Fontaine curve.  This latter derived category naturally contains the derived category of smooth representations of $G(F)$, and $\cO(S_\tau)$ roughly appears as functions on a neighborhood of the point $\tau$ in $LS_\lG$.  (See \cite[\S 2.4]{Fu24} for a rigorous explanation of this action.)
	
	 Since $G$ is semisimple, $S_\tau$ is a finite group. Under the hypothesis of Corollary \ref{genericstable}, its characters are in bijection with inner forms of $G(F)$.  
     Hansen has explained to us an argument that the character $\chi_0 \in \cO(S_\tau)$ corresponding to the quasi-split inner form $G(F)$ has the property that $\chi_0 \star \pi = \pi$.  Thus $\Theta_{\pi_0}$ is just the character of $\pi$.  This immediately implies the Corollary.

    \appendix
    \section{Inner Forms, by David Hansen}

  In this short appendix, we show that the main existence results of
this paper quickly extend to inner forms.
\begin{prop}
Let $G/F$ be a reductive group over a nonarchimedean local field
$F$, and let $\varphi:W_{F}\to\phantom{}^{L}G$ be an irreducible
parameter. Suppose there exists an irreducible supercuspidal $G(F)$-representation
$\pi$ with Fargues-Scholze parameter $\varphi$. Then for any basic
element $b\in B(G)_{\mathrm{bas}}$ with associated extended pure
inner form $G_{b}$, there is an irreducible supercuspidal $G_{b}(F)$-representation
$\sigma$ with Fargues-Scholze parameter $\varphi$.
\end{prop}

\begin{proof}
We will use some standard Fargues-Scholze notation. Fix $\pi$ and
$b$ as in the proposition. Choose some irreducible algebraic representation
$V\in\mathrm{Irr}\,\hat{G}$ whose central character matches the Kottwitz
point $\kappa(b)$ of $b$ under the Kottwitz bijection $X^{\ast}(Z(\hat{G})^{\Gamma})=B(G)_{\mathrm{bas}}$,
and let $T_{V}$ denote the associated Hecke operator acting on sheaves
on $\mathrm{Bun}_{G}$. Consider the sheaf $A=T_{V}i_{1!}\pi$. We
first observe that $A$ is nonzero. For this, note that $T_{V}$ is
conservative, as the identity functor is canonically a retract of
$T_{V^{\vee}}T_{V}\cong T_{V^{\vee}\otimes V}$ using the monoidality
of the association $V\rightsquigarrow T_{V}$ {[}FS, Section IX.2{]}.

Next, note that under the bijection $\pi_{0}\mathrm{Bun}_{G}=B(G)_{\mathrm{bas}}$
induced by {[}FS, Theorem III.2.7 and Theorem III.4.5{]}, $A$ is
supported on the connected component $\mathrm{Bun}_{G}^{\kappa(b)}$
indexed by $\kappa(b)$. This follows from {[}Zou, Lemma 5.3.2{]}.
Now the connected component $\mathrm{Bun}_{G}^{\kappa(b)}$ has a
unique open stratum $j:[\ast/\underline{G_{b}(F)}]\to\mathrm{Bun}_{G}^{\kappa(b)}$
by {[}FS, Theorem III.4.5{]} again, and all non-open strata correspond
to inner forms of proper Levis in the quasisplit form of $G$. In
particular, if $i$ is the inclusion of such a stratum, then any irreducible
representation occurring as a subquotient of $i^{\ast}A$ will have
$L$-parameter factoring through the $L$-group of a proper Levi by
{[}FS, Theorem IX.7.2{]}. But since excursion operators commute with
Hecke operators, the $L$-parameter of such a subquotient automatically
agrees with (a twist of) $\varphi$, which does not factor through
any Levi in the $L$-group of $G$ by assumption. Therefore $i^{\ast}A=0$
for all non-open strata. Since $A\neq0$, this implies that $j^{\ast}A\neq0$.
As such, we can choose \emph{some }nonzero irreducible $G_{b}(F)$-representation
$\sigma$ occurring as a subquotient of $H^{\ast}(j^{\ast}A)$. Any
such $\sigma$ is automatically supercuspidal by {[}FS, Corollary
IX.7.3{]}, and its Fargues-Scholze parameter automatically agrees
with the parameter of $\pi$.
\end{proof}
If $G$ has connected center, then every inner form of $G$ is of
the form $G_{b}$ for some basic isocrystal. In general, given some
quasisplit $G$ with a specified inner form $H$, we can choose a
z-embedding $1\to G\to\tilde{G}\to D\to1$ where $D$ is a product
of tori induced from separable extensions and $\tilde{G}$ has connected
center. By the last condition, we can choose some basic element $b\in B(\tilde{G})_{\mathrm{bas}}$
such that $\tilde{H}\simeq\tilde{G}_{b}$, where $\tilde{H}$ is the
inner form of $\tilde{G}$ induced by $H$ together with the identification
$H^{1}(F,G_{\mathrm{ad}})=H^{1}(F,\tilde{G}_{\mathrm{ad}})$. Suppose
there exists an irreducible supercuspidal $G(F)$-representation $\pi$
with Fargues-Scholze parameter $\varphi$. Then $\pi$ can be extended
(by the discussion after {[}Kal, Fact 5.5{]}) to an irreducible $\tilde{G}(F)$-representation
$\tilde{\pi}$. By {[}FS, Theorem I.9.6.(v){]}, the Fargues-Scholze
parameter of $\tilde{\pi}$ is a lift of the parameter of $\pi$ along
the canonical map $\phantom{}^{L}\tilde{G}\to\phantom{}^{L}G$. By
the Proposition above, we also see that there is some irreducible
supercuspidal $\tilde{H}(F)=\tilde{G}_{b}(F)$-representation $\tilde{\sigma}$
with the same Fargues-Scholze parameter as $\tilde{\pi}$. Then taking
$\sigma$ to be any irreducible constituent of $\tilde{\sigma}|_{H(F)}$,
we see that the Fargues-Scholze parameter agrees of $\sigma$ with
the parameter of $\pi$ by a second application of {[}FS, Theorem
I.9.6.(v){]}.


\begin{thebibliography}{99}
		
		
		\bibitem[Adl98]{Adl98}
		J.~D. Adler, Refined anisotropic $K$-types and supercuspidal representations, Pacific J. Math. {\bf 185} (1998), no.~1, 1--32
		
		
		\bibitem[Art96]{Art96}
		James Arthur,
		\newblock  On local character relations
		\newblock {\em Selecta Math.}, 2(4): 501-579, 1996.
		
		\bibitem[Art02]{Artstab}
		James Arthur, 
		\newblock A stable trace formula. I. General expansions, 
		\newblock {\em J. Inst. Math. Jussieu}, 1(2): 175--277, 2002.
		
		\bibitem[Art13]{Art13}
		James Arthur,
		\newblock {\em The endoscopic classification of representations}, volume~61 of
		{\em American Mathematical Society Colloquium Publications}.
		\newblock American Mathematical Society, Providence, RI, 2013.
		\newblock Orthogonal and symplectic groups.
		
		
		
			
		
		\bibitem[BHKT19]{Boc19}
		Gebhard B\"{o}ckle, Michael Harris, Chandrashekhar Khare, and Jack~A. Thorne,
		\newblock {$\hat G$}-local systems on smooth projective curves are potentially
		automorphic.
		\newblock {\em Acta Math.}, 223(1):1--111, 2019.

        \bibitem[BS68]{BS68}
        A. Borel and T.~A. Springer, Rationality properties of linear algebraic groups. II, Tohoku Math. J. (2) {\bf 20} (1968), 443--497
		
		
		\bibitem[Bour]{Bour}
		N. Bourbaki, 
		\newblock Groupes et Alg\`ebres de Lie, Chapitres 4-6.
		
		\bibitem[Ca95]{Ca95}
		William Casselman,
		\newblock Introduction to the theory of admissible representations of $p$-adic reductive groups.
		\newblock Notes available at \url{https://personal.math.ubc.ca/\%7Ecass/research/pdf/p-adic-book.pdf}
		
		
		
				
		
		
		
		
		
		\bibitem[Con]{Con15}
		Brian Conrad,
		\newblock Lifting global representations with local properties.
		\newblock Preprint, available at\url{http://math.stanford.edu/{$\sim$}conrad/papers/}.
		
		
		
		
		\bibitem[D84]{D84}
		Pierre Deligne,
		\newblock Les corps locaux de caract\'eristique $p$, limites de corps locaux de caract\'eristique $0$.
		\newblock in Representations of reductive groups over a local field 119--157, 1984
		\newblock Paris:  Hermann, Travaux en Cours.
		
		\bibitem[Dri83]{Dri83}
		Vladimir Drinfel'd, 
		\newblock Two-dimensional $\ell$-adic representations of the fundamental group of a curve over a finite field and automorphic forms on $\GL(2)$.		
		\newblock \emph{Amer. J. Math.} 105 (1983), no. 1, 85--114.
		
		\bibitem[Dr]{Drin}
		Vladimir Drinfel'd, 
		\newblock Centralizers of semisimple elements of a simply connected semisimple group.
		\newblock Notes available at \url{https://math.uchicago.edu/~drinfeld/langlands/Centralizers.pdf}
		
		
		
		
		
		\bibitem[FS]{FS}  
		Laurent Fargues and Peter Scholze,
		\newblock Geometrization of the local Langlands correspondence
		\newblock arXiv:2102.13459.
		
		
		\bibitem[Fin]{Fin}
		J. Fintzen, 
		\newblock Types for tame $p$-adic groups, 
		\newblock Ann. of Math. (2) {\bf 193} (2021), no.~1, 303--346; MR4199732
		
		\bibitem[Fu24]{Fu24}
		Chenji Fu,
		\newblock Stability of elliptic Fargues--Scholze $L$-packets
		\newblock https://arxiv.org/abs/2501.00652v1
		
		
		\bibitem[GR24]{GR24}
		Dennis Gaitsgory and Sam Raskin,
		\newblock Proof of the geometric Langlands conjecture V:  The multiplicity one theorem.
		\newblock https://arxiv.org/abs/2409.09856
		
		\bibitem[GR25]{GR25}
		Dennis Gaitsgory and Sam Raskin,
		\newblock Geometric Langlands in positive characteristic from characteristic $0$.
		\newblock https:/arXiv:2508.02237
		
		\bibitem[G21]{G21}
		Radhika Ganapathy,
		\newblock A Hecke algebra isomorphism over close local fields.
		\newblock   Available at \texttt{https://arxiv.org/pdf/2103.12363}.
		
		
		\bibitem[GLa]{GLa}
		Alain Genestier and Vincent Lafforgue,
		\newblock  Chtoucas restreints pour les groupes r\'eductifs et param\'etrisation de Langlands locale.
		\newblock Preprint. Available at \texttt{https://arxiv.org/abs/1709.00978}.
		
		\bibitem[GHS]{GHS}  Wee Teck Gan, Michael Harris, and Will Sawin,
		\newblock Local parameters of supercuspidal representations.
		\newblock {\em Forum of Math., Pi}, 12, e13, 2024.
		
		\bibitem[Hans]{Hansen}
		David Hansen, Beijing notes on the categorical local Langlands conjecture, avilable at \url{https://davidrenshawhansen.net/Beijing.pdf}

        \bibitem[HKW]{HKW}
        David Hansen, Tasho Kaletha, and Jared Weinstein
        \newblock On the Kottwitz conjecture for local shtuka spaces,
        \newblock {\em Forum of Math, Pi}, 10, e13, 2022.
		
		\bibitem[HL]{HL}  Guy Henniart and Bertrand Lemaire,
		\newblock Repr\'esentations des espaces tordus sur un groupe r\'eductif connexe $\mathfrak{p}$-adique.
		\newblock {\em Ast\'erisque},  386,  2017.

        \bibitem[Ka16]{Ka16}
		Tasho Kaletha,
		\newblock The local {L}anglands conjectures for non-quasi-split groups
		\newblock  {\em Families of automorphic forms and the trace formula}, 217--257, Simons Symp., Springer, 2016.
        
		\bibitem[Ka19a]{Ka19a}
		Tasho Kaletha,
		\newblock Regular supercuspidal representations
		\newblock  {\em J. Am. Math. Soc.} 32(4) 1071--1170, 2019.
		
		
		\bibitem[Ka19b]{Ka19b}
		Tasho Kaletha,
		\newblock Supercuspidal $L$-packets
		\newblock  arXiv:1912.03274v2 [math.RT].
		
		\bibitem[KP]{KP}
		Tasho Kaletha, Gopal Prasad, 
		\newblock  Bruhat-Tits theory a new approach.
        \newblock  Cambridge University Press, 2023.
        
		
		\bibitem[Ka86]{Ka86}
		David Kazhdan,
		\newblock Representations of groups over close local fields.
		\newblock {\em J. d'Analyse Math.}   47, 175--179, 1986. 
		
		
		
		
		
		
		
		\bibitem[Ko02]{Ko}  T. Konno,
		\newblock Twisted endoscopy and the generic packet conjecture. 
		\newblock {\em Israel J. Math.}, {\bf 129} 253--289, 2002.

        \bibitem[K82]{K82}
		R.~E. Kottwitz, Rational conjugacy classes in reductive groups, Duke Math. J. {\bf 49} (1982), no.~4, 785--806
		
		\bibitem[K84]{K84}
		R.~E. Kottwitz, Stable trace formula: cuspidal tempered terms, Duke Math. J. {\bf 51} (1984), no.~3, 611--650; MR0757954
		
		\bibitem[K86]{K86}  R. Kottwitz, 
		\newblock Base change for unit elements of Hecke algebras
		\newblock {\em Compositio Math.}, {\bf 60} 237--250 (1986).
		
		\bibitem[KS]{KS}
		R.~E. Kottwitz and D. Shelstad, Foundations of twisted endoscopy, Ast\'erisque No. 255 (1999), vi+190 pp.
		
		
		\bibitem[LLap15]{LLap15}
		J.-P. Labesse, E. Lapid,
		Characters of $G$ over local and global fields, appendix to
		E. Lapid, Z. Mao, A conjecture on Whittaker--Fourier coefficients of cusp forms, J. Number Theory {\bf 146} (2015), 448--505
		
		\bibitem[LL21]{LL21}
		Jean-Pierre Labesse and Bertrand Lemaire,
		\newblock La formule des traces tordue pour les corps de fonctions
		\newblock \url{arXiv:2102.02517v2 [math.RT]}
		
		\bibitem[Laf18]{Lafa}
		Vincent Lafforgue,
		\newblock Chtoucas pour les groupes r{\'e}ductifs et param{\'e}trisation de
		{L}anglands globale.
		\newblock {\em J. Amer. Math. Soc.} 31(3), 719--891,  2018. 
		
		

        \bibitem[LRS93]{LRS93}
		G. Laumon, M. Rapoport, and U. Stuhler,
		\newblock ${\mathcal{D}}$-elliptic sheaves and the {L}anglands
              correspondence
		\newblock {\em Invent. Math.}, 113(2):217--338, 1993.
        
		
		
		\bibitem[LH23]{LH23}  
		Siyan Daniel Li-Huerta,
		\newblock Local-global compatibility over function fields
		\newblock  	arXiv:2301.09711 [math.NT]
		
		\bibitem[LH24]{LH24}  
		Siyan Daniel Li-Huerta,
		\newblock On close fields and the local Langlands correspondence
		\newblock  	arXiv:2407.07063 [math.NT]
			
		\bibitem[MW87]{MW87}
		C. M\oe glin and J.-L. Waldspurger,
		\newblock Mod\`eles de Whittaker d\'eg\'en\'er\'es pour des groupes $p$-adiques
		\newblock Math. Z. {\bf 196} (1987), no.~3, 427--452
		
		
		\bibitem[MW16a]{MWstab1}
		C. M\oe glin and J.-L. Waldspurger, {\it Stabilisation de la formule des traces tordue. Vol. 1}, Progress in Mathematics, 316, Birkh\"auser/Springer, Cham, 2016
		
		\bibitem[MW16b]{MWstab2}
		C. M\oe glin and J.-L. Waldspurger, {\it Stabilisation de la formule des traces tordue. Vol. 2}, Progress in Mathematics, 317, Birkh\"auser/Springer, Cham, 2016
		
		\bibitem[MB89]{Mor89}
		Laurent Moret-Bailly,
		\newblock Groupes de {P}icard et probl\`emes de {S}kolem. {II}.
		\newblock {\em Ann. Sci. \'{E}cole Norm. Sup. (4)}, 22(2):181--194, 1989.
		
		\bibitem[MB90]{Mor90}
		Laurent Moret-Bailly,
		\newblock Extensions de corps globaux {\`a} ramification et groupe de Galois donn{\'e}s. 
		\newblock {\em C. R. Acad. Sci. Paris S{\'e}r. I Math. 311} (1990), no. 6, 273--276.
		
		
				
		\bibitem[Ngo]{Ngo}
		B.~C. Ng\^o, Le lemme fondamental pour les alg\`ebres de Lie, Publ. Math. Inst. Hautes \'Etudes Sci. No. 111 (2010), 1--169
		
		
		\bibitem[Ras]{Ras} 
		S. Raskin,
		\newblock{Conjectures of Arthur and Ramanujan for unramified automorphic forms:  announcement}, 
		\newblock to appear in the Proceedings of the 2026 International Congress of Mathematicians.

		
		\bibitem[Rich]{Rich}
		S. Riche,
		\newblock {Kostant section, universal centralizer, and a modular derived Satake equivalence}, 
		\newblock {\em Math. Z.}, {\bf 286} (2017), no.~1-2, 223--261
		
			
		\bibitem[Rod74]{Rod74}
		F. Rodier,
		\newblock{Mod\`ele de Whittaker et caract\`eres de repr\'esentations}
		\newblock in {\it Non-commutative harmonic analysis (Actes Colloq., Marseille-Luminy, 1974)}, pp. 151--171, Lecture Notes in Math., Vol. 466, Springer, Berlin-New York,
		
		\bibitem[ST21]{ST21}
		Will Sawin and Nicolas Templier,
		\newblock {\em On the {R}amanujan conjecture for automorphic forms over function fields {I}. {G}eometry}.
		\newblock To appear in Journal of the American Mathematical Society.
		
		\bibitem[Sc25]{Sc25} 
		P. Scholze,
		\newblock {\em  Geometrization of the local Langlands correspondence, motivically},
		\newblock  	arXiv:2501.07944 [math.AG]
		
		
		\bibitem[Sh90]{Sh90}  
		F. Shahidi,
		\newblock {\em A proof of Langlands' conjecture on Plancherel measures; complementary
			series for $p$-adic groups}
		\newblock  Ann. of Math. \. 132, (1990) no. \ 2,  273--330.
			
		
		\bibitem[Sol25]{Sol25}
		M. Solleveld, On principal series representations of quasi-split reductive $p$-adic groups, Pacific J. Math. {\bf 334} (2025), no.~2, 269--327.
		
		\bibitem[SS70]{Spr70}
		T.~A. Springer and R.~Steinberg,
		\newblock Conjugacy classes.
		\newblock In {\em Seminar on {A}lgebraic {G}roups and {R}elated {F}inite
			{G}roups ({T}he {I}nstitute for {A}dvanced {S}tudy, {P}rinceton, {N}.{J}.,
			1968/69)}, Lecture Notes in Mathematics, Vol. 131, pages 167--266. Springer,
		Berlin, 1970.

        \bibitem[Ste65]{Ste65}
        R. Steinberg, Regular elements of semisimple algebraic groups, Inst. Hautes \'Etudes Sci. Publ. Math. No. 25 (1965), 49--80;
			
		
		\bibitem[Va17]{Va17}
		Sandeep Varma,
		\newblock On descent and the generic packet conjecture. 
		\newblock {\em Forum Math.} 29, no. 1, 111--155, 2017.
			
		\bibitem[Wal97]{Waltransfer}
		J.-L. Waldspurger, Le lemme fondamental implique le transfert, Compositio Math. {\bf 105} (1997), no.~2, 153--236
		
		
		
		\bibitem[Yu09]{Yu09} J.-K. Yu,
        \newblock On the local Langlands correspondence for tori,
        \newblock in {\it Ottawa lectures on admissible representations of reductive $p$-adic groups}, 177--183, Fields Inst. Monogr., 26, Amer. Math. Soc., Providence, RI,

         \bibitem[Zou]{Zou}  K. Zou,
        \newblock The categorical form of Fargues' conjecture for tori,
         \newblock \url{arXiv:2202.13238v2 [math.RT] }
 for tori
	\end{thebibliography}
\end{document}